\documentclass[a4paper,10pt]{amsart}
\usepackage{amsmath,amsfonts,amssymb, amsthm, xypic}
\usepackage[all]{xy}

\def\g{\mathfrak{g}}
\def\n{\mathfrak{n}}
\def\C{\mathbb{C}}

\def\U+{\mathbf{U}_v^+(\mathcal{L}\mathfrak{g})}
\def\Z{\mathbb{Z}}
\def\Q{\mathbb{Q}}
\def\Ox{\mathcal{O}_X}
\def\O{\mathcal{O}}
\def\Hom{\mathrm{Hom}}
\def\Ext{\mathrm{Ext}}
\def\qed{$\hfill \checkmark$}
\def\N{\mathbb{N}}
\def\tto{\twoheadrightarrow}
\def\lto{\longrightarrow}
\def\a{\alpha}
\def\b{\beta}
\def\ab{\alpha+\beta}
\def\qlb{\overline{\mathbb{Q}_l}}

\newtheorem{theo}{\bf{Theorem}}
\newtheorem{lem}{Lemma}
\newtheorem{cor}{Corollary}
\newtheorem{prop}{Proposition}
\newtheorem{claim}{Claim}
\newtheorem{conj}{Conjecture}

\numberwithin{lem}{section}
\numberwithin{conj}{section}
\numberwithin{prop}{section}
\numberwithin{cor}{section}
\numberwithin{theo}{section}
\numberwithin{equation}{section}
\numberwithin{claim}{section}

\title{Canonical bases of loop
algebras via Quot schemes, I}
\author{Olivier Schiffmann}
\date{}

\begin{document}
\maketitle

\tableofcontents

\section{Introduction}

\vspace{.1in}

\paragraph{}One of the main outcomes of the theory of quantum groups has been the discovery, by Kashiwara
and Lusztig, of certain distinguished bases in quantized enveloping algebras
with very favorable properties~: positivity of structure constants, compatibility with all highest weight
integrable representations, etc. These canonical bases, which are in some ways analogous to the Kazhdan-Lusztig
bases of Hecke algebras, have proven to be a powerful tool in the study of the representation theory of (quantum)
Kac-Moody algebras, especially in relation to solvable lattice and vertex models in statistical mechanics
(see e.g. \cite{JMMO}, \cite{KMS}). They also sometimes carry representation-theoretic information themselves,
encoding character formulas and decomposition numbers (\cite{A}, \cite{LLT}, \cite{VV}).

The construction of the canonical basis $\mathbf{B}$ of a quantized enveloping algebra
$\mathbf{U}^+_v(\mathfrak{g})$ in \cite{L1} is based on Ringel's earlier discovery that $\mathbf{U}^+_v(\g)$
may be realized in the Hall algebra of the category of representations of a quiver $Q$ whose underlying graph
is the Dynkin diagram of $\g$ (\cite{Ri}). The set of isomorphism classes of representations of $Q$ of a given class 
$\mathbf{d} \in K_0(\mathrm{Rep}\;Q)$ is the set of orbits of a reductive group $G_{\mathbf{d}}$ on a vector space
$E_{\mathbf{d}}$; Lusztig realizes $\mathbf{U}^+_v(\g)$ geometrically as a convolution algebra of 
($G_{\mathbf{d}}$-equivariant, semisimple)
constructible sheaves on $E_{\mathbf{d}}$ (for all $\mathbf{d}$), and obtains the canonical
basis as the set of all simple perverse sheaves appearing in this algebra.

\vspace{.1in}

\paragraph{}The present work is an attempt to carry out the above constructions in the situation when $\g$ is no 
longer a Kac-Moody algebra but a loop algebra of a Kac-Moody algebra.  
An affinization of Ringel's construction was already
considered in \cite{S1}~: the Hall algebra of the category of coherent sheaves
on a weighted projective line (in the sense of \cite{GL}) provides a realization of a positive part of the (quantum)
affinization $\mathbf{U}_v(\mathcal{L}\mathfrak{g})$ of a Kac-Moody algebra $\g$ associated to a star Dynkin diagram.
In the case of a projective line with no weight, which was first considered by Kapranov \cite{Kap}, this gives the positive 
half of $\mathbf{U}_v(\widehat{\mathfrak{sl}}_2)$ in \textit{Drinfeld's} presentation. The most interesting situation occurs
when $\mathfrak{g}$ is itself an affine Lie algebra, in which case $\mathbf{U}_v(\mathcal{L}\g)$ is no longer a Kac-Moody
algebra but an elliptic (also called double-loop, or toroidal) algebra (see \cite{MRY}, \cite{Sai}).

\vspace{.1in}

In this paper, we propose a construction of a canonical basis $\mathbf{B}'$
(for the positive part $\mathbf{U}_v(\mathcal{L}\n)$) of
$\mathbf{U}_v(\mathcal{L}\g)$. Rather than categories of coherent sheaves on a weighted projective line we use the 
(equivalent) categories
of coherent sheaves on a smooth projective curve $X$, equivariant with respect to a fixed group $G \subset Aut(X)$
satisfying $X/G \simeq \mathbb{P}^1$. We consider a convolution algebra $\mathfrak{U}_{\mathbb{A}}$
of (equivariant, semisimple) constructible sheaves
on various Quot schemes, and take as a basis of $\mathfrak{U}_{\mathbb{A}}$
the set of all simple perverse sheaves appearing in this way. Note that
we need the (smooth part) of the \textit{whole} Quot scheme and not only its (semi)stable piece.
One difference between the coherent sheaves situation and the quiver situation is that in the former all regular 
indecomposables (torsion sheaves) lie ``at the top'' of the category 
(i.e they admit no nonzero homomorphism to a nonregular indecomposable). For instance, when
$X=\mathbb{P}^1$, this makes the picture essentially 
``of finite type'', although the corresponding quantum group is a quantum affine algebra (in particular, it is possible to
explicitly describe all simple perverse sheaves occuring in the canonical basis). On the other hand, unlike the
quiver case, the lattice of submodules has infinite depth in general, and it is necessary to consider some
inductive limit of Quot schemes (or the category of perverse sheaves on such an inductive limit). As a consequence, 
our canonical basis $\mathbf{B}'$ naturally lives in a completion $\widehat{\mathfrak{U}}_{\mathbb{A}}$ of
$\mathfrak{U}_{\mathbb{A}}$.

\vspace{.1in}

We conjecture that 
algebra $\widehat{\mathfrak{U}}_{\mathbb{A}}$ is isomorphic to a completion $\widehat{\mathbf{U}}_v(\mathcal{L}\n)$ of
$\mathbf{U}_v(\mathcal{L}\n)$ and prove this conjecture in the finite type case
(when $X=\mathbb{P}^1$, i.e. when $\mathbf{U}_v(\mathcal{L}\g)$ is a quantum affine algebra). For instance, we get
a canonical basis of (a completion of)
the positive half of $\mathbf{U}_v(\widehat{\mathfrak{sl}}_2)$ in Drinfeld's presentation.

\vspace{.1in}

Of course, when $\mathfrak{g}$ is of finite type, the loop algebra $\mathcal{L}\g$ is again a Kac-Moody algebra, and there
is already a canonical basis $\mathbf{B}$ for (the Kac-Moody positive part of) 
the quantum group $\mathbf{U}_v(\mathcal{L}\g)$ by
\cite{K} and \cite{L1}. Since $\mathbf{B}$ and $\mathbf{B}'$ are bases of two different Borel subalgebras, it is not 
possible to compare them directly. However, we conjecture that the basis $\mathbf{B}'$ is also compatible with
highest weight integrable representations, and we show in some examples (for $\mathfrak{\g}=\mathfrak{sl}_2$ and
the fundamental representation $\Lambda_0$) that $\mathbf{B}$ and $\mathbf{B}'$ 
project to the same bases in such a module. The action of $\mathbf{B}'$ on $\Lambda_0$ is particularly well-suited
for the ``semi-infinite'' description of $\Lambda_0$ given by Feigin and Stoyanovski, \cite{BS}~: the natural
filtration of $\mathbf{B}'$ (constructed from the Harder-Narasimhan filtration)
induces the filtration of $\Lambda_0$ by ``principal subspaces'' considered in \cite{BS}.
The coincidence of $\mathbf{B}$ and $\mathbf{B}'$ in $\Lambda_0$ suggests that $\mathbf{B}$ and $\mathbf{B}'$ are
incarnations of the same ``canonical basis'', which should be essentially unique.

\vspace{.1in}

More interestingly, we will show in \cite{S2} that $\widehat{\mathfrak{U}}_{\mathbb{A}}$ and 
$\widehat{\mathbf{U}}_v(\mathcal{L}\n)$ are isomorphic when $X$ is of genus one and $G=\Z/2\Z, \Z/3\Z, \Z/4\Z$ or
$\Z/6\Z$. This implies the existence
of a canonical basis in quantum toroidal algebras of type $D_4, E_6,E_7,E_8$, and provides a first example
of a canonical basis in a non Kac-Moody setting. We will also construct a canonical basis for 
$\mathbf{H}_\mathcal{E}=\widehat{\mathfrak{U}}_{\mathbb{A}}$ when $X=\mathcal{E}$ is an elliptic curve
and $G=Id$. The algebra $\mathbf{H}_{\mathcal{E}}$ is an elliptic analog of a Heisenberg algebra and its canonical basis
should have an interesting combinatorial description.

The methods of proof of both the present work and \cite{S2} are adapted from \cite{S1} and rely on the HN filtration
and on the theory of mutations. These methods do not seem to extend to the general case where $X$ is of genus at least
two. In \cite{L1}, wild quivers are dealt with by exploiting the freedom of choice of an orientation,
via a Fourier-Deligne transform. We could not find a proper substitute in the category of coherent
sheaves on a curve.
    
\vspace{.1in}

The plan of the paper is as follows. In Section~1 we recall the definitions of loop Kac-Moody algebras, of their
quantum cousins and of the cyclic Hall algebras. In Section~2 we consider the category $Coh_G(X)$ of $G$-equivariant
coherent sheaves on a curve $X$ and construct the relevant Quot schemes. Section~3 introduces inductive limits
of Quot schemes and the category of perverse sheaves on them. In Sections~4,5 we define the convolution algebra
$\widehat{\mathfrak{U}}_{\mathbb{A}}$ following the method in \cite{L1}. Our main theorem and conjecture
are in Section~5.4. Section~6 and Section~7 collect some technical results pertaining to the Harder-Narasimhan
stratification of Quot schemes. Our main theorem is proved in Sections~8 and 9. Section~10 relates
the algebras $\widehat{\mathfrak{U}}_{\mathbb{A}}$ and $\widehat{\mathbf{U}}_v(\mathcal{L}\n)$, and
gives some first properties of the canonical basis $\mathbf{B}'$. Section~11
provides some formulas for $\mathbf{B}'$ when $\mathcal{L}\g= \widehat{\mathfrak{sl}}_2$ as well as some
formulas for its action on the fundamental representation $\Lambda_0$. 

\vspace{.2in}

\section{Loop algebras and quantum groups.}

\vspace{.1in}

\paragraph{\textbf{1.1. Loop algebras.}} Let $\Gamma$ be a star shaped Dynkin
diagram with branches $J_1, \ldots, J_N$ of respective lengths $p_1-1, \ldots, p_N-1$.
Let $\star$ denote the central node and for $i=1, \ldots, N$ and $j=1, \ldots, p_i-1$, let
$(i,j)$ be the $j$th node of $J_i$, so that $\star$ and $(i,1)$ are adjacent for all $i$.
To $\Gamma$ is associated a generalized Cartan matrix
$A=(a_{\gamma\sigma})_{\gamma,\sigma \in \Gamma}$ and a Kac-Moody algebra
$\mathfrak{g}$. We will consider an affinization $\mathcal{L}\mathfrak{g}$ of $\g$
generated by some elements
$h_{\gamma,k}, e_{\gamma,k}, f_{\gamma,k}$ for $\gamma\in \Gamma$, $k \in \Z$ and
$\textbf{c}$ subject to the following set of relations :
\begin{equation*}
[h_{\gamma,k}, h_{\sigma,l}]=k \delta_{k,-l}a_{\gamma\sigma}\textbf{c},
\qquad [e_{\gamma,k},f_{\sigma,l}]=\delta_{\gamma,\sigma}h_{\gamma,k+l}+ k \delta_{k,-l} \textbf{c},
\end{equation*}
\begin{equation*}
[h_{\gamma,k}, e_{\sigma,l}]=a_{\gamma\sigma}e_{\sigma,l+k},\qquad [h_{\gamma,k}, f_{\sigma,l}]=
-a_{\gamma\sigma}f_{\sigma,k+l},
\end{equation*}
\begin{equation*}
[e_{\gamma,k+1},e_{\sigma,l}]=[e_{\gamma,k},e_{\sigma,l+1}],\qquad
[f_{\gamma,k+1},f_{\sigma,l}]=[f_{\gamma,k},f_{\sigma,l+1}]
\end{equation*}
\begin{equation*}
\begin{split}
&[e_{\gamma,k_1},[e_{\gamma,k_2},[\ldots [e_{\gamma,k_n},e_{\sigma,l}]\cdots ]=0 \qquad
\mathrm{if\;} n=1-a_{\gamma\sigma},\\
&[f_{\gamma,k_1},[f_{\gamma,k_2},[\ldots [f_{\gamma,k_n},f_{\sigma,l}]\cdots ]=0 \qquad
\mathrm{if\;} n=1-a_{\gamma\sigma}.
\end{split}
\end{equation*}
The Lie algebra $\g$ is embedded in $\mathcal{L}\g$ as the subalgebra generated by
$e_{\gamma,0}, f_{\gamma,0}, h_{\gamma,0}$ for $\gamma \in \Gamma$.
When $\g$
is a simple Lie algebra (resp. an affine
Kac-Moody algebra), $\mathcal{L}\g$ is the corresponding affine Lie algebra $\widehat{\g}$ (resp. the
corresponding toroidal algebra with infinite dimensional center, see \cite{S1}).

We write $\Delta$ (resp. $Q$)  for the root system
(resp. the root
lattice) of $\g$, so that the root system of $\mathcal{L}\g$ is
$\widehat{\Delta}=(\Delta + \Z \delta)\cup \Z^* \delta$ and the weight lattice is $\widehat{Q}=Q \oplus \Z \delta$. 
Here $\delta$ is the extra imaginary root.
The lattices $Q$ and $\widehat{Q}$ are both equipped with standard symmetric bilinear forms
$(\,,\,)$ with values in $\Z$.
Write $\alpha_\gamma$ for the simple root corresponding to a vertex $\gamma \in  \Gamma$.

\vspace{.05in}

Let $\mathfrak{g}_1 \simeq \mathfrak{sl}_{p_1}, \ldots, \g_N \simeq \mathfrak{sl}_{p_N}$ be the
subalgebras of $\g$ associated to the branches $J_1, \ldots, J_N$. For $i=1, \ldots, N$ the map
$e_{\gamma,l} \mapsto e_\gamma t^l, f_{\gamma,l} \mapsto f_\gamma t^l, h_{\gamma,l} \mapsto h_\gamma t^l$
gives rise to an isomorphism
between the subalgebra $\mathcal{L}\g_i$ of $\mathcal{L}\g$ generated by
 $\{e_{\gamma,l}, f_{\gamma,l}, h_{\gamma,l}\:|\;\gamma \in J_i , l \in \Z\}$ and $\widehat{\mathfrak{sl}}_{p_i}$.
We denote by $\widehat{\n}_i \subset \widehat{\mathfrak{sl}}_{p_i}$ the standard nilpotent subalgebra, generated by
$e_{\gamma}, \gamma \in J_i$ and $f_{\theta}t$, where $f_{\theta}$ is a lowest root vector of $\g_i$.
Finally, we call $\mathcal{L}\n$ the subalgebra of $\mathcal{L}\g$ generated by $\widehat{\n}_i$ for
$i=1, \ldots, N$,
$\{e_{\star,l}\;|\; l \in \Z\}$ and $\{h_{\star,r}\;|r \geq 1\}$. The set of weights
occuring in $\mathcal{L}\n$ defines a positive cone
$\widehat{Q}^+\subset \widehat{Q}$. A root $\alpha=c_{\star} \alpha_{\star} + \sum_{\gamma \in \Gamma \backslash \{\star\}}
c_{\gamma}\alpha_\gamma + c_{\delta}\delta \in \widehat{Q}$ belongs to $\widehat{Q}^+$ if and only if
either $c_{\star} >0$ or $c_{\star}=0$, $c_{\delta}\geq 0$ and $c_{\gamma} \geq 0$ for all $\gamma$
if $c_{\delta}=0$.

\vspace{.1in}

\paragraph{\textbf{1.2. Hall algebras.}} Fix $n \in \N$. The cyclic quiver of type $A^{(1)}_{n-1}$ is the oriented
graph with set of vertices $\Z/n\Z$ and set of arrows
$\{(i,i-1)\;|i \in \Z/n\Z\}$. A nilpotent representation of
$A^{(1)}_{n-1}$ over a field $k$ consists by definition of a pair
$(V, x)$ where $V=\bigoplus_{i \in \Z/n\Z} V_i$ is a
$\Z/n\Z$-graded $k$-vector space and $x=(x_i)$ belongs to the
space
$$\mathcal{N}^{(n)}_V=\{(x_i) \in \prod_i \mathrm{Hom}(V_i,V_{i-1})\;|\;x_{i-N}\cdots x_{i-1}x_{i}=0\;\text{for
\;any\;}i\text{\;and\;}N \gg 0\}.$$
The set of all nilpotent representations of $A^{(1)}_{n-1}$ forms an abelian category of
global dimension one. We briefly recall the definition of the Hall algebra $\mathbf{H}_{n}$.
Let us take $k$ to be a finite field with $q$ elements.
For any $\Z/n\Z$-graded vector space $V$, let $\C_{G}(\mathcal{N}^{(n)}_V)$ be the set of
$G_V$-invariant functions $\mathcal{N}^{(n)}_V \to \C$
where the group $G_V=\prod_i GL(V_i)$ naturally acts on $\mathcal{N}^{(n)}_V$ by conjugation. For each
$\mathbf{d} \in \N^{\Z/n\Z}$, let us fix an $\Z/n\Z$-graded vector space $V_{\mathbf{d}}$ of dimension $\mathbf{d}$
and set $\mathcal{N}^{(n)}_{\mathbf{d}}=\mathcal{N}^{(n)}_{V_\mathbf{d}}$. Given
$\mathbf{d}',\mathbf{d}'' \in \N^{\Z/n\Z}$ such that $\mathbf{d}'+\mathbf{d}''=\mathbf{d}$, consider the diagram
$$\xymatrix{
\mathcal{N}^{(n)}_{\mathbf{d}''}\times \mathcal{N}^{(n)}_{\mathbf{d}'} &
 E' \ar[l]_-{p_1} \ar[r]^-{p_2} & E'' \ar[r]^-{p_3} &
\mathcal{N}^{(n)}_{\mathbf{d}}}$$
where the following notation is used\\
-$E'$ is the set of tuples $(x,V',a,b)$ such that $x \in \mathcal{N}^{(n)}_{{\mathbf{d}}}$, $V' \subset
 V_{\mathbf{d}}$ is an $x$-stable
subrepresentation of $V_{\mathbf{d}}$, and $a: V' \stackrel{\sim}{\to} V_{\mathbf{d}'},\;
b\;:V_{\mathbf{d}}/V' \stackrel{\sim}{\to}
V_{\mathbf{d}''}$ are isomorphisms,\\
-$E''$ is the set of tuples $(x,V')$ as above,\\
-$p_1(x,V',a,b)=( b_{*}x_{|V_{\mathbf{d}}/V'},a_{*}x_{|V'})$ and $p_2,p_3$ are the natural projections. \\

Given $f \in \C_{G}(\mathcal{N}^{(n)}_{\mathbf{d}'})$ and $g \in \C_G(\mathcal{N}^{(n)}_{\mathbf{d}''})$, set
\begin{equation}\label{E:1.2}
f \cdot g= q^{\frac{1}{2}\{ \mathbf{d}',\mathbf{d}''\}} (p_3)_{!}h \in \C_G(\mathcal{N}^{(n)}_{\mathbf{d}}),
\end{equation}
where $h\in \C(E'')$ is the unique function satisfying
$p_2^*(h)=p_1^*(fg)$ and $\{ \mathbf{d}',\mathbf{d}''\}=\sum_{i}
d'_id''_i - \sum_{i} d'_{i}d''_{i+1}$. Set $\mathbf{H}_{n} =
\bigoplus_{\mathbf{d}} \C_G(\mathcal{N}^{(n)}_{\mathbf{d}})$. It
is known that $\mathbf{H}_{n}$ equipped with the above product, is
an associative algebra, whose structure constants are integral
polynomials in $q^{\pm \frac{1}{2}}$. Hence we may consider it as
a $\C[v,v^{-1}]$-algebra where $v^{-2}=q$.

\vspace{.05in}

For $i \in \Z/n\Z$ put $\epsilon_i=(\delta_{ij}) \in \N^{\Z/n\Z}$.
If $n >1$ the variety $\mathcal{N}^{(n)}_{l\epsilon_i}$ consists
of a single point, and we let
$E_i^{(l)}=1_{\mathcal{N}^{(n)}_{l\epsilon_i}} \in \mathbf{H}_n$
stand for its characteristic function. By Ringel's theorem
(\cite{Ri}), the elements $E_i^{(l)}$ for $i \in \Z/n\Z$ and $l
\geq 1$ generate a subalgebra isomorphic to the
$\mathbb{C}[v,v^{-1}]$-form of
$\mathbf{U}^+_v(\widehat{\mathfrak{sl}}_{n})$. Now set
$\mathbf{d}=(l,\ldots, l)$. The set of $G_{V_{\mathbf{d}}}$-orbits
in $\mathcal{N}^{(n)}_{\mathbf{d}}$ for which $x_2, x_{3},\ldots,
x_{n}$ are all isomorphisms is naturally indexed by the set of
partitions of $l$ (see \cite{S1}, Section 4). Denoting by
$1_{O_{\lambda}}$ the characteristic function of the orbit
corresponding to a partition $\lambda$, we put
\begin{equation}\label{E:Hall1}
\mathbf{h}_{l}=\frac{1}{l}\sum_{\lambda, |\lambda|=l}n(l(\lambda)-1)1_{O_{\lambda}}, \qquad
n(l)=\prod_{i=1}^l (1-v^{-2i}).
\end{equation}
By \cite{S1}, Lemma~4.3, $\mathbf{H}_n$ is generated by
$E_i^{(l)}$ and $\mathbf{h}_{l}$ for $i \in \Z/n\Z$ and $l \geq
1$, and it can naturally be viewed as a quantum deformation of
$\mathbf{U}^+(\widehat{\mathfrak{gl}}_{n})$.

\vspace{.1in}

Finally, we also equip $\mathbf{H}_{n}$ with a coproduct $\Delta: \mathbf{H}_n \otimes \mathbf{H}_n
\to \mathbf{H}_n$. Let us fix, for each triple $\mathbf{d}, \mathbf{d}', \mathbf{d}''$ satisfying
$\mathbf{d}' + \mathbf{d}''=\mathbf{d}$, an exact sequence
$$\xymatrix{ 0 \ar[r] & V_{\mathbf{d'}} \ar[r]^-{a} & V_{\mathbf{d}} \ar[r]^-{b} & V_{\mathbf{d}''} \ar[r] & 0}$$
and let $F \subset \mathcal{N}_{\mathbf{d}}^{(n)}$ be the subset of representations fixing $a(V_{\mathbf{d}})$.
Consider the diagram
$$\xymatrix{ \mathcal{N}^{(n)}_{\mathbf{d}} & F \ar[l]_-{\iota} \ar[r]^-{\kappa} &
\mathcal{N}^{(n)}_{\mathbf{d}''} \times
\mathcal{N}^{(n)}_{\mathbf{d}'}}$$ where $\iota$ denotes the
embedding and
$\kappa(x)=(b_{*}(x_{|V_{\mathbf{d}}/V_{\mathbf{d}'}}),
a_*(x_{|V_{\mathbf{d}'}}))$. For any $f$ in
$\C_{G}(\mathcal{N}^{(n)}_{\mathbf{d}})$ we put
$$\Delta(f) =\sum_{\mathbf{d}', \mathbf{d}''} \Delta_{\mathbf{d}'',\mathbf{d}'}(f), \qquad
\Delta_{\mathbf{d}'', \mathbf{d}'}(f)=v^{\{\mathbf{d}',\mathbf{d}''\}}\kappa_!\iota^*(f) \in
\C_G(\mathcal{N}^{(n)}_{\mathbf{d}''})\otimes \C_{G}(\mathcal{N}^{(n)}_{\mathbf{d}'}).$$

\begin{prop}[\cite{S0}]\label{P:HallS0}
We have $\mathbf{H}_n \simeq \mathbf{U}^+_v(\widehat{\mathfrak{sl}}_n)
\otimes \C[v,v^{-1}][z_1,z_2, \ldots]$, where $z_i \in
\mathbf{H}_n$ are central elements of degree $(i, \cdots, i)$
satisfying $\Delta(z_i)=z_i \otimes 1 + 1 \otimes z_i$.\end{prop}

\vspace{.1in}

\paragraph{\textbf{Remarks.}} i) The convention $v^{-2}=q$ differs from the one used in \cite{S1}, by the change of variables
$v \mapsto v^{-1}$ (this is also valid for Section~1.3.). \\
ii)The set of
$G_{V_{\mathbf{d}}}$-orbits in $\mathcal{N}_{\mathbf{d}}^{(n)}$ is
independent of the field $k$, and we may thus define an algebra and coalgebra
structure on the space $\mathbf{H}_{n}$ for any field.

\vspace{.1in}

\paragraph{\textbf{1.3. Quantum groups.}} The enveloping algebra $\mathbf{U}(\mathcal{L}\g)$ admits a well-known
deformation which is due to Drinfeld (see \cite{Dr}). We will be
concerned here with a certain deformation of its positive part
$\mathbf{U}(\mathcal{L}\n)$. For each branch $J_i$ of $\Gamma$ we
consider one copy of the algebra $\mathbf{H}_{p_i}$ and we denote
by $E_{(i,j)}$ resp. $\mathbf{h}_{i,l}$ the generators of
$\mathbf{H}_{p_i}$ considered in Section~1.2, corresponding to the
$j$th vertex of $A^{(1)}_{p_i-1}$ or to a positive integer $l$. We
let $\mathbf{U}'_v(\mathcal{L}\mathfrak{n})$ be the
$\C(v)$-algebra generated by $\mathbf{H}_{p_i}$ for $i=1, \ldots
N$, $H_{\star,(l)}$,
 for $l \in \N^*$ and $E_{\star,k}$ for
$k \in \Z$, subject to the following set of relations :
\begin{enumerate}
\item[i)] For all $i$, $H_{\star,(l)}=v^{lp_i}\mathbf{h}_{i,l}$, and
$[E_{(i,j)},E_{(i',j')}]=0$ if $i \neq i'$.
\item[ii)] We have
\begin{equation*}
[H_{\star,(l)},E_{\star,t}]=((v^l+v^{-l})E_{\star,l+t})/l,
\end{equation*}
\begin{equation*}
v^2E_{\star,t_1+1}E_{\star,t_2}-E_{\star,t_2}E_{\star, t_1+1}
= E_{\star,t_1}E_{\star,t_2+1}-v^2E_{\star,t_2+1}E_{\star,t_1},
\end{equation*}
\item[iii)]
\begin{equation*}
a_{\star,\gamma}=0 \Rightarrow [E_{\star,t},E_{\gamma}]=0, \qquad
E_{\star,t}E_{(i,0)}-vE_{(i,0)}E_{\star,t}=0,
\end{equation*}
\begin{equation*}
p_i >2 \Rightarrow [E_{(i,0)}, E_{(i,1)}E_{\star,t}-vE_{\star,t}
E_{(i,1)}]=0,
\end{equation*}
\item[iv)] For $\gamma=(i,1)$ with $i=1, \ldots, N$, and $l \geq 1$ set $E_{\gamma,l}=
l[E_{\gamma},H_{\star,(l)}]$. Then for any $l,r_1,r_2 \in \N$ and
$n,t_1,t_2 \in \Z$,
\begin{equation*}
\mathrm{Sym}_{r_1,r_2}
\big\{E_{\gamma,r_1} E_{\gamma,r_2}E_{\star,n}
-(v+v^{-1})E_{\gamma,r_1}E_{\star,n}E_{\gamma,r_2} +E_{\star,n}
E_{\gamma,r_1} E_{\gamma,r_2}\big\}=0,
\end{equation*}
\begin{equation*}
\mathrm{Sym}_{t_1,t_2} \big\{E_{\star,t_1}E_{\star,t_2}E_{\gamma,l}
-(v+v^{-1})E_{\star,t_1}E_{\gamma,l}E_{\star,t_2}+E_{\gamma,l}
E_{\star,t_1}E_{\star,t_2}\big\}=0,
\end{equation*}
\item[v)] For $\gamma=(i,1)$ with $i=1, \ldots, N$ and $l,r_1,r_2 \in \N$, $t \in \Z$ we have
\begin{equation*}
E_{\star, t+1}E_{\gamma,l}-vE_{\gamma,l}E_{\star,t+1}=
vE_{\star,t}E_{\gamma,l+1}-E_{\gamma,l+1}E_{\star,t},
\end{equation*}
\begin{equation*}
v^2E_{\gamma,r_1+1}E_{\gamma,r_2}-E_{\gamma,r_2}
E_{\gamma,r_1+1}=E_{\gamma,r_2+1}E_{\gamma,r_1}
-v^2E_{\gamma,r_1}E_{\gamma,r_2+1}.
\end{equation*}
\end{enumerate}
We also let $\mathbf{U}_v(\mathcal{L}\n)$ be the
$\C[v,v^{-1}]$-subalgebra of $\mathbf{U}_v(\mathcal{L}\n)$
generated by $\mathbf{H}_{p_i}$ for $i=1, \ldots, N$ and the
divided powers $E_{\star,t}^{(n)}=E_{\star,t}^n/[n]!$ (here as usual we put
$[l]=(v^l-v^{-l})/(v-v^{-1})$ and $[l]!=[2]\cdots[l]$). Lastly, it
is convenient to introduce another set of elements of
$\mathbf{U}_v(\mathcal{L}\n)$~: define $\xi_l$ for $l \geq 1$ by
the formal relation $1+ \sum_{l\geq 1} \xi_l s^l=exp(\sum_{l \geq
1} H_{\star,(l)}s^l)$. The sets $\{\xi_l\}$ and
$\{H_{\star,(l)}\}$ span the same subalgebra over $\C[v,v^{-1}]$.

\vspace{.1in}

\paragraph{\textbf{Remark.}} The relations are slightly renormalized from those in \cite{S1}.
In particular, the element $H_{\star,(l)}$ corresponds to
$H_{\star,l}/[l]$ in \cite{S1}. The elements $\xi_l$ are sometimes
denoted $\tilde{P}_l$ in the litterature.

\vspace{.2in}

\section{Coherent sheaves and Quot schemes}

\vspace{.1in}

\paragraph{\textbf{2.1. Equivariant coherent sheaves.}} Let $(X,G)$ be a pair
consisting of a smooth
projective curve $X$
over a field $k$ and a finite group $G\subset Aut(X)$ such that $X/G
\simeq \mathbb{P}^1$.
The category
of $G$-equivariant coherent sheaves on $X$ will be denoted by
$Coh_G(X)$. It is an abelian
category of global dimension one. The quotient map $\pi:X \to \mathbb{P}^1$
is ramified at
points, say $\lambda_1, \ldots,
\lambda_N \in \mathbb{P}^1$ with respective indices of ramification
$p_1, \ldots,
p_N$. We assume that all the points $\lambda_i$ lie in $X(k)$.
Let $\Gamma=\Gamma_{X,G}$ be the star-shaped Dynkin diagram with
branches of length
$p_1-1,\ldots, p_N-1$ (if $(X,G)=(\mathbb{P}^1,Id)$ then $\Gamma$ is
of type $A_1$).
We also put $\Lambda=\{\lambda_1,\ldots, \lambda_N\} \subset \mathbb{P}^1$.

\vspace{.1in}

Any $\mathcal{F} \in Coh_G(X)$ has a canonical torsion subsheaf
$\tau(\mathcal{F})$ and
locally free quotient $\nu(\mathcal{F})$ and there is a (noncanonical) isomorphism
$\mathcal{F} \simeq \tau(\mathcal{F}) \oplus
\nu(\mathcal{F})$.

\vspace{.1in}

Let us first describe the torsion sheaves in $Coh_G(X)$.
Let $\Ox$ be the structure sheaf of $X$ with the trivial $G$-structure. For
each closed point $x \not\in \Lambda$ there exists, up to
isomorphism, a unique simple torsion sheaf $\O_x$ with support in
$\pi^{-1}(x)$, and
there holds
$$\mathrm{dim\;Ext}(\O_x,\O_x)=\mathrm{dim\;Ext}
(\O_x, \mathcal{O}_X)=\mathrm{dim\;Hom}(\Ox, \O_x)
=1.$$
In particular, the full subcategory $\mathcal{T}_x$ of torsion
sheaves supported at $x$
is equivalent to the category of nilpotent representations (over the
residue field $k_x$ at $x$)
of the cyclic quiver $A_0^{(1)}$ with one vertex and one loop.

On the other hand, for each ramification point $\lambda_i \in \Lambda$ of index
$p_i$, there exists
$p_i$ simple torsion sheaves $S_i^{(j)}$, $j \in \Z/p_i\Z$ with
support in $\pi^{-1}(\lambda_i)$. We have
$$\mathrm{dim\;Ext}(S_i^{(j)},S_i^{(l)})=\delta_{j-1,l},$$
$$\mathrm{dim\;Ext}
(S_i^{(j)}, \mathcal{O}_X)=\delta_{j,1},\qquad
\mathrm{dim\;Hom}(\Ox, S_i^{(l)}) =\delta_{j,0}.$$ Hence, the
category $\mathcal{T}_{\lambda_i}$ is equivalent to the category
of nilpotent representations (over $k$) of the cyclic quiver
$A_{p_i-1}^{(1)}$. Under this equivalence, the sheaf $S_i^{(j)}$
goes to the simple module $(V,x)$ where $V_j=k$ and $V_h=\{0\}$ if
$j \neq h$. We denote by $\aleph=\{(i,j)\;|i =1, \ldots N, j \in
\Z/p_i\Z\}$ the set of simple exceptional torsion sheaves.

\vspace{.1in}

The Picard group of $(X,G)$ is also easy to describe. Let
$\mathcal{L}_{\delta}$
(resp. $\mathcal{L}_i$) be the nontrivial extension of $\Ox$ by a
generic simple
torsion sheaf $\O_x$
for $x \not\in \Lambda$ (resp. by the simple torsion sheaf
$S_i^{(1)}$). Then $\mathcal{L}_{\delta}$ and $\mathcal{L}_i$ for
$i=1, \ldots, N$
generate $Pic(X,G)$ with the relations $\mathcal{L}_i^{\otimes p_i} \simeq
\mathcal{L}_{\delta}$. Thus
$$Pic(X,G)/\Z[\mathcal{L}_{\delta}] \simeq \prod_i \Z/p_i\Z.$$
In particular, $Pic(X,G)$ is of rank one.

\vspace{.1in}

The category $Coh_G(X)$ inherits Serre duality from the
category $Coh(X)$. Let $\omega_X$ be the canonical sheaf of $X$ with
its natural
$G$-structure. Then for any $\mathcal{F}, \mathcal{G} \in Coh_G(X)$
there is a (functorial)
isomorphism $\Ext(\mathcal{F},\mathcal{G})^* \simeq \Hom(\mathcal{G},
\mathcal{F} \otimes \omega_X )$. It is known that (see \cite{GL})
\begin{equation}\label{E:0}
\omega_X\simeq
\mathcal{L}_{\delta}^{\otimes N-2} \otimes \bigotimes_i \mathcal{L}_i^{-1}
\end{equation}
from which it follows that $H^0(\omega_X)=0$.

\vspace{.1in}

Write $K(X)$ for the Grothendieck group of $Coh_G(X)$.
This group is equipped with the Euler form
$\langle([M],[N]\rangle=\mathrm{dim\;Hom}\;(M,N)-\mathrm{dim\;Ext}(M,N)$
and its symmetrized
version $(M,N)=\langle M,N\rangle+ \langle N,M\rangle$. The set of
classes $\alpha \in K(X)$
for which there exists some sheaf $\mathcal{F}$ with
$[\mathcal{F}]=\alpha$ forms a cone of
$K(X)$ denoted $K^+(X)$. Then (see \cite{S1}, Prop~ 5.1)

\vspace{.05in}

\begin{lem}\label{L:00}
There is a canonical isomorphism $h:K(X) \simeq \widehat{Q}$ compatible
with the symmetric forms $(\,,\,)$. It maps $K^+(X)$ to the positive
weight lattice
$\widehat{Q}^+$.
\end{lem}

\vspace{.1in}

Under this isomorphism, we have
$h([\Ox])=\alpha_\star$, $h([\O_x])=\delta$ for any closed point $x
\not\in \Lambda$
of degree one, $h([S_i^{(j)}])=\alpha_{(i,j)}$ for $j \neq 0$ and
$h([S_i^{(0)}])=\delta-\sum_{j=1}^{p_i-1} \alpha_{(i,j)}$. In particular, the
Picard group $Pic(X,G)$ is in correspondence with the set of roots of
$\mathcal{L}\g$
of the form $\alpha_\star+\sum_{i,j} l_{i,j} \alpha_{(i,j)}$.
Define a finite set
$S=\{\alpha_\star\} \cup \{\alpha_\star+
\sum_{k=1}^j\alpha_{(i,k)}\;|(i,j) \in \Gamma\backslash \{\star\}\}$.
For $s \in S$ we denote by
$\mathcal{L}_s$ the line bundle corresponding to $s$. Observe that by
(\ref{E:0}) we have
\begin{equation}\label{E:00}
\Ext(\mathcal{L}_s, \mathcal{L}_{s'})=0\qquad \text{for\;any\;}s, s' \in S.
\end{equation}
For any $\mathcal{F} \in Coh_G(X)$ and $n \in \Z$ we put
$\mathcal{F}(n)
=\mathcal{F} \otimes \mathcal{L}_{\delta}^{\otimes n}$.

\vspace{.1in}

\paragraph{}Lastly, we define the degree map $deg: K(X) \to
K(X)/\Z[\Ox], \alpha \mapsto \alpha + \Z [\Ox]$ and write
$deg(\alpha) \geq deg(\beta) $ if $deg(\alpha)-deg(\beta) \in
\Sigma_{i,j} \mathbb{N} [S_i^{(j)}] + \Z[\Ox]$. 
We set $\omega=deg[\omega_X]$. We
will also need a degree map taking values in $\mathbb{Q}$~: if
$\alpha \in K(X)$ and $deg(\alpha)=\sum_{i,j} n_{i,j}\alpha_{(i,j)}$
we put $|\alpha|=\sum_i (n_{i,j}/p_i) \in \mathbb{Q}$. We extend
both of these notations to a coherent sheaf $\mathcal{F}$ by
setting $deg(\mathcal{F})=deg([\mathcal{F}])$ and similarly for
$|\;|$.

\vspace{.1in}

\paragraph{\textbf{Remarks.}} i)
The category $Coh_G(X)$ admits an equivalent description as
the category
of \textit{$D$-parabolic coherent sheaves} on $\mathbb{P}^1$ where $D=\sum_i
(p_i-1)\lambda_i$ is the ramification divisor of the quotient map $X \to
\mathbb{P}^1$ (see
\cite{Bi}) . It is also equivalent to the category of coherent sheaves on a
weighted projective line in the sense of Geigle and Lenzing (see
\cite{GL} where
this category was first singled out and studied in details; in fact
the notion of a weighted
projective line is slightly more general as it also covers the case of diagrams
of type $A_{2n}$).

\vspace{.05in}

ii) Given any positive divisor $D =\sum_{i=1}^N (p_i-1)\lambda_i$
with $N \geq 3$
there exists a
smooth projective curve $X$ (over $\C$) and a group
$G$ of automorphisms of $X$ such that $X \to X/G \simeq \mathbb{P}^1$ has
ramification locus
$D$. This can be proved as follows. From Galois theory, it is enough
to find a finite
  group $G$
generated by elements $y_i$ of order $p_i$ such that $y_1 \cdots y_N=Id$.
Consider matrices in $SL(2, \mathbb{F})$
$$y_i=\left(\begin{matrix}\zeta_{p_i}& x_i \\ 0 &
\zeta_{p_i}^{-1}\end{matrix}\right)
\qquad \text{for\;}i=1, \ldots, N-2,$$
$$y_{N-1}=\left(\begin{matrix} \zeta_{p_{N-1}} & 0\\ x_{N-1} &
\zeta_{p_{N-1}}^{-1}
\end{matrix}
\right),\qquad  y_N^{-1}=y_1 \cdots y_{N-1},$$
where $\zeta_l$ is a primitive $l$th root of unity in some big enough
finite field
$\mathbb{F}$. Each $y_i$, $i=1, \ldots, N-1$ is of order $p_i$, and the product
$y_1 \cdots y_{N-1}$ has determinant one and trace equal
to a certain nonconstant polynomial in the $x_i$'s. In particular, replacing
$\mathbb{F}$ by a finite extension if necessary, we can choose $x_i$, $i=1,
\ldots , N-1$ such that $Tr(y_N^{-1})=\zeta_{p_N} + \zeta_{p_N}^{-1}$. But then
$y_N$ is diagonalizable and of order $p_N$ as desired.

\vspace{.1in}

\paragraph{\textbf{2.2. Equivariant Quot schemes.}}
  We will say that $\mathcal{F} \in Coh_G(X)$ is generated
by a collection of sheaves $\{\mathcal{T}_i\}_{i \in I}$ if the natural map
$$\phi:\;\bigoplus_i \Hom(\mathcal{T}_i,\mathcal{F}) \otimes \mathcal{T}_i
\to \mathcal{F}$$ is surjective.

\begin{lem}\label{L:1} For any $\mathcal{E} \in Coh_G(X)$ there
exists $l \in \Z$ such that for any
$m \leq l$ we have
\begin{enumerate}
\item[i)] $\mathcal{E}$ is generated by $\{\mathcal{L}_s(m)\}_{s \in S}.$
\end{enumerate}
Moreover, this also implies that
\begin{enumerate}
\item[ii)] $\Ext(\mathcal{L}_s(m), \mathcal{E})=0$ for any $s \in S$.
\end{enumerate}
\end{lem}

\noindent
\textit{Proof.} Assume $\mathcal{E}'$ and $\mathcal{E}''$ both
satisfy $i)$ and $ii)$ with
corresponding integers $l', l''$. Let $0 \to \mathcal{E}' \to
\mathcal{E}
\to \mathcal{E}'' \to 0$ be a short exact sequence and fix $m \leq
inf\{l',l''\}$.
Applying $\Hom(\mathcal{L}_s(m), \cdot)$ and using $ii)$ we obtain
an exact sequence
$$0 \to \Hom(\mathcal{L}_s (m), \mathcal{E}') \to  \Hom(\mathcal{L}_s
(m), \mathcal{E})
\to  \Hom(\mathcal{L}_s (m), \mathcal{E}'') \to 0$$
and $\Ext(\mathcal{L}_s(m), \mathcal{E})=0$. The image of the natural map
$$\phi: \bigoplus_s \Hom(\mathcal{L}_s(m), \mathcal{E})\otimes
\mathcal{L}_s(m) \to
\mathcal{E}$$ contains
$\mathcal{E}'$ and maps surjectively on $\mathcal{E}''$ so that
$\phi$ is surjective.
Hence if $i), ii)$ hold for two sheaves then it also holds for any
extension between them.

\vspace{.1in}

Any torsion sheaf is obtained as an extension of the simple sheaves
$\mathcal{O}_x$
for $x \not\in \Lambda$ and $S_i^{(j)}$ for $(i,j) \in \aleph$.
Moreover, for any $x \not\in \Lambda$ there is an exact sequence
$0 \to \Ox(-1) \to \Ox \to \mathcal{O}_x \to 0$ and likewise for any
$(i,j)$ there is
an exact sequence $0 \to \mathcal{L}_{s} \to \mathcal{L}_{s'} \to
S_i^{j} \to 0$ for suitable
$s,s' \in S$. Thus the statement of the Lemma holds for all torsion
sheaves. Note that
$ii)$ holds for an arbitrary torsion sheaf and arbitrary $m \in \Z$.
Similarly, by \cite{GL}, Proposition 2.6. any vector bundle is an
extension of line bundles and,
in turn, for any line bundle $\mathcal{L}$
there exists $l \in \N$ such that $\mathcal{L}(l)$ is an extension of
$\Ox$ by a torsion
sheaf. In particular, any vector bundle is generated by
$\{\mathcal{L}_s(m)\}$ for some $m$.
Conversely, if a vector bundle $\mathcal{E}$ is generated by
$\{\mathcal{L}_s(m)\}$ then it can be obtained as an extension of
line bundles $\mathcal{V}_i$
for $i=1, \ldots, r$ with $deg( \mathcal{V}_i) \geq deg
(\mathcal{L}_s(m))$ for any $s \in S$.
But then using (\ref{E:00}) we deduce that $\Ext(\mathcal{L}_s(m),
\mathcal{V}_i)=0$ for all
$s \in S$ and the Lemma follows. \qed

\begin{lem}\label{L:2} Let $\mathcal{F} \in Coh_G(X)$ be generated by
$\{\mathcal{L}_s(n_0)\}_{s \in S}$
and fix $\alpha \in K^+(X)$. There exists an integer $n
=q(\alpha,\beta,n_0) \in \Z$
depending only on $n_0, \alpha$ and $\beta=[\mathcal{F}]$ such that any
subsheaf $\mathcal{G} \subset \mathcal{F}$ with
$[\mathcal{G}]=\alpha$ satisfies
\begin{enumerate}
\item[i)] $\mathcal{G}$ is generated by $\{\mathcal{L}_s(m)\}_{s \in S}$,
\item[ii)] $\Ext(\mathcal{L}_s(m), \mathcal{G})=0$ for any $s \in S$
\end{enumerate}
for any $m \leq n$.
\end{lem}

\noindent
\textit{Proof.} If $\mathcal{F}$ is generated by
$\{\mathcal{L}_s(n_0)\}_{s \in S}$ then
so is $\nu(\mathcal{F})$. In particular, $\nu(\mathcal{F})$ posseses
a filtration
$$0 \subset \mathcal{V}_1 \subset \cdots \subset
\mathcal{V}_r=\nu(\mathcal{F})$$
with successive quotients $\mathcal{V}_i/\mathcal{V}_{i-1}$ being line bundles
of degree at least $deg(\mathcal{L}_{\star}(n_0))=n_0\delta$. Thus we have
$deg(\tau(\mathcal{F}))=deg(\mathcal{F})-deg(\nu(\mathcal{F})) \leq
deg(\mathcal{F}) -rank(\mathcal{F}) n_0 \delta$. A similar reasoning shows that
\begin{equation}\label{E:1}
\Hom(\mathcal{L}_{\star}(l), \nu(\mathcal{F})) \neq 0 \Rightarrow
l\delta
\leq deg(\mathcal{F}) - (rank(\mathcal{F})-1) n_0\delta.
\end{equation}
Now fix $\mathcal{G} \subset \mathcal{F}$, $[\mathcal{G}]=\alpha$. Then
$deg(\tau(\mathcal{G})) \leq deg(\tau(\mathcal{F}))$ and therefore
$deg(\nu(\mathcal{G}))
\geq deg(\alpha) - deg(\tau(\mathcal{F})) \geq deg(\alpha) +
rank(\mathcal{F})n_0\delta -
deg(\mathcal{F})$. Consider a filtration
$$0 \subset \mathcal{G}_1 \subset \cdots \subset \mathcal{G}_s =\mathcal{G}$$
of $\mathcal{G}$ with the property that
$\mathcal{G}_i/\mathcal{G}_{i-1}$ is a line subbundle
of $\mathcal{G}/\mathcal{G}_{i-1}$ of maximal possible degree.

\vspace{.05in}

\noindent
\begin{claim} We have $deg(\mathcal{G}_i/\mathcal{G}_{i-1}) \leq
deg(\mathcal{F})
-(rank(\mathcal{F})-1)n_0\delta + (i-1)(3\delta -\omega)$.
\end{claim}
\noindent
\textit{Proof of claim.} We have $deg(\mathcal{G}_1) \leq deg(\mathcal{F}) -
(rank(\mathcal{F})-1)n_0\delta$ by (\ref{E:1}). Now let $0 \to \mathcal{H}' \to
\mathcal{H} \to \mathcal{H}'' \to 0$ be a short exact sequence where
$\mathcal{H}',
\mathcal{H}''$ are line bundles. Assume that
\begin{equation}\label{E:2}
deg(\mathcal{H}'') > deg(\mathcal{H}') + 3\delta
-\omega.
\end{equation}
Applying $\Hom(\mathcal{H}'(1),\cdot)$ we obtain a sequence
$$0 \to \Hom(\mathcal{H}'(1),\mathcal{H}) \to
\Hom(\mathcal{H}'(1),\mathcal{H}'')
\to \Ext(\mathcal{H}'(1),\mathcal{H}').$$
By Serre duality and (\ref{E:2}) we have
$$\mathrm{dim}\;\Ext(\mathcal{H}'(1), \mathcal{H}')=
\mathrm{dim}\;\Hom(\mathcal{H}', \mathcal{H}'(1) \otimes \omega_X) <
\mathrm{dim}\;
\Hom(\mathcal{H}'(1),\mathcal{H}'').$$
Thus $\Hom(\mathcal{H}'(1),\mathcal{H}) \neq 0$.
In particular, if (\ref{E:2}) holds then $\mathcal{H}'$ is not a
maximal degree line subbundle of $\mathcal{H}$. The claim easily follows by induction.\qed

\vspace{.05in}

 From the claim and from $\sum_i
deg(\mathcal{G}_i/\mathcal{G}_{i-1})=deg\;\alpha$
we deduce that
\begin{equation*}
\begin{split}
deg(\mathcal{G}_i/\mathcal{G}_{i-1}) \geq deg(\alpha) -(rank(\alpha)
-1)&\big\{deg(\beta)
-(rank(\beta) -1)n_0\delta \\
&\qquad + (rank(\alpha) -1)(3\delta -\omega)\big\}.
\end{split}
\end{equation*}
Let $Q(\alpha,\beta,n_0)$ denote the right hand side of the above
inequality. Let
$q(\alpha,\beta,n_0)$ be the greatest integer such that $q(\alpha,\beta,n_0) <
Q(\alpha,\beta,n_0)-\omega$. Then each
$\mathcal{G}_i/\mathcal{G}_{i-1}$ is generated
by $\{\mathcal{L}_s(q(\alpha,\beta,n_0))\}_{s \in S}$ and
$\Ext^1(\mathcal{L}_s(q(\alpha,\beta,n_0)),\mathcal{G}_i/\mathcal{G}_{i-1})
=0$ for
all $s \in S$. The Lemma follows.\qed

\vspace{.1in}

\paragraph{}We assume
from now on the ground field $k$ to be algebraically
closed. Let $\mathcal{E} \in Coh_G(X)$ and $\alpha \in K^+(X)$.
Consider the functor from the category of smooth schemes over $k$
to the category of sets defined as follows.
\begin{equation*}
\begin{split}
\underline{\mathrm{Hilb}}_{\mathcal{E},\alpha}(\Sigma)=&\{\phi:
\mathcal{E} \boxtimes \O_\Sigma
\twoheadrightarrow \mathcal{F}\;|
\mathcal{F}\;is\;a\;G-equivariant,\;coherent,\;\Sigma-flat\\
& \qquad \qquad \qquad sheaf,\;
\mathcal{F}_{|\sigma}\;is\;
of\;class\;\alpha\;for\;all\;closed\;points\;\sigma \in \Sigma\}
\end{split}
\end{equation*}
In the above, we consider two maps $\phi,\phi'$ equivalent if their
kernels coincide.

\vspace{.05in}

\begin{prop}\label{P:1} The functor
$\underline{\mathrm{Hilb}}_{\mathcal{E},\alpha}$ is represented
by a projective scheme $\mathrm{Hilb}_{\mathcal{E},\alpha}$.\end{prop}

\noindent
\textit{Proof.} This result is known in the (equivalent) situation of parabolic
bundles on a curve (see \cite{P} and Remark 1.1). We will only sketch
a construction
of $\mathrm{Hilb}_{\mathcal{E},\alpha}$. Fix $n_0 \in \Z$ such that
$\mathcal{E}$ is generated
by $\{\mathcal{L}_s(n_0)\}_{s \in S}$ and choose an integer
$n < \inf\{n_0,q([\mathcal{E}]-\alpha,[\mathcal{E}],
n_0)\}$. To a map $\phi: \mathcal{E} \twoheadrightarrow \mathcal{F}$
satisfying $[\mathcal{F}]=\alpha$
we associate the collection of
subspaces $\Hom(\mathcal{L}_s(n),\mathrm{Ker}\;\phi) \subset
\Hom(\mathcal{L}_s(n),\mathcal{E})$,
for $s \in S$. By Lemma~\ref{L:2} we have
$\Ext(\mathcal{L}_s(n),\mathrm{Ker}\;\phi)=0$ for all $s$, so that
$\mathrm{dim}\;\Hom(\mathcal{L}_s(n),\mathrm{Ker}\;\phi)=\langle
[\mathcal{L}_s(n)], [\mathcal{E}]-
\alpha \rangle$. This defines a map to a Grassmannian
$$ \iota:\; \mathrm{Hilb}_{\mathcal{E},\alpha} \to \prod_s
\mathrm{Gr}\big(\langle [\mathcal{L}_s],
[\mathcal{E}] -\alpha \rangle, \langle [\mathcal{L}_s],
[\mathcal{E}] \rangle\big).$$
Conversely, define a map $\psi: \prod_s \mathrm{Gr}\big(\langle
[\mathcal{L}_s],
[\mathcal{E}] -\alpha \rangle, \langle [\mathcal{L}_s],
[\mathcal{E}]\rangle \big)
  \to \{\mathcal{G} \subset \mathcal{E}\}$
by assigning to $(V_s)_s$ the image of the composition
$$\bigoplus_s V_s \otimes \mathcal{L}_s \to \bigoplus_s
\Hom(\mathcal{L}_s(n),\mathcal{E})
\otimes \mathcal{L}_s \twoheadrightarrow \mathcal{E}.$$
By Lemma~\ref{L:2}, every subsheaf of $\mathcal{E}$ of class
$[\mathcal{E}]-\alpha$ is generated by
$\{\mathcal{L}_s(n)\}_{s \in S}$ and therefore $\psi \circ \iota=
Id$. Thus $\iota$ is injective.

\vspace{.05in}

The image of $\iota$ can be explicitly described. Let $H$ be the subset of
the Grassmannian
$\prod_s \mathrm{Gr}\big(\langle [\mathcal{L}_s],[\mathcal{E}]
-\alpha \rangle, \langle [\mathcal{L}_s]
, [\mathcal{E}]\rangle \big)$
consisting of elements $(V_s)_s$ such that the image of the composition
\begin{equation*}
\begin{split}
\bigoplus_s\Hom(\mathcal{L}_{s'}(n'),\mathcal{L}_s(n)) \otimes V_s &\to
\bigoplus_s\Hom(\mathcal{L}_{s'}(n'),\mathcal{L}_s(n)) \otimes
\Hom(\mathcal{L}_s(n), \mathcal{E})\\
&\twoheadrightarrow \Hom(\mathcal{L}_{s'}(n'), \mathcal{E})
\end{split}
\end{equation*}
is of dimension $\langle [\mathcal{L}_{s'}(n')], [\mathcal{E}]-\alpha
\rangle$ for all $s' \in S$
and $n'<n$. It is clear that
$\iota(\mathrm{Hilb}_{\mathcal{E},\alpha}) \subset H$. Conversely,
if $(V_s)_s \in H$ then for all $n'<n$ and $s' \in S$ we have
$\mathrm{dim}\;\Hom(\mathcal{L}_{s'}(n'), \psi((V_s)_s))=\langle
[\mathcal{L}_{s'}(n')], [\mathcal{E}]
-\alpha \rangle$. On the other hand, by Lemma~\ref{L:1},
$\langle [\mathcal{L}_{s'}(n')], \psi((V_s)_s)\rangle =\langle
[\mathcal{L}_{s'}(n')], [\mathcal{E}]-\alpha \rangle$ for $n' \ll 0$.
This implies that $[\psi((V_s)_s)]=[\mathcal{E}]-\alpha$. Thus
$\iota(\mathrm{Hilb}_{\mathcal{E},\alpha})=H$.

\vspace{.1in}

The fact that $H$ can be equipped with the structure of a projective
algebraic variety and that it satisfies
the universal property is proved in the same manner as in \cite{LP},
Chap. 4.\qed

\vspace{.1in}

\paragraph{}Let $\alpha \in K^+(X)$ and $n \in \Z$. For $s \in S$
we put $d_s(n,\alpha)=\langle[\mathcal{L}_s(n)],\alpha\rangle$
and
$$\mathcal{E}_{n}^{\alpha}= \bigoplus_{s \in S}
\mathcal{L}_s(n)\otimes k^{d_s(n,\alpha)}.$$
Consider the subfunctor
$\underline{\mathrm{Hilb}}^0_{\mathcal{E}_{n}^{\alpha},\alpha}$
of $\underline{\mathrm{Hilb}}_{\mathcal{E}_{n}^{\alpha},\alpha}$ given by the
assignement
\begin{equation*}
\begin{split}
\Sigma \mapsto \big\{(\phi: \mathcal{E}_{n}^{\alpha} \otimes \O_\Sigma
&\twoheadrightarrow \mathcal{F}) \in
\underline{\mathrm{Hilb}}_{\mathcal{E}_n^\alpha,\alpha}\;|\;
\forall\;\sigma \in \Sigma, \forall\; s \in S\;\\
&\phi_{*,|\sigma}: \Hom(\mathcal{L}_s(n),\mathcal{L}_s(n)\otimes
k^{d_s(n,\alpha)})
\stackrel{\sim}{\to}
\Hom(\mathcal{L}_s(n),\mathcal{F}_{|\sigma})\big\}.
\end{split}
\end{equation*}

This functor is represented by an open subvariety $Q_n^{\alpha} \subset
\mathrm{Hilb}_{\mathcal{E}_{n}^{\alpha},\alpha}$.
The group $\mathrm{Aut}(\mathcal{E}_{n}^{\alpha})$
naturally acts on $\mathrm{Hilb}_{\mathcal{E}_{n}^{\alpha},\alpha}$
and $Q_n^\alpha$.
This group is non semisimple, and we will only consider the action of
its semisimple part
$$G_{n}^{\alpha}= \prod_s \mathrm{Aut}(\mathcal{L}_s (n)\otimes
k^{d_s(n,\alpha)}) \simeq
\prod_s GL(d_{s}(n,\alpha)).$$

\begin{lem}\label{L:3} There is a one-to-one correspondence
$\mathcal{F} \leftrightarrow
{O}^{\mathcal{F}}_n$ between
isomorphism classes of coherent sheaves of class $\alpha$
generated by $\{\mathcal{L}_s(n)\}$ and
$G_{n}^{\alpha}$-orbits on $Q_n^{\alpha}$. Moreover, the stabilizer
of any point $x
\in {O}^{\mathcal{F}}_n$ is isomorphic to $\mathrm{Aut}(\mathcal{F})$.
\end{lem}
\noindent
\textit{Proof.} Let $\phi: \mathcal{E}_{n}^{\alpha} \twoheadrightarrow
\mathcal{F}$ and $\phi': \mathcal{E}_{n}^{\alpha} \tto \mathcal{F}'$
be two points
in $Q_n^{\alpha}$. It is clear that if they are in the same
$G_{n}^{\alpha}$-orbit then $\mathcal{F} \simeq \mathcal{F}'$. On the
other hand, let
$u: \mathcal{F} \stackrel{\sim}{\to} \mathcal{F}'$ be an isomorphism.
This induces
isomorphisms $u_*:\Hom(\mathcal{L}_s(n),\mathcal{F})\stackrel{\sim}{\to}
\Hom(\mathcal{L}_s(n),\mathcal{F}')$. Define an element $(g_s) \in
G_{n}^{\alpha}$ as
the composition
$$g_s:k^{d_{s}(n,\alpha)} \stackrel{\phi_*}{\lto}
\Hom(\mathcal{L}_s(n),\mathcal{F})\stackrel{u_{*}}{\lto}
\Hom(\mathcal{L}_s(n),\mathcal{F}')\stackrel{(\phi'_*)^{-1}}{\lto}
k^{d_s(\alpha,n)}.$$
It is easy to see that $(g_s)$ conjugates $\phi$ to $\phi'$. This
proves the first
statement in the Lemma. Now let $(g_s)$ be an element in the
stabilizer of $\phi$. By
definition, $(g_s)(\mathrm{Ker}\;\phi)=\mathrm{Ker}\;\phi$, hence
$(g_s)$ induces
an automorphism of $\mathcal{F}$. This gives a map
$\theta:\mathrm{Stab}\;\phi \to \mathrm{Aut}(\mathcal{F})$.
The inverse map is constructed as follows. As above, an
automorphism $u \in \mathrm{Aut}(\mathcal{F})$ induces automorphisms
$u_* \in \mathrm{Aut}(\Hom(\mathcal{L}_s(n), \mathcal{F}))$. Set
$\chi(u)= (\phi_*^{-1} u_* \phi_*)$. Then
$\chi \circ \theta=Id$ and $\theta \circ \chi=Id$, and the conclusion follows.\qed

\begin{lem}\label{L:4} The variety $Q_n^\alpha$ is smooth.\end{lem}
\noindent
\textit{Proof.} By an equivariant version of \cite{LP}, Corrolary~8.10, a point
$\phi: \mathcal{E}_{n}^{\alpha} \tto \mathcal{F}$ is smooth if
$\Ext(\mathrm{Ker}\;\phi,
\mathcal{F})=0$. Observe that if $\phi \in Q_n^\alpha$ then for all
$s \in S$ we have
$\langle [\mathcal{L}_s(n)], [\mathcal{F}]
\rangle=\mathrm{dim}\;\Hom(\mathcal{L}_s(n),
\mathcal{F})$, hence $\Ext(\mathcal{L}_s(n),\mathcal{F})=0$. Applying
$\Hom(\cdot,
\mathcal{F})$ to the exact sequence $0 \to \mathrm{Ker}\;\phi \to
\mathcal{E}_{n,\alpha}
\to \mathcal{F} \to 0$ yields an exact sequence
$$\Ext(\mathcal{F},\mathcal{F}) \to \Ext(\mathcal{E}_{n}^{\alpha},\mathcal{F})
\to \Ext(\mathrm{Ker}\;\phi,\mathcal{F}) \to 0$$
from which we conclude $\Ext(\mathrm{Ker}\;\phi,\mathcal{F})=0$ as desired.\qed

\vspace{.2in}

\section{Transfer functors}

\vspace{.1in}

In this section we set up the construction of the ``inductive limit'' of the varieties $Q_n^\a$
as $n$ tends to $-\infty$. More precisely, we define the category of (equivariant) perverse sheaves
and complexes on such an inductive limit.

\vspace{.1in}

\paragraph{\textbf{3.1. Notations.}} We use notations of \cite{L1}, Chapter 8
regarding perverse sheaves. In particular, for an algebraic
variety $X$ defined over $k$ we denote by $\mathcal{D}(X)$ (resp.
$\mathcal{Q}(X)$) the derived category of
$\overline{\Q_l}$-constructible sheaves (resp. the category of
semisimple $\overline{\Q_l}$-constructible complexes of geometric
origin) on $X$. The Verdier dual of a complex $P$ is denoted by
$D(P)$. If $G$ is a connected algebraic group acting on $X$ then
we denote by $\mathcal{Q}_G(X)$ the category of semisimple
$G$-equivariant complexes. Recall that if $H
\subset G$ is an algebraic subgroup and if $Y$ is an $H$-variety
then there are canonical (inverse) equivalences
$ind_H^G:\mathcal{Q}_H(Y)[\mathrm{dim}\;G/H] \to \mathcal{Q}_G(G
\underset{H}{\times} Y)$ and $res_G^H:\mathcal{Q}_G(G
\underset{H}{\times} Y) \to\mathcal{Q}_H(Y)[\mathrm{dim}\;G/H]$ which commute with Verdier 
duality. These are characterized (up to isomorphism) as follows~: let
$$\xymatrix{
Y & G \times Y \ar[l]_-{g} \ar[r]^-{f}& G \underset{H}{\times} Y}$$
be the natural diagram; if $P \in \mathcal{Q}_H(Y)[\mathrm{dim}\; G/H]$ then
$f^*ind_H^G(P)\simeq g^*P$. Conversely, if $Q \in \mathcal{Q}_G(G
\underset{H}{\times} Y)$ then
$res_G^H(Q)=i^*Q$ where $i: Y \to G \underset{H}{\times} Y$ is the
canonical embedding.
Finally, if $Y \subset X$ then $\overline{Y}$ stands for the Zariski closure of
$Y$ in $X$.

\vspace{.2in}

\paragraph{\textbf{3.2. Transfer functors.}} Fix $\alpha \in
K^+(X)$. We will define, for each pair of integers $n < m$,
an exact functor $\Xi_{n,m}^\alpha: \mathcal{Q}_{G_{n}^{\alpha}}
(Q_n^{\alpha}) \to \mathcal{Q}_{G_{m}^{\alpha}}(Q_m^\alpha)$.
We first fix certain notations and data. For all $n$ we put
$$\mathcal{E}^{\alpha,s}_{n,n+1}=\bigoplus_{s'}
\Hom(\mathcal{L}_s(n),\mathcal{L}_{s'}(n+1))\otimes k^{d_{s'}(n+1,\alpha)}$$
and we set $\mathcal{E}_{n,n+1}^\alpha=\bigoplus_s
\mathcal{E}^{\alpha,s}_{n,n+1} \otimes \mathcal{L}_s(n)$. There is a
canonical (evaluation) surjective map $can_{n}:
\mathcal{E}^{\alpha}_{n,n+1} \tto \mathcal{E}^\a_{n+1}$.
We also fix a collection of surjective maps $u_{n}^s:
\mathcal{E}^{\alpha,s}_{n,n+1} \tto \mathcal{E}^{\alpha,s}_n$ and set
$u_{n}=\bigoplus_s u^s_{n}:\;\mathcal{E}_{n,n+1}^\a \tto \mathcal{E}^\a_n$.

\vspace{.1in}

Consider the following
open subvariety of $Q_n^\alpha$
$$Q_{n,\geq n+1}^\alpha=\{(\phi: \mathcal{E}_{n}^{\alpha} \tto
\mathcal{F}) \in Q_n^\alpha|
\; \mathcal{F}\;\text{is\;generated\;by\;} \{\mathcal{L}_s(n+1)\}_{s \in S}\}$$
and let us denote by $j: Q_{n,\geq n+1}^{\alpha} \hookrightarrow
Q_n^\alpha$ the
embedding. We also define an open subvariety of
$\mathrm{Hilb}_{\mathcal{E}_{n,n+1}^\alpha,\alpha}$ by
\begin{equation*}
\begin{split}
R_{n,n+1}^\alpha=\{\phi: \mathcal{E}_{n,n+1}^\alpha \tto \mathcal{F}\;|&\;
\phi_{*}: \mathcal{E}^{\alpha,s}_{n,n+1}
\tto
\Hom(\mathcal{L}_s(n),\mathcal{F})\;\text{is\;surjective}\\
&\text{for\;all\;}s \in S\;\text{and}\; \mathcal{F}\;\text{is\;generated\;by}\;
\{\mathcal{L}_s(n+1)\}_s\}.
\end{split}
\end{equation*}

\vspace{.1in}

Composing with $can_n$ induces an embedding
$\iota_{n,n+1}:
Q_{n+1}^{\alpha} \hookrightarrow
\mathrm{Hilb}_{\mathcal{E}_{n,n+1}^\alpha,\alpha}.$
Similarly, composing with $u_{n}$ gives an embedding
$v_{n,n+1}:\;Q^\alpha_{n,\geq n+1} \hookrightarrow R_{n,n+1}$.

\vspace{.1in}

Set $G_{n,n+1}^\alpha=\prod_s GL(\mathcal{E}^{\alpha,s}_{n,n+1})$. Let
$P_s \subset GL(\mathcal{E}^{\alpha,s}_{n,n+1})$
be the stabilizer of $\mathrm{Ker}\;u_n^s$. Finally, put
$P=\prod_s P_s
\subset G_{n,n+1}^\alpha \subset \mathrm{Aut}(\mathcal{E}_{n,n+1}^\alpha)$.
The group $\mathrm{Aut}(\mathcal{E}_{n,n+1}^\alpha)$
(and thus $G_{n,n+1}^\alpha$) naturally acts on
$\mathrm{Hilb}_{\mathcal{E}_{n,n+1}^\alpha,\alpha}$.
The group $P_s$ acts on $Q_n^\alpha$
via the quotient induced by $u_n^s$
$$P_s \tto GL(d_s(n,\alpha))
\subset G_{n}^\a.$$
By construction, there is a natural $G_{n,n+1}^\alpha$-equivariant isomorphism
\begin{equation}\label{E:4}
G_{n,n+1}^\alpha \underset{P}{\times} Q_{n,\geq n+1}^\alpha \simeq
R_{n,n+1}^\alpha.
\end{equation}

\vspace{.1in}

Note that there is a canonical embedding $GL(d_{s'}(n+1,\alpha))
\hookrightarrow
GL(\mathcal{E}^{\alpha,s}_{n,n+1})$, giving rise to an embedding
$$G_{n+1}^\alpha = \prod_{s'}GL(d_{s'}(n+1,\alpha))
\hookrightarrow \prod_s GL(\mathcal{E}^{\alpha,s}_{n,n+1}) =G_{n,n+1}^\alpha.$$

\begin{lem}\label{L:5} The map $\iota_{n,n+1}$ induces a canonical
$G_{n,n+1}^\alpha$-equivariant
surjective morphism
$$\theta_{n,n+1}^\alpha: G_{n,n+1}^\alpha \underset{G_{n+1}^\alpha}{\times}
Q_{n+1}^\alpha \tto R_{n,n+1}^\alpha .$$\end{lem}

\noindent
\textit{Proof.} Let us first show that the image of $\iota_{n,n+1}$
belongs to $R_{n,n+1}^\alpha$.

\vspace{.05in}

\begin{claim} If $\mathcal{F}$ is generated by
$\{\mathcal{L}_s(m)\}_{s \in S}$ with $m >n$ then the natural map
$$a:\bigoplus_{s'}\Hom(\mathcal{L}_s(n), \mathcal{L}_{s'}(m)) \otimes
\Hom(\mathcal{L}_{s'}(m),\mathcal{F})
\to \Hom(\mathcal{L}_s(n),\mathcal{F})$$ is surjective for every $s \in S$.
\end{claim}
\noindent
\textit{Proof of claim.}  We argue by induction on the rank.
Since $\mathcal{T}(1) \simeq \mathcal{T}$ for any torsion sheaf,
it is enough to prove the claim for vector bundles.
The statement is clear for line bundles. Let $\mathcal{F}_1 \subset
\mathcal{F}$ be a subsheaf isomorphic to $\mathcal{L}_{\a_\star}(m)$.
Applying $\Hom(\mathcal{L}_s(n), \cdot)$ to the exact
sequence $0 \to \mathcal{F}_1 \to \mathcal{F} \to
\mathcal{F}/\mathcal{F}_1 \to 0$
yields a sequence
$$0 \to \Hom(\mathcal{L}_s(n),\mathcal{F}_1) \to
\Hom(\mathcal{L}_s(n),\mathcal{F})
\to \Hom(\mathcal{L}_s(n), \mathcal{F}/\mathcal{F}_1) \to
\Ext(\mathcal{L}_s(n),
\mathcal{F}_1).$$
Using Serre duality, we have
\begin{equation*}
\begin{split}
\mathrm{dim}\;\Ext(\mathcal{L}_s(n),\mathcal{L}_{\a_\star}(m))
&=\dim\;\Hom(\mathcal{L}_{a_\star}(m), \mathcal{L}_s(n) \otimes \omega_X) \\
&\leq
\dim\;\Hom(\mathcal{L}_s(n),\mathcal{L}_s(n) \otimes \omega_X)=H^0(\omega_X)=0.
\end{split}
\end{equation*}
Now, it is clear that $\mathrm{Im}\;a \supset
\Hom(\mathcal{L}_s(n),\mathcal{F}_1)$.
Hence it is enough to show that $\mathrm{Im}\;a$ surjects onto
$\Hom(\mathcal{L}_s(n),\mathcal{F}/\mathcal{F}_1)$.
We have a commutative diagram of canonical maps
$$
\xymatrix{
\bigoplus_{s'} \Hom(\mathcal{L}_s(n),\mathcal{L}_{s'}(m)) \otimes \Hom(
\mathcal{L}_{s'}(m), \mathcal{F}) \ar@{->>}[d]^{c} \ar[r]^-{a} &
\Hom(\mathcal{L}_s(n),\mathcal{F}) \ar[d]^{c'}\\
\bigoplus_{s'} \Hom(\mathcal{L}_s(n),\mathcal{L}_{s'}(m)) \otimes \Hom(
\mathcal{L}_{s'}(m), \mathcal{F}/\mathcal{F}_1) \ar@{->>}[r]^-{a'} &
\Hom(\mathcal{L}_s(n),\mathcal{F}/\mathcal{F}_1)
}
$$
The map $c$ is surjective by the computation above and, by
induction hypothesis, so is $a'$. Thus $c' \circ a$ is onto as desired.\qed

\vspace{.1in}

Next, fix $(\phi: \mathcal{E}_{n+1}^{\alpha} \tto \mathcal{F}) \in
Q_{n+1}^{\alpha}$
and consider the following commutative diagram,
$$
\xymatrix{
\bigoplus_{s'}\Hom(\mathcal{L}_s(n),\mathcal{L}_{s'}(n+1)) \otimes \Hom(
\mathcal{L}_{s'}(n+1), \mathcal{F})
\ar@{->>}[r]^-{a} &
\Hom(\mathcal{L}_s(n),\mathcal{F})\\
\bigoplus_{s'}\Hom(\mathcal{L}_s(n),\mathcal{L}_{s'}(n+1)) \otimes
k^{d_{s'}(n+1,\alpha)}
\ar[u]^{1 \otimes \phi_*}_\sim
\ar[ur]^{\iota(\phi)_{*}}}
$$
By the claim $a$ is surjective and by hypothesis $\phi_{*}$ is an isomorphism.
Hence $\iota(\phi)_{*}$ is also surjective, and $\iota_{n,n+1}$ maps
$Q_{n+1}^\alpha$ into $R_{n,n+1}^\alpha$. Observe that $\iota_{n,n+1}$ is
$G_{n+1}^\alpha$-equivariant so that $\theta_{n,n+1}^\alpha$ is
well-defined and
$G_{n,n+1}^\alpha$-equivariant. Finally, $G_{n,n+1}^\alpha$-orbits
are parametrized
on both sides by isomorphism classes of coherent sheaves of class
$\alpha$, generated
by $\{\mathcal{L}_s(n+1)\}_{s \in S}$. Hence $\theta_{n,n+1}^\alpha$
is onto.\qed

\vspace{.1in}

\begin{lem}\label{L:6} The map $\theta_{n,n+1}^\alpha$ is smooth with
connected fibers of dimension
$$\mathrm{dim} (P/G_n^\alpha)=\sum_s
c_s(n,n+1,\alpha)(c_s(n,n+1,\alpha)-d_s(n,\alpha)).$$\end{lem}

\noindent
\textit{Proof.} Let $\mathcal{F}$ be a coherent sheaf of class $\alpha$,
generated by $\{\mathcal{L}_s(n+1)\}_{s \in S}$.
By (\ref{E:4}), the dimension of the orbit corresponding to
$\mathcal{F}$ in $R_{n,n+1}^\alpha$
is equal to $(\mathrm{dim}\;G_{n,n+1}^\alpha-\mathrm{dim}\;P) +
(\mathrm{dim}\;G_{n}^\alpha
-\mathrm{dim}\;
\mathrm{Aut}(\mathcal{F}))$. On the other hand, the dimension of the
corresponding orbit in
$G_{n,n+1}^\alpha \underset{G_{n+1}^\alpha}{\times} Q_{n+1}^\alpha$
is equal to $(\mathrm{dim}\;G_{n,n+1}^\alpha
-\mathrm{dim}\; G_{n+1}^\alpha) +
(\mathrm{dim}\;G_{n+1}^\alpha-\mathrm{dim\;Aut}(\mathcal{F}))$. This
gives
the dimension of the fiber. Smoothness and connectedness are clear.\qed

\vspace{.1in}

\paragraph{\textbf{Definition.}} Introduce the functor
\begin{align*}
\widetilde{\Xi}_{n,n+1}^\alpha:\;
\mathcal{Q}_{G_n^\alpha}(Q_n^\alpha) &\to
\mathcal{Q}_{G_{n+1}^\alpha}(Q_{n+1}^\alpha)\\
F &\mapsto res_{G_{n,n+1}^\alpha}^{G_{n+1}^\alpha} \circ
(\theta_{n,n+1}^{\alpha})^* \circ
ind_P^{G_{n,n+1}^\alpha}\circ j^*(F)
\end{align*}
and put $\Xi_{n,n+1}^\alpha=\widetilde{\Xi}_{n,n+1}^\alpha
[\mathrm{dim}\;G_{n+1}^\alpha-\mathrm{dim}\;G_n^\alpha]$.
We extend this definition to arbitrary $n<m$ by setting
$$\widetilde{\Xi}_{n,m}:= \widetilde{\Xi}_{m-1,m} \circ \cdots \circ
\widetilde{\Xi}_{n,n+1},\qquad
{\Xi}_{n,m}:= {\Xi}_{m-1,m} \circ \cdots \circ {\Xi}_{n,n+1}.$$

\vspace{.1in}

Note that $\widetilde{\Xi}_{n,n+1}^\alpha(F)=(\iota_{n,n+1})^*
\circ ind_P^{G_{n,n+1}} \circ j^*(F)$.
 From the properties of $ind_P^{G_{n,n+1}}$ one deduces that
$G=\widetilde{\Xi}_{n,n+1}(F)$ is characterized, up to isomorphism, by
the property that there exists an element $H \in
\mathcal{Q}_{G_{n,n+1}^\alpha}(R_{n,n+1}^\a)$ such that $v_{n,n+1}^*H=j^*F$
and $\iota_{n,n+1}^*H=G$.

Using Lemma~\ref{L:6} and the fact that $j$ is an open embedding, we
see that $\Xi_{n,n+1}^\alpha$ is exact, preserves perversity
and commutes with Verdier duality.

\vspace{.1in}

The functors $\widetilde{\Xi}_{n,m}$ can be given a direct definition
as follows~: let us put
$$\mathcal{E}^{\alpha,s}_{n,m}=\bigoplus_{\underline{s}}
\big(\bigotimes_{j=n}^{m-1} \mathrm{Hom}(\mathcal{L}_{s_j}(j),
\mathcal{L}_{s_{j+1}}(j+1))\big) \otimes k^{d_{s_m}(m,\alpha)},$$
$$\mathcal{E}^\alpha_{n,m}=\bigoplus_{s} \mathcal{L}_s(n) \otimes
\mathcal{E}^{\alpha,s}_{n,m},$$
where $\underline{s}=\{(s_j)_{j=n+1}^m\;|s_j \in S\}$. We have a
canonical (evaluation) map $can_{n,m}: \mathcal{E}^\alpha_{n,m} \tto
\mathcal{E}^\alpha_m$
and a composition $u_{n,m}:=u_n \circ \cdots \circ u_{m-1}:
\mathcal{E}^\alpha_{n,m}\tto \mathcal{E}^\a_n$.
Putting
\begin{equation*}
\begin{split}
R^\alpha_{n,m}=\{\phi: \mathcal{E}^\alpha_{n,m} \tto \mathcal{F}\;|\;&
\phi_*^s: \mathcal{E}^{\alpha,s}_{n,m}
\tto
\Hom(\mathcal{L}_s(n),\mathcal{F})\;\text{is\;surjective}\\
&\text{for\;all\;}s \in S\;\text{and}\; \mathcal{F}\;\text{is\;generated\;by}\;
\{\mathcal{L}_s(m)\}_s\}.
\end{split}
\end{equation*}
we thus obtain embeddings $\iota_{n,m}: Q_m^\alpha \to R_{n,m}^\alpha$
and $V_{n,m}: Q^\alpha_{n,\geq m} \to R_{n,m}^\alpha$. Finally, the group
$G_{n,m}^\alpha=\prod_s \mathrm{Aut}(\mathcal{E}^{\a,s}_{n,m})$
naturally acts on $R_{n,m}^\alpha$. If $F \in
\mathcal{Q}_{G_n^\alpha}(Q_n^\alpha)$ then
$G=\widetilde{\Xi}^\alpha_{n,m}(F)$ is characterized by the condition
that there exists $H \in \mathcal{Q}_{G_{n,m}^\alpha}(R_{n,m}^\a)$
such that $v_{n,m}^*H=j^*F$
and $\iota_{n,m}^*H=G$.

\vspace{.1in}

By construction, we have

\begin{lem}\label{L:7} Let $\alpha \in K^+(X)$ and $n<m<l$. There is
a canonical natural transformation
$ \Xi_{m,l}^\alpha \circ \Xi_{n,m}^\alpha \stackrel{\sim}{\to}
\Xi_{n,l}^\alpha$.
\end{lem}

\vspace{.1in}

\paragraph{\textbf{3.3. The category} $\mathbb{Q}^\alpha$.}
We will now define a triangulated category
$\mathbb{Q}^\alpha$ as a projective limit of the system
$(\mathcal{Q}_{G_n^\alpha}(Q_n^\alpha),
\Xi_{n,m})_{n,m \in \Z}$. Let $\mathbb{Q}^\alpha$ be the additive
category with \\
-$\mathcal{O}bj(\mathbb{Q}^\alpha)$ is the set of
collections $(F_n, r_{n,m})_{n <m \in \Z}$ where $F_n \in
\mathcal{Q}_{G^\alpha_n}(Q_n^\alpha)$ and
$r_{n,m}: \Xi_{n,m}(F_n) \stackrel{\sim}{\to} F_m$ satisfy the conditions
$$r_{n,l}=r_{m,l} \circ \Xi_{m,l}(r_{n,m})\qquad \text{for\;every\;}n<m<l.$$
\noindent
-For any two objects $\mathbb{F}=(F_n,r_{n,m}),
\mathbb{F}'=(F'_n,r'_{n,m})$, we have
$$\Hom_{\mathbb{Q}^\alpha}(\mathbb{F},\mathbb{F}')=\{(h_n) \in
\prod_n \Hom(F_n,F'_n)\;| h_m=r'_{n,m}
\circ \Xi_{n,m}(h_n) \circ r_{n,m}^{-1}\;\; \text{if\;}n<m\}.$$

\vspace{.05in}

The collection of translation functors $F \mapsto F[1]$ of
$\mathcal{Q}_{G_n^\alpha}(Q_n^\alpha))$
give rise to
an automorphism $\mathbb{F} \mapsto \mathbb{F}[1]$ of $\mathbb{Q}^\alpha$.
Let us call a triangle
$\mathbb{F} \to \mathbb{G} \to \mathbb{H} \to \mathbb{F}[1]$
distinguished if all the
corresponding triangles $F_n \to G_n \to H_n \to F_n[1]$ are.

\begin{lem}\label{L:8} The category $\mathbb{Q}^\alpha$, equipped
with the automorphism
$\mathbf{T}$ and the above collection of distinguished triangles, is
a triangulated category.\end{lem}

\vspace{.1in}

\paragraph{\textbf{Remarks.}} i) Assume that $\alpha$ is the class of
a torsion sheaf. Then it
it easy to see that $Q_n^\alpha \simeq Q_m^\alpha$ for any two $n,m
\in \Z$. Hence in this case
$\mathbb{Q}^\alpha \simeq \mathcal{Q}_{G_n^\alpha}(Q_n^\alpha)$.\\
ii) Now assume that $\alpha$ is the class of a sheaf of rank at least
one. Then from Section~2,
one deduces that $Q_n^\alpha$ is empty for $n$ big enough.

\vspace{.1in}

We will still call objects of $\mathbb{Q}^\alpha$ semisimple
complexes. Note that since Verdier duality commutes with the
transfer functors $\Xi_{n,m}^\a$, it gives rise to an involutive
functor on $\mathbb{Q}^\a$, which we still call Verdier duality
and denote again by $D$. The category $\mathbb{Q}^\a$ is closed
under certain infinite sums~: let us call a countable collection
$\{\mathbb{F}^u\}_{u \in U}$ of objects of $\mathbb{Q}^\alpha$
\textit{admissible} if for any $n \in \Z$ the set of $u \in U$ for
which $F_n^u \neq 0$ is finite. It is clear that if
$\{\mathbb{F}^u\}_{u \in U}$ is admissible then the direct sum
$\bigoplus_u \mathbb{F}^u$ is a well-defined object of
$\mathbb{Q}^\alpha$. Let us call an object $\mathbb{F}$
\textit{simple} if $F_n$ is a simple perverse sheaf for all $n$.
Note that if $\mathbb{F}, \mathbb{H}$ are simple and $F_m \simeq
H_m \neq 0$ for any $m \in \Z$ then $\mathbb{F}$ and $\mathbb{H}$
are isomorphic. Indeed, for any $n<m$, $\Xi_{n,m}$ induces an
equivalence $\mathcal{Q}_{G_{n}^\alpha}(Q_{n,\geq m})
\stackrel{\sim}{\to} \mathcal{Q}_{G_m^\alpha}(Q_m^\alpha)$ thus
$j^*_{n,m}(F_n) \simeq j^*_{n,m}(G_n)$. But $F_n$ and $H_n$ being
simple perverse sheaves, they are determined up to isomorphism by
their restriction to the open set $Q^\alpha_{n,\geq m}$.

\vspace{.1in}

\begin{cor}\label{C:01} For any object $\mathbb{H} \in
\mathbb{Q}^\alpha$ there exists an admissible collection
of simple objects $\{\mathbb{F}^u\}_{u \in U}$ and integers $d_u$
such that $\mathbb{H}\simeq\bigoplus_u
\mathbb{F}^u[d_u]$.\end{cor}

\vspace{.2in}

\section{Induction and restriction functors}

\vspace{.1in}

\paragraph{\textbf{4.1. Induction functor.}} As in \cite{L1}, we
consider, for all $n \in \Z$ and
$\alpha, \beta \in K^+(X)$ a diagram~:
\begin{equation}\label{E:diagind}
\xymatrix{Q_n^\beta \times Q_n^\alpha& E' \ar[l]_-{p_1} \ar[r]^-{p_2}
& E'' \ar[r]^-{p_3} &
Q_n^{\alpha +\beta}}
\end{equation}
where we used the following notations:\\
-$E'$ is the variety of tuples $(\phi, (V_s, a_s, b_s)_{s \in S})$ with
\begin{enumerate}
\item[-] $(\phi: \mathcal{E}^{\alpha+\beta}_n \tto \mathcal{F}) \in
Q_n^{\alpha+\beta}$,
\item[-] $V_s \subset k^{d_s(n,\alpha+\beta)}$ is a subspace of
dimension $d_s(n,\alpha)$ such
that
$$[\phi(\bigoplus_s V_s \otimes \mathcal{L}_s(n))]=\alpha,$$
\item[-] $a_s: V_s \stackrel{\sim}{\to} k^{d_s(n,\alpha)}$, $b_s:
k^{d_s(n,\alpha+\beta)}/V_s
\stackrel{\sim}{\to} k^{d_s(n,\beta)}$,
\end{enumerate}
-$E''$ is the variety of tuples $(\phi,(V_s)_{s \in S})$ as above,\\
-$p_1(\phi,(V_s,a_s,b_s)_{s \in S})=
( b_*\phi_{|\mathcal{E}^{\alpha+\beta}_n/\underline{V}},
a_* \phi_{|\underline{V}})$, where
$\underline{V}=\bigoplus_s
V_s \otimes \mathcal{L}_s(n)$,\\
-$p_2,p_3$ are the projections.

\vspace{.1in}
\begin{lem}\label{L:9} The map $p_1$ is smooth.
\end{lem}
\noindent
\textit{Proof.} Let $X$ be the variety of tuples $(V_s, a_s, b_s)_s$
as above (a
$G_n^{\alpha+\beta}$-homogeneous space). The map $p_1$ factors as
$$E' \stackrel{p_1'}{\to} X \times Q_n^\beta \times Q_n^\alpha
\stackrel{p_1''}{\to} Q_n^\b \times
Q_n^\a,$$
where $p_1'(\phi, (V_s,a_s,b_s)_s)=((V_s, a_s, b_s)_s,b_*
\phi_{|\mathcal{E}_n^{\alpha+\beta}/\underline{V}},a_*
\phi_{|\underline{V}})$ and $p_1''$ is
the projection. We claim that $p_1'$
is a vector bundle of rank
$\langle[\mathcal{E}^\beta_n]-\beta,\alpha\rangle$. Indeed,
the fiber of $p_1'$ over a point $((\phi_2: \mathcal{E}^\beta_n \tto \mathcal{F}_2),
(\phi_1: \mathcal{E}_n^\alpha \tto
\mathcal{F}_1))$ is canonically
isomorphic to $\Hom(\mathrm{Ker}(\phi_2),
\mathcal{F}_1)$. Let us set $K_i=\mathrm{Ker}(\phi_i)$ for $i=1,2$.
 From the sequence $0 \to K_2
\to \mathcal{E}^\beta_n \to \mathcal{F}_2 \to 0$ we have
$$\Ext(\mathcal{F}_2, \mathcal{E}^\beta_n) \to
\Ext(\mathcal{E}^\beta_n,\mathcal{E}^\beta_n)
\to
\Ext(K_2, \mathcal{E}^\beta_n) \to 0.$$
But since $\Ext(\mathcal{L}_s(n), \mathcal{L}_{s'}(n))=0$ for all $s,s'$,
\begin{equation}\label{E:9}
\Ext(\mathcal{E}^\beta_n, \mathcal{E}^\beta_n)=\Ext(K_2,
\mathcal{E}^\beta_n)=0.
\end{equation}
Similarly, applying
$\Hom(K_2, \cdot)$ to $0 \to K_1 \to \mathcal{E}^\alpha_n \to
\mathcal{F}_1 \to 0$ yields a sequence
$\Ext(K_2,K_1) \to \Ext(K_2, \mathcal{E}^\alpha_n) \to \Ext(K_2,
\mathcal{F}_1)\to 0$. Using
(\ref{E:9}) we now obtain $\Ext(K_2, \mathcal{F}_1)=0$ and hence
$\mathrm{dim}\;\Hom(K_2, \mathcal{F}_1)
=\langle[K_2],[\mathcal{F}_1]\rangle=\langle[\mathcal{E}^{\beta}_n]-\beta,\alpha\rangle$.\qed

\vspace{.1in}

From the above proof it follows that $E'$ is smooth.
Note that $p_2$ is a principal $G_n^\b \times G_n^\a$-bundle,
and hence $E''$ is also
smooth. Unfortunately, $p_3$ is not
proper in general, as one can readily see when $X=\mathbb{P}^1$ and
$G=Id$. However, the following weaker property will be sufficient for us.

\vspace{.1in}

\begin{lem}\label{L:10} Fix $m \in \Z$ such that $n \leq
q(\alpha,\alpha+\beta,m)$. Then $p_3: p_3^{-1}
(Q_{n,\geq m}^{\alpha+\beta}) \to Q_{n,\geq m}^{\alpha+\beta}$ is
proper.\end{lem}
\noindent
\textit{Proof.}  Let $(\phi: \mathcal{E}^{\alpha+\beta}_n \tto \mathcal{F}) \in
Q_{n, \geq m}^{\alpha+\beta}$. Since $n \leq
q(\alpha,\alpha+\beta,m)$, Lemma~\ref{L:2} implies that
any subsheaf $\mathcal{G} \subset \mathcal{F}$ of class $\alpha$ is
generated by $\{\mathcal{L}_s(n)\}$
and satisfies $\Ext(\mathcal{L}_s(n), \mathcal{G})=0$ for all $s \in
S$. Therefore, we have a canonical
commutative square
$$\xymatrix{
\bigoplus_s \Hom(\mathcal{L}_s(n), \mathcal{G}) \otimes
\mathcal{L}_s(n) \ar@{->>}[r]
\ar@{->}[d] & \mathcal{G} \ar@{^{(}->}[d]\\
\bigoplus_s \Hom(\mathcal{L}_s(n), \mathcal{F}) \otimes
\mathcal{L}_s(n) \ar@{->>}[r]
& \mathcal{F}}
$$
and composing with the isomorphism $\phi_*
:\Hom(\mathcal{L}_s(n),\mathcal{F}) \stackrel{\sim}{\to}
k^{d_s(n,\alpha+\beta)}$ yields a collection of subspaces $\{V_s
\subset k^{d_s(n,\alpha+\beta)}\}$
such that $(\phi, (V_s)) \in p_3^{-1}(\phi)$. Conversely, any such
collection of subspaces $\{V_s\}$
arises in this way. In other terms, $p_3^{-1}(Q_{n, \geq
m}^{\alpha+\beta})$ represents the functor
$\underline{\mathrm{Hilb}}^0_{\mathcal{E}_n^{\alpha+\beta},\alpha+\beta,\beta}$
from the category of smooth
schemes over $k$ to the category of sets which assigns to $\Sigma$
the set of all $ \mathcal{G}_{\alpha+\beta} \subset \mathcal{G}_{\beta} \subset
\mathcal{E}^{\alpha+\beta}_n \boxtimes \mathcal{O}_{\Sigma}$ where
$\mathcal{G}_\beta,
\mathcal{G}_{\alpha+\beta}$ are $\Sigma$-flat, $G$-equivariant coherent sheaves
such that
\begin{enumerate}
\item[-]
$[\mathcal{E}^{\alpha+\beta}_n/\mathcal{G}_{\beta|\sigma}]=\beta$
(resp.
$[\mathcal{E}^{\alpha+\beta}_n/\mathcal{G}_{\alpha+\beta|\sigma}]=\alpha+\beta$
and is generated by $\{\mathcal{L}_s(m)\}$) for any closed point
$\sigma \in \Sigma$,
\item[-] The projection map $\phi: \mathcal{E}_{n}^{\alpha+\beta} \tto
\mathcal{E}^{\alpha+\beta}_n/\mathcal{G}_{\alpha+\beta}$ induces isomorphisms
$$\phi_{*|\sigma}:\; k^{d_s(n,\alpha+\beta)} \stackrel{\sim}{\to}
\Hom(\mathcal{L}_s(n),\mathcal{E}_{n}^{\alpha+\beta}/\mathcal{G}_{\alpha+\beta|\sigma})$$
for any closed point $\sigma.$
\end{enumerate}
The map $p_3$ is induced by the natural forgetful map $( \mathcal{G}_{\alpha},
\mathcal{G}_{\alpha+\beta}) \mapsto \mathcal{G}_{\alpha + \beta}$.
Hence the fiber at
a point $(\phi: \mathcal{E}^{\alpha+\beta}_n \tto \mathcal{F})$
represents the usual
Quot functor $\underline{\mathrm{Hilb}}_{\mathcal{F},\beta}$, and
from Proposition~\ref{P:1}
we deduce that $p_3$ is proper.\qed

\vspace{.1in}

\paragraph{}If $n \leq q(\alpha,\alpha+\beta,m)$ we set
$$\widetilde{\mathrm{Ind}}_{n,m}^{\b,\a}=\widetilde{\Xi}^{\alpha+\beta}_{n,m}
\circ p_{3!} p_{2\flat}p_1^*:
\mathcal{Q}_{G_n^\b}(Q_n^\b) \times
\mathcal{Q}_{G_{n}^\a}(Q_n^\a) \to
\mathcal{D}(Q_m^{\alpha+\beta})$$
and
$\mathrm{Ind}^{\b,\a}_{n,m}=\widetilde{\mathrm{Ind}}^{\b,\a}_{n,m}[\mathrm{dim}
(G_m^{\alpha+\beta}) - \mathrm{dim}(G_n^\b \times G_n^\a)
-\langle\beta,\alpha\rangle]$. With this
notation, $\mathrm{Ind}^{\b,\a}_{n,m}$ commutes with Verdier
duality. Note that from
Lemma~\ref{L:10} and the (equivariant version of the) Decomposition
theorem \cite{BBD} it follows that
$\mathrm{Ind}_{n,m}^{\b,\a}$ actually takes values in
$\mathcal{Q}_{G_m^{\alpha+\beta}}
(Q_m^{\alpha+\beta})$ (recall that all perverse sheaves considered
are of geometric origin).

\vspace{.1in}

\begin{lem}\label{L:11} i) Let $n \leq q(\alpha,\alpha+\beta,m) \leq
m \leq l$. There is a canonical natural
transformation
$$\Xi^{\alpha+\beta}_{m,l} \circ \mathrm{Ind}_{n,m}^{\b,\a}
\simeq \mathrm{Ind}^{\b,\a}_{n,l}.$$
ii) Let $n \leq m \leq q(\alpha,\alpha+\beta,l) \leq l$. There is a
canonical natural transformation
$$\mathrm{Ind}_{m,l}^{\b,\a} \circ (\Xi_{n,m}^\b \times
\Xi_{n,m}^\a) \simeq
\mathrm{Ind}_{n,l}^{\b,\a}.$$
\end{lem}
\noindent
\textit{Proof.} It is enough to prove both formulas with
$\widetilde{\Xi}$ and $\widetilde{\mathrm{Ind}}$.
Statement i) follows directly from Lemma~\ref{L:7} and the
definitions. Now let $n,m,l$ be as in ii). From
Lemma~\ref{L:7}, it is enough to show the commutativity of the
following diagrams of functors
\begin{equation}\label{E:10}
\xymatrix{
\mathcal{Q}_{G_n^\b}(Q_n^\b) \times
\mathcal{Q}_{G_n^\a}(Q_n^\a) \ar[r]^-{p_1^*}
\ar[d]^{\widetilde{\Xi}_{n,m}} & \mathcal{Q}_{G_n^{\alpha+\beta}}(E'_n)
\ar[r]^-{p_{2\flat}} & \mathcal{Q}_{G_n^{\alpha+\beta}}(E''_n)
\ar[r]^-{j_{n,l}^*p_{3!}}&
\mathcal{Q}_{G_n^{\alpha+\beta}}(Q_{n, \geq l}^{\alpha+\beta})
\ar[d]^{\widetilde{\Xi}_{n,m}}\\
\mathcal{Q}_{G_m^\b}(Q_m^\b) \times
\mathcal{Q}_{G_m^\a}(Q_m^\a) \ar[r]^-{p_1^*}&
  \mathcal{Q}_{G_m^{\alpha+\beta}}(E'_m)
\ar[r]^-{p_{2\flat}} & \mathcal{Q}_{G_m^{\alpha+\beta}}(E''_m)
\ar[r]^-{j_{m,l}^*p_{3!}}&
\mathcal{Q}_{G_m^{\alpha+\beta}}(Q_{m, \geq l}^{\alpha+\beta})}
\end{equation}
Since $n,m \leq q(\alpha, \alpha+\beta,l)$, we have, as in the proof
of Lemma~\ref{L:10},
canonical identifications
$$E''_n \simeq \mathrm{Hilb}^0_{\mathcal{E}^{\alpha+\beta}_n,
\alpha+\beta, \beta}, \qquad
E''_m \simeq \mathrm{Hilb}^0_{\mathcal{E}^{\alpha+\beta}_m,
\alpha+\beta, \beta}.$$
This allows one to define transfer functors $\Xi_{n,m}:
\mathcal{Q}_{G_n^{\alpha+\beta}}(E''_n)
\to \mathcal{Q}_{G_m^{\alpha+\beta}}(E''_m)$ and $\Xi_{n,m}:
\mathcal{Q}_{G_n^{\alpha+\beta}}(E'_n)
\to \mathcal{Q}_{G_m^{\alpha+\beta}}(E'_m)$ in the same manner as in
Section~3, and therefore
complete the diagram (\ref{E:10}) with two middle vertical arrows.
The commutativity
of each ensuing square follows from standard base change arguments.
We leave the details to the
reader.\qed

\vspace{.1in}

\paragraph{\textbf{Definition.}} Let $\mathbb{F} \in
\mathbb{Q}^\alpha, \mathbb{H} \in \mathbb{Q}^\beta$.
By Lemma~\ref{L:11} ii),
for each fixed $m \in \Z$, the complexes
$\mathrm{Ind}^{\b,\a}_{n,m}(H_n\boxtimes
F_n)$ are all \textit{canonically} isomorphic for $n \leq
q(\alpha,\alpha+\beta,m)$.
Furthermore,
by Lemma~\ref{L:11} i),
for each $l>m$ the complexes $\Xi_{m,l}
(\mathrm{Ind}_{n,m}^{\b,\a}(H_n \boxtimes
F_n))$ and $\mathrm{Ind}_{n,l}^{\b,\a}(H_n \boxtimes F_n)$ are
\textit{canonically} isomorphic.
The collection of complexes $\{\mathrm{Ind}_{n,m}^{\b,\a}(H_n
\boxtimes F_n), \; n \leq
q(\alpha,\alpha+\beta,m)\}$ thus gives rise to an object of
$\mathbb{Q}^{\alpha+\beta}$ which is unique up to isomorphism. We
denote this functor by
$$\mathrm{Ind}^{\b,\a}:\; \mathbb{Q}^\b \times \mathbb{Q}^\a \to
\mathbb{Q}^{\alpha+\beta}$$
and call it the induction functor.

\vspace{.1in}

\begin{lem}[associativity of induction]\label{L:12}
  For each triple $\alpha, \beta, \gamma \in K^+(X)$ there are
(canonical) natural transformations $\mathrm{Ind}^{
\beta+\gamma,\a} \circ \mathrm{Ind}^{\gamma,\b}
\simeq \mathrm{Ind}^{,\gamma,\alpha+\beta} \circ
\mathrm{Ind}^{\b,\a}$.\end{lem}
\noindent
\textit{Proof.} It is enough to show that, for any $n\ll m \ll l$
we have a canonical natural transformation
\begin{equation}\label{E:11}
\widetilde{\mathrm{Ind}}_{m,l}^{\gamma,\alpha+\beta} \circ
(\widetilde{\Xi}_{n,m}^\gamma \times \widetilde{\mathrm{Ind}}_{n,m}^{\b,\a}) \simeq
\widetilde{\mathrm{Ind}}_{m,l}^{\beta+\gamma,\alpha} \circ
(\widetilde{\mathrm{Ind}}_{n,m}^{\gamma,\beta} \times \widetilde{\Xi}_{n,m}^\alpha).
\end{equation}
Consider the functor
\begin{equation*}
\begin{split}
\mathrm{I}_n^{\beta,\a}:\;\mathcal{Q}_{G_n^\b}(Q_n^\b)
\times \mathcal{Q}_{G_n^\a}(Q_n^\a) &\to
\mathcal{D}(Q_n^{\alpha+\beta})\\
P &\mapsto p_{3!}p_{2\flat}p_1^*(P)
\end{split}
\end{equation*}
induced by the diagram (\ref{E:diagind}).
Using Lemma~\ref{L:11} we have
\begin{equation*}
\begin{split}
\widetilde{\mathrm{Ind}}_{m,l}^{\gamma,\alpha+\beta} \circ
(\widetilde{\Xi}_{n,m}^\gamma \times \widetilde{\mathrm{Ind}}_{n,m}^{\b,\a})=&
\widetilde{\mathrm{Ind}}_{m,l}^{\gamma,\alpha+\beta} \circ (\widetilde{\Xi}_{n,m}^\gamma
\times \widetilde{\Xi}_{n,m}^{\alpha+\beta}) \circ (Id \times \mathrm{I}^{\b,\a})\\
=&\widetilde{\mathrm{Ind}}_{n,l}^{\gamma,\alpha+\beta} \circ
(Id \times \mathrm{I}^{\b, \a})\\
=& \widetilde{\Xi}_{n,l}^{\alpha+\beta+\gamma} \circ
\mathrm{I}_n^{\gamma,\alpha+\beta} \circ
(Id \times \mathrm{I}^{\b, \a}),
\end{split}
\end{equation*}
and a similar expression for the right-hand side of (\ref{E:11}).
Hence it is enough to prove that
$\mathrm{I}_n^{\gamma,\alpha+\beta} \circ
(Id \times \mathrm{I}^{\b, \a}) \simeq
\mathrm{I}_n^{\beta+\gamma,\a} \circ
(\mathrm{I}^{\gamma,\b}\times Id)$.
Define the following triple analogue of diagram (\ref{E:diagind})
$$
\xymatrix{
Q_n^\gamma \times Q_n^\beta \times Q_n^\a & F' \ar[l]_-{q_1}
\ar[r]^-{q_2} & F'' \ar[r]^-{q_3} &
Q_n^{\alpha +\beta+\gamma}}$$
where\\
-$F'$ is the variety of tuples $(\phi, (W_s,V_s, a_s, b_s)_{s \in S})$ with
\begin{enumerate}
\item[-] $(\phi: \mathcal{E}^{\alpha+\beta+\gamma}_n \tto
\mathcal{F}) \in Q_n^{\alpha+\beta+\gamma}$,
\item[-] $W_s \subset V_s \subset k^{d_s(n,\alpha+\beta+\gamma)}$ is
a chain of subspaces of respective
dimensions $d_s(n,\alpha), d_s(n,\alpha+\beta)$ such
that
$$[\phi(\bigoplus_s W_s \otimes \mathcal{L}_s(n))]=\alpha,\qquad
[\phi(\bigoplus_s V_s \otimes \mathcal{L}_s(n))]=\alpha+\beta,$$
\item[-] $a_s: W_s \stackrel{\sim}{\to} k^{d_s(n,\alpha)}$, $b_s:
V_s/W_s \stackrel{\sim}{\to}
k^{d_s(n, \beta)}, c_s:
k^{d_s(n,\alpha+\beta+\gamma)}/V_s
\stackrel{\sim}{\to} k^{d_s(n,\gamma)}$,
\end{enumerate}
-$F''$ is the variety of tuples $(\phi,(W_s,V_s)_{s \in S})$ as above,\\
-$q_1(\phi,(W_s,V_s,a_s,b_s,c_s)_{s \in S})=(c_*
\phi_{|\mathcal{E}^{\alpha+\beta}_n/\underline{V}},
b_{*} \phi_{|\underline{V}/\underline{W}},a_{*} \phi_{|\underline{W}})$, where
$$\underline{W}=\bigoplus_s W_s \otimes \mathcal{L}_s(n),
\qquad\underline{V}=\bigoplus_s
V_s \otimes \mathcal{L}_s(n),$$ \\
-$q_2,q_3$ are the projections.

\vspace{.1in}

Standard arguments now show that
$$\mathrm{I}_n^{\gamma,\alpha+\beta} \circ
(Id \times \mathrm{I}^{\b, \a})=q_{3!}q_{2\flat}q_1^*=
\mathrm{I}_n^{\beta+\gamma,\a} \circ
(\mathrm{I}^{\gamma,\b}\times Id)$$
as wanted. \qed

\vspace{.1in}

Using Lemma~\ref{L:12}, we may now define an iterated induction map
$$\mathrm{Ind}^{\alpha_r, \ldots, \alpha_1}:\; \mathbb{Q}^{\alpha_r}
\times \cdots \times \Q^{\alpha_1}
\to \Q^{\alpha_1 + \cdots + \alpha_r}.$$

\vspace{.2in}

\paragraph{\textbf{4.2. Restriction functors.}} Let $n \in \Z$,
$\alpha, \beta \in K^+(X)$. For $s \in S$ we fix
a subspace $V_s \subset k^{d_s(n,\alpha+\beta)}$ of dimension $d_s(n,
\alpha)$ and we fix isomorphisms
$a_s: V_s \stackrel{\sim}{\to} k^{d_s(n,\alpha)}, b_s:
k^{d_s(n,\alpha+\beta)}/V_s \stackrel{\sim}{\to}
k^{d_s(n,\beta)}$.
We now consider the diagram
\begin{equation}\label{E:12}
\xymatrix{
Q_n^{\alpha+\beta} & F  \ar[l]_-{i} \ar[r]^-{\kappa} & Q_n^\b
\times Q_n^\a}
\end{equation}
where\\
- $F$ is the subvariety of $Q_n^{\alpha+\beta}$ consisting of $(\phi:
\mathcal{E}^{\alpha+\beta}_n
\tto \mathcal{F})$ such that
$$[\phi(\bigoplus_s V_s \otimes \mathcal{L}_s(n))]=\alpha$$
and $i: F \hookrightarrow Q_n^{\alpha+\beta}$ is the (closed) embedding,\\
- $\kappa(\phi)=( b_*
\phi_{|\mathcal{E}^{\alpha+\beta}_n/\underline{V}},a_{*} \phi_{|\underline{V}})$,
where we put $\underline{V}=\bigoplus_s V_s \otimes \mathcal{L}_s(n)$.

Observe that by the proof of Lemma~\ref{L:9}, $\kappa$ is a vector
bundle of rank $\langle [\mathcal{E}^\beta_n]-\beta,
\alpha \rangle$. We define a functor
$$\widetilde{\mathrm{Res}}^{\b, \a}_n=\kappa_{!}i^*:
\mathcal{Q}_{G_n^{\alpha+\beta}}(Q_n^{\alpha+\beta})
\to \mathcal{D}(Q_n^{\b} \times Q_n^\a)$$
and set
$\mathrm{Res}^{\b,\a}_n=\widetilde{\mathrm{Res}}_n^{\b,\a}
[-\langle \beta,\alpha\rangle]$.

\vspace{.1in}

\begin{lem}\label{L:13} Let $P \in
\mathcal{Q}_{G_n^{\ab}}(Q_n^{\ab})$ and assume that
$\mathrm{Res}_n^{\b,\a}(P) \in \mathcal{Q}_{G_n^\b}(Q_n^\b) \boxtimes
\mathcal{Q}_{G_n^\a}(Q_n^\a)$.
Then, for any $n <m \in \Z$ there is a canonical isomorphism of functors
  $(\Xi_{n,m}^\b \times \Xi_{n,m}^\a) \circ
\mathrm{Res}_n^{\b,\a}(P) \simeq \mathrm{Res}_m^{\b,\a} \circ
\Xi_{n,m}^{\ab} (P).$
\end{lem}
\noindent
\textit{Proof.} Given the definitions of $\Xi_{n,m}$ and
$\mathrm{Res}$, it is enough
to construct a natural transformation
\begin{equation*}
(\widetilde{\Xi}_{n,m}^{\b} \times \widetilde{\Xi}_{n,m}^{\a}) \circ
\widetilde{\mathrm{Res}}_n^{\b,\a} \simeq
\widetilde{\mathrm{Res}}_m^{\b,\a} \circ
\widetilde{\Xi}_{n,m}^{\ab}\big[2\sum_s
d_s(m,\a)d_s(m,\b)-d_s(n,\a)d_s(n,\b)\big].
\end{equation*}
Also we may assume that $m=n+1$.
Consider the inclusion diagram
\begin{equation}\label{E:14}
\xymatrix{
Q_n^{\ab}& F_n \ar[l]_-{i} \ar[r]^-{\kappa} & Q_n^\b \times Q_n^\a\\
Q_{n, \geq m}^{\ab} \ar[u]^-{j_1} & F_{n, \geq m} \ar[u]^-{j_2} \ar[l]_-{i'}
\ar[r]^-{\kappa'} & Q_{n,\geq m}^\b \times Q_{n,\geq m}^\a \ar[u]^-{j_3}}
\end{equation}

\begin{claim}\label{C:1}Both squares in (\ref{E:14}) are cartesian.
In particular, we have
\begin{equation}\label{E:15}
j_3^* \widetilde{\mathrm{Res}}_n^{\b,\a}=\kappa'_!{i'}^*j_1^*.
\end{equation}
\end{claim}
\noindent
\textit{Proof of claim.} This is clear for the first square. Now let
$\mathcal{F},\mathcal{G}$
be coherent sheaves generated by $\{\mathcal{L}_s(m)\}$ and let $0
\to \mathcal{F} \to
\mathcal{H} \to \mathcal{G} \to 0$ be an extension. By Lemma~2.1,
$\Ext(\mathcal{L}_s(m), \mathcal{F})
=\Ext(\mathcal{L}_s(m),\mathcal{G})=0$ for all $s \in S$. Then the
proof of Lemma~2.1
shows that $\mathcal{H}$ is also generated by $\{\mathcal{L}_s(m)\}$,
and hence that the
second square is also cartesian.\qed

Let $V_s^m \subset k^{d_s(m, \ab)}, a_s^m, b_s^m$ for $s \in S$ be as
in the definition of
$\mathrm{Res}_m^{\b,\a}$. Let $P_m^{\b,\a} \subset G_m^{\ab}$ be the
stabilizer of
$$\underline{V}^m=\bigoplus_{s'} V_{s'}^m \otimes \mathcal{L}_{s'}(m)$$
and let $P_{n,m}^{\b,\a}\subset G_{n,m}^{\ab}$ be the stabilizer of
$$\underline{V}_{n,m}=\bigoplus_{s,s'} \Hom(\mathcal{L}_s(n),
\mathcal{L}_{s'}(m) \otimes
V_{s'}^m) \otimes \mathcal{L}_s(n) \subset \mathcal{E}_{n,m}^{\ab}.$$
There is a natural affine fibration $P_{n,m}^{\b,\a} \tto G_{n,m}^\b
\times G_{n,m}^\a$.
For simplicity we set $G_{n,m}^{\b,\a}=G_{n,m}^\b \times G_{n,m}^\a$. Now
consider the following diagram
$$\xymatrix{
Q_{n, \geq m}^{\ab} & F_{n, \geq m} \ar[l]_-{i'} \ar[r]^-{\kappa'} &
Q_{n, \geq m}^\b \times
Q_{n, \geq m}^\a\\
G_{n,m}^{\ab} \times Q_{n,\geq m}^{\ab} \ar[u]^-{p_1} & G_{n,m}^{\b,\a} \times
F_{n, \geq m} \ar[u]^-{p_2} \ar[l]_-{i_1} \ar[r]^-{\kappa_1} & G_{n,m}^{\b,\a}
\times (Q_{n,\geq m}^\b
\times Q_{n,\geq m}^\a) \ar[u]^-{p_3}\\
& P_{n,m}^{\b,\a} \times F_{n, \geq m} \ar[ul]_-{i_1'} \ar[u]^-{p'_2}
\ar[ur]^-{\kappa'_1} &}
$$
where all squares are commutative. The rightmost square being
cartesian, we have
$p_3^* \kappa'_! {i'}^* =
\kappa_{1!}p_2^*{i'}^*=\kappa_{1!}i_1^*p_1^*.$ Furthermore, since
$p'_2$ is an affine fibration,
\begin{equation*}
\begin{split}
\kappa_{1!}&= \kappa_{1!}p'_{2*}{p'_2}^*\\
&=
\kappa_{1!}p'_{2!}{p'_{2}}^*[2\mathrm{dim}(P^{\b,\a}_{n,m}/G_{n,m}^{\b,\a})]\\
&=\kappa'_{1!}{p'_2}^*[2\mathrm{dim}(P^{\b,\a}_{n,m}/G_{n,m}^{\b,\a})]
\end{split}
\end{equation*}
so that
\begin{equation}\label{E:16}
p_3^*\kappa'_{!}{i'}^*=\kappa'_{1!}{i_1'}^*p_1^*[2\mathrm{dim}(P^{\b,\a}_{n,m}/G_{n,m}^{\b,\a})].
\end{equation}

Next, there are compatible surjective maps
\begin{equation}\label{E:19}
\xymatrix{
0 \ar[r] &\underline{V}_{n,m} \ar[r] \ar[d]_-{u^{\a}} & \mathcal{E}_{n,m}^{\ab}
\ar[r] \ar[d]_-{u^{\ab}_{n,m}} & \underline{W}_{n,m} \ar[r]
\ar[d]_-{u^{\b}}& 0\\
0 \ar[r] & \underline{V}_m \ar[r] & \mathcal{E}_m^{\ab} \ar[r] &
\underline{W}_m \ar[r] &0}
\end{equation}
where we put
$\underline{W}_{n,m}=\mathcal{E}_{n,m}^{\a+\b}/\underline{V}_{n,m}$
and
$\underline{W}_{m}=\mathcal{E}_{m}^{\a+\b}/\underline{V}_{m}$. Let
$T_{\ab} \subset
G_{n,m}^{\ab}$ (resp. $T_\a \subset G_{n,m}^\a$, resp. $T_\b \subset
G_{n,m}^\b$) be the
stabilizer of $\mathrm{Ker}\; u^{\ab}_{n,m}$ (resp. of
$\mathrm{Ker}\;u^{\a}$, resp. of
$\mathrm{Ker}\;u^{\b}$) and put $T_{\b,\a}=T_{\ab} \cap
P_{n,m}^{\b,\a}$. By construction,
we have a natural affine fibration $T_{\b,\a} \tto T_\b \times
T_{\a}$. We also fix
a section $T_{\b} \times T_\a \to T_{\b,\a}$. There is a projection
diagram with maps induced by (\ref{E:19})

$$\xymatrix{
G_{n,m}^{\ab} \times Q_{n,\geq m}^{\ab} \ar[d]_-{q_1} & P_{n,m}^{\b,\a} \times
F_{n, \geq m} \ar[d]_-{q_2} \ar[l]_-{i'_1} \ar[r]^-{\kappa'_1} &
G_{n,m}^{\b,\a}
\times (Q_{n,\geq m}^\b
\times Q_{n,\geq m}^\a) \ar[d]_-{q_3}\\
G_{n,m}^{\ab} \underset{T_{\ab}}{\times} Q_{n, \geq m}^{\ab} & P^{\b,\a}_{n,m}
\underset{T_{\b} \times T_{\a}}{\times} F_{n, \geq m} \ar[l]_-{i_2}
\ar[r]^-{\kappa_2}
\ar[d]_-{q'_2}&
G_{n,m}^{\b,\a} \underset{T_\b \times T_\a}{\times} (Q_{n,\geq m}^\b \times
Q_{n, \geq m}^\a)\\
& P_{n,m}^{\b,\a} \underset{T_{\b,\a}}{\times} F_{n, \geq m}
\ar[ul]_-{i_2'} \ar[ur]^-{\kappa'_2}
  &}
$$
The rightmost square being cartesian, we have
$q_3^*\kappa_{2!}i_2^*=\kappa'_{1!}q_2^*i_2^*=
\kappa'_{1!}{i'_1}^*q_1^*$. Furthermore, $q'_2$ is a vector bundle so
\begin{equation*}
\begin{split}
\kappa_{2!}i_2^*&=\kappa'_{2!}q'_{2!}{q_2'}^*{i'_2}^*\\
&=\kappa'_{2!}q'_{2*}{q_2'}^*{i'_2}^*[-2\mathrm{dim}(T_{\b,\a}/T_{\b}
\times T_{\a})]\\
&=\kappa'_{2!}{i'_2}^*[-2\mathrm{dim}(T_{\b,\a}/T_{\b} \times T_{\a})]
\end{split}
\end{equation*}
and finally
\begin{equation}\label{E:20}
q_3^* \kappa'_{2!}{i_2}^*[-2\mathrm{dim}(T_{\b,\a}/T_{\b} \times
T_{\a})]=\kappa'_{1!}
{i_1'}^*q_1^*.
\end{equation}
Now, if $P \in \mathcal{Q}_{G_n^{\ab}}(Q_{n, \geq m}^{\ab})$ then
$P'=ind_{T_{\ab}}^{G_{n,m}^{\ab}}(P)$
is characterized (up to isomorphism) by the relation $q_1^*P'=p_1^*P$
(and similarly for
$ind_{T_{\b,\a}}^{P_{n,m}^{\b,\a}}$, etc...). Thus, putting
(\ref{E:16}) and (\ref{E:20}) together
gives
\begin{equation}\label{E:21}
\begin{split}
ind_{T_{\b}\times T_{\a}}^{G_{n,m}^{\b,\a}} \kappa'_!
{i'}^*&=\kappa'_{2!}{i'_2}^*
ind_{T_{\ab}}^{G_{n,m}^{\ab}}\big[2
(\mathrm{dim}(P_{n,m}^{\b,\a}/T_{\b,\a} )- \mathrm{dim}
(G_{n,m}^{\b,a}/T_{\b} \times T_{\a}))\big]\\
&=\kappa'_{2!}{i'_2}^*
ind_{T_{\ab}}^{G_{n,m}^{\ab}}\big[2\sum_s
d_s(n,\a)(d_{s}(n,m,\beta)-d_s(n,\b))].
\end{split}
\end{equation}
As in Section~3 one sees that $P_{n,m}^{\b,\a}
\underset{T_{\b,\a}}{\times} F_{n, \geq m}$ is isomorphic
to the subvariety of $R_{n,m}^{\ab}$ defined by
\begin{equation*}
\begin{split}
R_{n,m}^{\b,\a}=\{(\phi: \mathcal{E}_{n,m}^{\ab} \tto &\mathcal{F})
\in R_{n,m}^{\ab}\;|\\
&\big[\phi(\bigoplus_{s,s'} \Hom(\mathcal{L}_s(n),
\mathcal{L}_{s'}(m)) \otimes V_{s'}^m \otimes
\mathcal{L}_s(n))]=\alpha\}
\end{split}
\end{equation*}
and that the natural map
$$\bigoplus_{s,s'}\Hom(\mathcal{L}_s(n), \mathcal{L}_{s'}(m)^{\oplus
d_{s'}(m,\ab)}) \otimes
\mathcal{L}_s(n)\tto \mathcal{E}_m^{\ab}$$
induces embeddings $F_m \hookrightarrow R^{\b,\a}_{n,m}, P_m^{\b,\a}
\hookrightarrow
P_{n,m}^{\b,\a}$ as well as a surjective morphism $\theta_2: P_{n,m}^{\b,\a}
\underset{P_{m}^{\b,\a}}{\times} F_m \tto R_{n,m}^{\b,\a}$.

\vspace{.05in}

Next observe the diagram
\begin{equation}\label{E:22}
\xymatrix{
G_{n,m}^{\ab} \underset{T_{\ab}}{\times} Q_{n, \geq m}^{\ab} & P^{\b,\a}_{n,m}
\underset{T_{\b,\a}}{\times} F_{n, \geq m} \ar[l]_-{i'_2} \ar[r]^-{\kappa'_2} &
G_{n,m}^{\b,\a} \underset{T_\b \times T_\a}{\times} (Q_{n,\geq m}^\b \times
Q_{n, \geq m}^\a)\\
G_{n,m}^{\ab} \underset{G_m^{\ab}}{\times} Q_{m}^{\ab}
\ar[u]^-{\theta_1} & P^{\b,\a}_{n,m}
\underset{P^{\b,\a}_m}{\times} F_{m} \ar[l]_-{i_3} \ar[r]^-{\kappa_3}
\ar[u]^-{\theta_2}&
G_{n,m}^{\b,\a} \underset{G_m^{\b} \times G^{\a}_m}{\times} (Q_{m}^\b \times
Q_{m}^\a)\ar[u]^-{\theta_3}}
\end{equation}

For simplicity, let us label $A_{11}, A_{12}, A_{21}, A_{22}$ the
vertices of the
rightmost square in (\ref{E:22}). A direct computation shows that the
natural map
$\theta'_2: A_{21} \to A_{11} \underset{A_{12}}{\times} A_{22}$ is a
vector bundle
of rank $r=\sum_s (d_s(n,m,\a)-d_s(n,\a))d_s(n,m,\b)$. Hence a
computation similar to
(\ref{E:16}) or (\ref{E:20}) gives
\begin{equation}\label{E:23}
\theta_3^* \kappa'_{2!}{i'_2}^*=\kappa_{3!}i_3^*\theta_{1}^*[2r].
\end{equation}
Finally, again as in (\ref{E:21}), we have an isomorphism of functors
\begin{equation}\label{E:24}
\kappa_{4!}i_4^* res_{G_{n,m}^{\ab}}^{G_m^{\ab}} =
res_{G_{n,m}^{\b,\a}}^{G_m^{\b} \times G_m^{\a}}
\kappa_{3!}i_3^*\big[2(\mathrm{dim}(P_{n,m}^{\b,\a}/P_{m}^{\b,\a})-\mathrm{dim}(G_{n,m}^{\b,\a}/G_m^{\b}
\times G_m^{\a}))\big],
\end{equation}
where
$$\xymatrix{
Q_m^{\ab} & F_m \ar[l]_-{i_4} \ar[r]^-{\kappa_4} & Q_m^\b \times Q_m^\a}$$
is the diagram for $\mathrm{Res}_{m}^{\b,\a}$.

Combining (\ref{E:15}), (\ref{E:21}), (\ref{E:23}) and (\ref{E:24})
yields (using simplified
notations)
\begin{equation*}
\begin{split}
(\widetilde{\Xi}_{n,m}^\b &\times \widetilde{\Xi}_{n,m}^\a) \circ
\widetilde{\mathrm{Res}}_n^{\b,\a}\\
&=res \circ \theta^* \circ ind \circ \kappa_! {i}^* j^*\\
&=res \circ \theta^* \kappa_! i^* \circ ind \circ j^*[2\sum_s
d_s(n,\a)(d_s(n,m,\beta)-d_s(n,\b))]\\
&=res \circ \kappa_! i^* \theta^* \circ ind \circ j^*[2\sum_s
d_s(n,m,\a)d_s(n,m,\beta)-d_s(n,\a)
d_s(n,\b))]\\
&=\kappa_! i^* \circ res \circ \theta^* \circ ind \circ j^*[2\sum_s
d_s(m,\a)d_s(m,\beta)-d_s(n,\a)
d_s(n,\b))]\\
&=\widetilde{\mathrm{Res}}_n^{\b,\a} \circ \widetilde{\Xi}_{n,m}^{\ab}
[2\sum_s d_s(m,\a)d_s(m,\beta)-d_s(n,\a)d_s(n,\b))]
\end{split}
\end{equation*}
as desired. \qed

\vspace{.2in}

\section{The algebra $\widehat{\mathfrak{U}}_{\mathbb{A}}$}

\vspace{.1in}

\paragraph{\textbf{5.1. Some simple Quot schemes.}}
  We start by describing the varieties $Q_n^\alpha$ in some important
examples. We will need the
following pieces of notation. For $i=1, \ldots, N$ there exists a
unique additive function
$deg_i: Pic(X,G) \to \Z/p_i\Z$ satisfying
$$deg_i(\alpha_{(i,j)})=1,\qquad
deg_i(\delta)=deg_i(\alpha_{(k,j)})=0\qquad \text{if\;} k \neq i$$
(recall that $\alpha_{(i,j)}$ is the class of the simple torsion
sheaf $S_i^{(j)}$
while $\delta$ is that of a generic simple torsion sheaf). For $s \in
S$ we have
$\langle [\mathcal{L}_s], \alpha_{(i,j)} \rangle=1$ if $j=deg_i(s)$ and
$\langle [\mathcal{L}_s],
\alpha_{(i,j)} \rangle=0$ otherwise.

For $i \in \{1, \ldots, ,N\}$ and $\eta=\sum_j \eta_j
\alpha_{(i,j)} \in K^+(X)$ we denote as in Section~1 by
$\mathcal{N}^{(p_i)}_\eta$ the space of all nilpotent
representations of the cyclic
quiver $A_{p_i-1}^{(1)}$ of dimension $(\eta_j) \in \N^{\Z/(p_i)\Z}$, on which
the group $G_{\eta}=\prod_j GL(\eta_j)$ naturally acts.
For $l \in
\N$ and $\eta=\sum_j \eta_j \alpha_{(i,j)}$
as above we have embeddings
\begin{align*}
\varphi_{l\delta}: GL(l) &\to \prod_s
\mathrm{Aut}(\mathcal{L}_s^{\oplus d_s(0,l\delta)})=\prod_s
\mathrm{Aut}(\mathcal{L}_s^{\oplus l})\\
   (x) &\mapsto (x,x, \ldots, x),\\
\varphi^i_{\eta}: G_{\eta} &\to \prod_s
\mathrm{Aut}(\mathcal{L}_s^{\oplus d_s(0,\eta)})=\prod_s \mathrm{Aut}
(\mathcal{L}_s^{\oplus \eta_{deg_i(s)}}),\\
(x_k) &\mapsto \prod_s x_{deg_i(s)}.
\end{align*}

Finally, in general let $l_1, \ldots, l_r \in \N$ and let $\eta_1,
\ldots, \eta_s $ be as above (for
corresponding values of the index $i_1, \ldots, i_s$). Setting
$\alpha=(\sum_j l_j)\delta + \sum_k \eta_k$, we denote
by
$$\varphi_{l_1\delta} * \cdots * \varphi_{l_r\delta} *
\varphi^{i_1}_{\eta_1} * \cdots * \varphi^{i_s}_{\eta_s}:\;
\prod_{j=1}^r GL(l_j) \times \prod_{k=1}^s G_{\eta_k} \to G_0^\alpha$$
the composition of the embedding
$$\prod_{j=1}^r \varphi_{l_j\delta} \times \prod_{k=1}^s
\varphi^{i_k}_{\eta_k}:\;
\prod_{j=1}^r GL(l_j) \times \prod_{k=1}^s G_{\eta_k} \to \prod_s
\big(\prod_{j=1}^r \mathrm{Aut}(\mathcal{L}_s^{\oplus l_j}) \times
\prod_{k=1}^s \mathrm{Aut}(\mathcal{L}_s^{\oplus d_s(0,\eta_k)})\big)$$
with the inclusion of the Levi factor
$$ \prod_s \big(\prod_{j=1}^r \mathrm{Aut}(\mathcal{L}_s^{\oplus l_j}) \times
\prod_{k=1}^s \mathrm{Aut}(\mathcal{L}_s^{\oplus d_s(0,\eta_k)})\big)
\to G_0^\alpha.$$
Such an embedding is clearly well-defined up to inner automorphism of
$G_0^\alpha$.

\vspace{.1in}

We now proceed with the description of the Quot schemes $Q_n^\alpha$
in some simple cases. Recall that if
$\alpha$ is the class of a torsion sheaf then $Q_n^\alpha \simeq
Q_0^\alpha$ for any $n \in \Z$. For simplicity,
we will simply write $Q^\alpha, G^\alpha,\cdots$ for $Q^\alpha_0,
G^\a_0$, etc. Also, by
the support of a sheaf $\mathcal{F} \in Coh_G(X)$, we mean the
support of the sheaf $\pi_*(\mathcal{F})$, where
$\pi: X \to \mathbb{P}^1$ is the quotient map.

\vspace{.1in}

\noindent
-Assume that $\alpha=\sum_j \eta_j\alpha_{(i,j)}$ for some $i$,
and that $\eta_j=0$ for at least one $j$. Then any sheaf
$\mathcal{F}$ with $[\mathcal{F}]=\alpha$
is supported at $\lambda_i$. In this case we have
$$Q^\alpha \simeq G^\alpha
\underset{\varphi_{\alpha}(G_\alpha)}{\times} \mathcal{N}^{(p_i)}_\alpha.$$
In particular, $Q^\alpha =\{pt\}$ if $\alpha=\alpha_{(i,j)}$ with $j
\neq 0$ and $Q^\alpha=PGL(M)$ where
$M=1+\sum_{k \neq i}(p_k-1)$ if $\alpha=\alpha_{(i,0)}$.

\vspace{.1in}

\noindent
-Next, assume that $\alpha=\delta$. We have $\mathcal{E}^\delta =\bigoplus_s
\mathcal{L}_s$, so that $G^\delta \simeq \prod_s \C^*$.
The assignment $(\phi: \mathcal{E}^\delta \tto \mathcal{F}) \mapsto
supp(\mathcal{F})$ induces a flat morphism
$\rho: Q^\delta \to \mathbb{P}^1$ and we have
\begin{align*}
\rho^{-1}(t) &\simeq G^\delta/\varphi_\delta(\C^*), \qquad
\mathrm{if}\; t \not\in \Lambda,\\
\rho^{-1}(\lambda_i) &\simeq G^\delta
\underset{\varphi^i_{\delta}((\C^*)^{\Z/p_i\Z})}{\times} \mathcal{N}^{(p_i)}_\delta.
\end{align*}
Observe that $\mathcal{N}^{(p_i)}_\delta$ has $p_i$ open orbits, and hence
$p^{-1}(\lambda_i)$ has $p_i$
irreducible components.
In the simplest case $(X,G)=(\mathbb{P}^1, Id)$ we have
$Q^\delta \simeq \mathbb{P}^1$.

\vspace{.1in}

\noindent
-More generally, assume $\alpha=l \delta$ for some
$l \in \N$. We have $\mathcal{E}^{l\delta}=\bigoplus_s
\mathcal{L}_s^{\oplus l}$
and $G^{l\delta} \simeq \prod_s GL(l)$.
Let us consider the fibers of the support map $\rho_l: Q^{l\delta}
\to S^l \mathbb{P}^1$
(here the support is counted with the multiplicity given
by the class in $K(X)$). If
$t \in \mathbb{P}^1\backslash\Lambda$ then
$$
\rho_l^{-1}(t, \cdots, t)\simeq
G^{l\delta} \underset{\varphi_{l\delta}(GL(l))}{\times} \mathcal{N}_l,$$
where $\mathcal{N}_l= \mathcal{N}^{(1)}_l\subset \mathfrak{gl}_l$ is the nilpotent cone, while
$$\rho_l^{-1}(\lambda_i, \cdots, \lambda_i) \simeq G^{l\delta}
\underset{\varphi^i_{l\delta}(G_{l\delta})}
{\times} \mathcal{N}^{(p_i)}_{l\delta}.$$
The fiber at a general point looks like a product of fibers of
the previous type for smaller values of $l$. Namely, if
$\underline{t}=(t_1^{l_1},\ldots, t_r^{l_r}, \lambda_1^{n_1},
\cdots, \lambda_N^{n_N}) \in S^l \mathbb{P}^1$ with $\sum l_j + \sum
n_k=l$ then
$$\rho_l^{-1}(\underline{t})=G^{l\delta} \underset{H}{\times}
\bigg(\prod_j \mathcal{N}_{l_j} \times \prod_k
\mathcal{N}^{(p_k)}_{n_k\delta}\bigg)$$
where $H= \varphi_{l_1\delta}* \cdots * \varphi_{l_r\delta} *
\varphi_{n_1\delta}^1 *
\cdots * \varphi^N_{n_N\delta}(\prod_j GL(l_j) \times \prod_k
G_{n_k\delta})$.\\

As an example of a variety $Q^{l\delta}$, let us again assume that
$(X,G)=(\mathbb{P}^1, Id)$.
We then have $Q^{l\delta} \simeq Gr(l, 2l)$.
Via the embedding $\mathfrak{gl}_l \to Gr(l,2l)$, $g \mapsto
\text{graph\;of\;}g$, the support map restricts
to $\rho_l: \mathfrak{gl}_l \to S^l\C$, $g \mapsto
\{\text{eigenvalues\;of\;}g\}$.

\vspace{.1in}

\noindent
We will need the following subvarieties of $Q^{\alpha}$, defined for
an arbitrary torsion class $\alpha$.
Recall that if $\gamma =\sum_{(i,j) \in \aleph} \eta_{i,j}
\alpha_{(i,j)}$ then
$|\gamma|=\sum_{i,j} \frac{\eta_{i,j}}{p_i}$. If $\mathcal{F}$ is a
torsion sheaf then $\mathcal{F}=\mathcal{F}' \oplus
\mathcal{F}_{\Lambda}$ where $supp(\mathcal{F}') \subset
\mathbb{P}^1\backslash \Lambda$ and $supp(\mathcal{F}_\Lambda)
\subset \Lambda$.We also put
$$U^{\alpha}_{\geq d}=\{(\phi: \mathcal{E}^{\alpha}_0 \tto
\mathcal{F})\;|\;
|[\mathcal{F}_\Lambda]| \geq d\}.$$
It is clear that $U^{\alpha}_{\geq d}$ defines a decreasing filtration
of $Q^{\alpha}$, and that
$$U^{\alpha}:=U^{\alpha}_{\geq 0} \backslash U^{\alpha}_{\geq 1}=
\rho^{-1}(S^l(\mathbb{P}^1\backslash \Lambda))$$
is nonempty if and only if $\alpha=l\delta$ for some $l$, in which
case it is an open set in $Q^{l\delta}$. We construct a
stratification of $U^{l\delta}$ as follows. Let
$\underline{\lambda}=(\lambda^{(1)}, \ldots ,\lambda^{(r)})$ be an
$r$-tuple of partitions such that $\sum_{j,k} \lambda^{(k)}_j=l$. We put
\begin{equation*}
\begin{split}
U^{l\delta}_{\underline{\lambda}}=\{(\phi: \mathcal{E}^{l\delta} \tto
\mathcal{F}\;|
\mathcal{F}\simeq \bigoplus_{k=1}^r &\bigoplus_{j}
\mathcal{O}_{x_k}^{(\lambda^{(k)}_j)}\\
& \text{for\;some\;distinct\;}
x_1, \ldots, x_r \in \mathcal{P}^1\backslash\Lambda\}
\end{split}
\end{equation*}
where we denote by $\mathcal{O}_x^{(n)}$ the indecomposable torsion
sheaf supported at $x$ of length $n$. Observe that
$U^{l\delta}_{((1), \ldots, (1))}$ is a smooth open subvariety of
$Q^{l\delta}$.\\

 From the explicit description above, it follows that the fibers of
the map $\rho$ are contractible.
Thus, the embedding of
$U^{l\delta}_{((1), \ldots, (1))}$ in
$\rho_l^{-1}(S^l(\mathbb{P}^1)\backslash \{t_i=t_j\})$
gives a projection
$\pi_1(U^{l\delta}_{((1),\ldots,(1))})\tto
\pi_1(S^l(\mathbb{P}^1)\backslash \{t_i=t_j\})=B_l$,
where $B_l$ denotes the braid
group of order $l$ of $\mathbb{P}^1$. Hence there is a canonical projection
$\pi_1(U^{l\delta}_{((1),\ldots,(1))})
\tto \mathfrak{S}_l$, where $\mathfrak{S}_l$ is the symmetric group on
$l$ letters.\\

Continuing with the example above of $(X,G)=(\mathbb{P}^1,Id)$ we see
that $U^{l\delta}_{\underline{\lambda}}$ is the set of elements $x \in
\mathfrak{gl}_n$ whose Jordan decomposition is given by
$\underline{\lambda}$. In particular,
$U^{l\delta}_{((1),\ldots,(1))} \cap \mathfrak{gl}_l$ is the set of
semisimple elements with distinct
eigenvalues.

\vspace{.1in}

\noindent -Finally, let $\alpha$ be the class of a line bundle
$\mathcal{L}$ of degree $\alpha_0 \in K^+(X)/\Z[\mathcal{O}_X]$.
If $\alpha_0 \not\geq n\delta$ then $\mathcal{E}_n^\alpha=0$ and
$Q_n^\alpha$ is empty. Otherwise, $Q_n^\alpha$ contains a unique
open (dense) $G_n^\alpha$-orbit corresponding to $\mathcal{L}$ and
infinitely many smaller dimensional orbits corresponding to
sheaves of the form $\mathcal{L}' \oplus \mathcal{T}$ where
$\mathcal{L}'$ is a line bundle and $\mathcal{T}$ is a torsion
sheaf. More precisely, there exists a stratification
$Q_n^\a=\bigsqcup_{[\mathcal{L}']} Q_n^\a([\mathcal{L}'])$ where
the sum ranges over the (finite) set of classes of line bundles
$\mathcal{L}'$ of degree $d$ satisfying $n\delta \leq d \leq
[\mathcal{L}]$, and
$$Q_n^\a([\mathcal{L}'])=\{\phi: \mathcal{E}_n^\a \tto
\mathcal{F}\;|\; \nu(\mathcal{F}) \simeq \mathcal{L}'\}.$$ Each
$Q_n^\a([\mathcal{L}'])$ is in turn described as follows. We put
$\a'=[\mathcal{L}']$ and $\beta=\alpha-\a'$. For each $s \in S$,
let us fix a subspace $V_s \subset \mathcal{E}^\a_{n,s}$ of
dimension $\langle [\mathcal{L}_s(n)], \beta \rangle$ together
with isomorphisms $a_s:\; V_s \simeq \mathcal{E}_{n,s}^\beta$;
$b_s: \mathcal{E}^\a_{n,s}/V_s \simeq \mathcal{E}^{\a'}_{n,s}$.
Finally, we put $\underline{V}=\bigoplus_s V_s \otimes
\mathcal{L}_s(n)$ and
$$S_n^\a(\mathcal{L}')=\{ (\phi: \mathcal{E}^\a_{n} \tto
\mathcal{F}) \in Q_n^\a([\mathcal{L}'])\;|\; p \circ
\phi(\underline{V})=0\},$$ where $p: \mathcal{F} \tto
\mathcal{F}/\tau(\mathcal{F})$ is the canonical surjection. Then,
denoting by $P \subset G_n^\a$ the stabilizer of $\underline{V}$,
we have $Q_n^\a([\mathcal{L}']) \simeq G_n^\a \underset{P}{\times}
S_n^\a(\mathcal{L}')$ and the morphism
\begin{align*}
S_n^\a(\mathcal{L}') & \to Q_n^\b \times
Q_n^{\a'}([\mathcal{L}'])\\
\phi & \mapsto (a_* \phi_{|\underline{V}},
b_*\phi_{|\mathcal{E}_n^\a/\underline{V}})
\end{align*}
is an affine fibration. Note that $Q_n^{\a'}([\mathcal{L}'])
\simeq G_n^{\a'}/k^*$ consists of a single $G_n^{\a'}$-orbit, and
thus the description of $Q_n^\b$ for a torsion sheaf $\beta$
essentially carries over to each strata $Q_n^\a([\mathcal{L}'])$
of $Q_n^\a$.

\vspace{.2in}

\paragraph{\textbf{5.2. A class of simple perverse sheaves.}} In this
section, we construct a full subcategory
$\mathbb{U}^\alpha$ of $\mathbb{Q}^\alpha$ for each $\alpha \in K^+(X)$.
Set
$$\mho=\{\alpha_{(i,j)}, (i,j) \in \aleph\}
\cup \{\delta\} \cup \{[\mathcal{L}], \mathcal{L} \in Pic(X,G)\}
\subset K^+(X).$$
For any $\alpha \in K^+(X)$, the constant sheaf
$\qlb_{Q_n^\alpha}[\mathrm{dim}\;Q_n^\alpha]$ belongs
to $\mathcal{Q}_{G_n^\alpha}(Q_n^\alpha)$ and for each $n \leq m$ we
have $\Xi_{n,m}^\alpha
(\qlb_{Q_n^\alpha}[\mathrm{dim}\;Q_n^\alpha])=\qlb_{Q_m^\alpha}[\mathrm{dim}\;Q_m^\alpha]$.
This gives rise to
a well-defined simple object $\mathbf{1}_\alpha
=(\qlb_{Q_n^\alpha}[\mathrm{dim}\;Q_n^\alpha])_{n \in \Z}$ of
$\Q^\alpha$. 

\vspace{.1in}

Let $\mathcal{P}^\alpha$ be the set of all simple objects of
$\Q^\a$ appearing (possibly with a shift) in an induction product
$\mathrm{Ind}^{\alpha_r, \ldots, \alpha_1}( \mathbf{1}_{\alpha_r}
\boxtimes \cdots \boxtimes \mathbf{1}_{\alpha_1})$ with $\alpha_1,
\ldots, \alpha_r \in \mho$ and $\a_1+ \cdots + \a_r=\a$. We define
$\mathbb{U}^\a$ as the full subcategory of $\mathbb{Q}^{\a}$
consisting of all admissible sums $\bigoplus_h \mathbb{P}_h[d_h]$
with $\mathbb{P}_h \in \mathcal{P}^\a$ and $d_h \in \Z$.

\vspace{.1in}

We also define a category $\mathbb{U}^\b \hat{\boxtimes}
\mathbb{U}^\a$ as follows. First, using the collection of functors
$\Xi^\b_{n,m} \times \Xi^{\a}_{n,m}: \mathcal{Q}_{G_{n}^\b
\times G_{n}^\a}(Q_n^\b \times Q_n^\a)
\to \mathcal{Q}_{G_{m}^\b \times G_{m}^\a}(Q_m^\b \times Q_m^\a)$ we
may define a category
$\mathbb{Q}^{\b,\a}$ by the same method as in Section~3.3. We let
$\mathbb{U}^\b \hat{\boxtimes}
\mathbb{U}^\a$ be the full subcategory of $\mathbb{Q}^{\b,\a}$
consisting of objects which are
admissible sums of objects of the form $\mathbb{P}_\b \boxtimes
\mathbb{P}_{\a}$ with
$\mathbb{P}_\b \in \mathbb{U}^\b$ and $\mathbb{P}_\a \in
\mathbb{U}^\a$.

\vspace{.1in}

\begin{lem}\label{L:14} For any $\alpha, \beta \in K^+(X)$, the
induction and restriction functors give rise to functors
$$\mathrm{Ind}^{\b,\a}:\; \mathbb{U}^\b \boxtimes
\mathbb{U}^\a \to \mathbb{U}^{\alpha+\beta},.$$
$$\mathrm{Res}^{\b,\a}:\; \mathbb{U}^{\ab} \to \mathbb{U}^\b
\hat{\boxtimes} \mathbb{U}^\a.$$
\end{lem}
\noindent
\textit{Proof.} The statement concerning the induction functor is
clear from the definitions and from Section~4.1.
We prove the second statement. Let us first show that
for any $\mathbb{P}=(P_n)_{n \in \Z}$ in $\mathbb{U}^{\ab}$ and any
$n \in \Z$ we have
$\mathrm{Res}_n^{\b,\a}(P_n) \in
\mathcal{Q}_{G_n^\b}(Q_n^\b) \boxtimes
\mathcal{Q}_{G_n^\a}(Q_n^\a)$.
It is clearly enough to prove this for
$\mathbb{P}=\widetilde{\mathrm{Ind}}^{\gamma_r, \ldots,
\gamma_1}(\mathbf{1}_{\gamma_r}
\boxtimes\cdots\boxtimes \mathbf{1}_{\gamma_1})$ with $\sum
\gamma_i=\gamma:=\alpha+\beta$, and $\gamma_i \in \mho$. Set
\begin{equation*}
\begin{split}
E''_l=\{(\phi: \mathcal{E}_l^\gamma \tto \mathcal{F},
\underline{V}_1&\subset \cdots \subset \underline{V}_r
=\mathcal{E}_l^\gamma)\;|\\
&\;\underline{V}_i/\underline{V}_{i-1} \simeq \mathcal{E}_l^{\gamma_i},
[\phi(\underline{V}_i)/\phi(\underline{V}_{i-1}]=\gamma_i\}
\end{split}
\end{equation*}
and let $p_3: E''_l \to Q_l^\gamma$ be the projection. By
Lemma~\ref{L:9}, the variety $E''_l$ is smooth.
Using the definitions (see Lemma~\ref{L:12}), we have $P_n \simeq
\Xi_{l,n} p_{3!}(\qlb_{E''_l})$ for $l \ll n$.
Now, let us fix a subsheaf $\underline{V}^\alpha \subset
\mathcal{E}_l^\gamma$ as well as identifications
$\underline{V}^\alpha \stackrel{\sim}{\to} \mathcal{E}_l^\alpha$,
$\mathcal{E}^\gamma_l/\underline{V}^\alpha
\stackrel{\sim}{\to} \mathcal{E}^\beta_l$, and let us
consider the two varieties
\begin{align*}
F''_l=&\{\big(\phi: \mathcal{E}_l^\gamma \tto \mathcal{F},
\underline{V}_1\subset \cdots \subset \underline{V}_r
=\mathcal{E}_l^\gamma) \in E''_l\;|\; [\phi(\underline{V}^\alpha)]=\alpha\},\\
F_l=&\{ (\phi: \mathcal{E}_l^\gamma \tto \mathcal{F}) \in
Q_l^\gamma\;|\;[\phi(\underline{V}^\alpha)]=\alpha\}.
\end{align*}
These fit together in a diagram
\begin{equation}\label{E:50}
\xymatrix{
E''_l \ar[d]_-{p_3} & F''_l \ar[l]_-{i'} \ar[d]_-{p'_3}&\\
Q_l^\gamma & F_l \ar[l]_-{i} \ar[r]^-{\kappa} & Q_l^\b \times Q_l^\a}
\end{equation}
The (singular) variety $F''_l$ admits a smooth stratification
constructed as follows. Fix $\mathbf{a}=(\alpha_1, \ldots,
\alpha_r)$ and $\mathbf{b}=(\beta_1, \ldots \beta_r)$ such that
$\alpha_i+\beta_i=\gamma_i$ for all $i$
and $\sum_i \alpha_i=\alpha$. The subvariety
\begin{equation*}
\begin{split}
F''_{l}(\mathbf{b},\mathbf{a})=\{&\big(\phi: \mathcal{E}_l^\gamma
\tto \mathcal{F}, \underline{V}_1\subset \cdots
\subset \underline{V}_r=\mathcal{E}_l^\gamma) \in E''_l\;|\;\\
&(\underline{V}_i \cap \underline{V}^\alpha) /
(\underline{V}_{i-1} \cap \underline{V}^\alpha) \simeq
\mathcal{E}^{\alpha_i}_n,
[\phi(\underline{V}_i \cap
\underline{V}^\alpha)/\phi(\underline{V}_{i-1}\cap
\underline{V}^\alpha)]=\alpha_i\}
\end{split}
\end{equation*}
is smooth by Lemma~\ref{L:9}, and we clearly have $F''_l =
\bigcup_{\mathbf{a},\mathbf{b}} F''_l(\mathbf{b},\mathbf{a})$.
Furthermore, the restriction of $\kappa p'_3$ to
$F''_l(\mathbf{b},\mathbf{a})$ factors as a composition
$$\xymatrix{F''_l(\mathbf{b},\mathbf{a})
\ar[r]^-{\kappa(\mathbf{b},\mathbf{a})} & E''_l(\mathbf{b})
\times E''_l(\mathbf{a}) \ar[r]^-{p''_3 \times p''_3}& Q_l^\b
\times Q_l^\a},$$
where: $E''_l(\mathbf{b})$ (resp. $E''_l(\mathbf{a})$) is defined as
$E''_l$ by replacing
$\gamma$ and $\gamma_i$ by $\b$ and $\b_i$ (resp. by $\a$
and $\a_i$); $\kappa(\mathbf{b},\mathbf{a}):\;
(\phi,\underline{V}_i) \mapsto ((\phi_{|\mathcal{E}^\gamma_l/\underline{V}^\alpha}, (\underline{V}_i
+ \underline{V}^\alpha/\underline{V}^\alpha)_i),(\phi_{|\underline{V}_{\alpha}},
(\underline{V}_i \cap \underline{V}^\alpha)_i))$
is a vector bundle; $p''_3$ is the natural projection. Let
(\ref{E:50}$)_{\geq n}$ be the
restriction of (\ref{E:50}) to the open subvarieties $E''_{l,\geq n},
Q^\gamma_{l,\geq n}$, etc...
By Claim~\ref{C:1}, the square in (\ref{E:50}$)_{\geq n}$ is still
cartesian. Without risk of confusion,
let us denote by $j_{l,n}$ all the open embeddings
$Q_{l,\geq n}^\gamma \to Q_l^\gamma, F_{l,\geq n} \to F_l,$ etc...
Then $j_{l,n}^*p_{3!}(\qlb_{E''_l})=p_{3!}
(\qlb_{E''_{l,\geq n}})$ and $\kappa_!i^* p_{3!}(\qlb_{E''_{l,\geq
n}})=\kappa_{!}p'_{3!}{i'}^*
(\qlb_{E''_{l,\geq n}})=\kappa_!p'_{3!}(\qlb_{F''_{l,\geq n}})$.
Note that the restriction
of $\kappa(\mathbf{b},\mathbf{a})$ to $F''_{l,\geq n}$ is still a
vector bundle and that, for $l \ll n$, the
restriction of $p''_3$ to $E''_{l,\geq n}(\mathbf{a})$ and
$E''_{l,\geq n}(\mathbf{b})$ is proper. Using
\cite{L1}, 8.1.6, we deduce that
$$\kappa_!p'_{3!}(\qlb_{F''_{l,\geq n}})\simeq
\bigoplus_{\mathbf{a},\mathbf{b}} j_{l,n}^*(p''_{3!}\times p''_{3!})
(\qlb_{E''_{l}(\mathbf{b})}\boxtimes
\qlb_{E''_{l}(\mathbf{a})})[-2r_{\mathbf{a},\mathbf{b}}]$$
where $r_{\mathbf{a},\mathbf{b}}$ is the rank of
$\kappa(\mathbf{b},\mathbf{a})$.
Now, applying Lemma~\ref{L:12} we obtain
\begin{equation*}
\begin{split}
\widetilde{\mathrm{Res}}_n^{\b,\a} P_n
&=\widetilde{\mathrm{Res}}_n^{\b,\a} \Xi_{l,n} p_{3!}(\qlb_{E''_l})\\
&=\Xi_{l,n} \kappa_{!}i^* p_{3!}(\qlb_{E''_l})\\
&=\Xi_{l,n} \big(\bigoplus_{\mathbf{b},\mathbf{a}} (p''_{3!}\times p''_{3!})
(\qlb_{E''_{l}(\mathbf{b})}\boxtimes
\qlb_{E''_{l}(\mathbf{a})})[-2r_{\mathbf{a},\mathbf{b}}]\big)
\end{split}
\end{equation*}
Finally, notice that $\Xi_{l,n} p''_{3!}(\qlb_{E''_l(\mathbf{a})})=
\widetilde{\mathrm{Ind}}^{\alpha_r,
\ldots,\alpha_1}_{l,n}(\qlb_{Q_l^{\a_r}} \boxtimes \cdots \boxtimes
\qlb_{Q_l^{\a_1}})$. We conclude the proof of the Lemma by showing
the following result.

\begin{lem}\label{L:15} Let $\alpha \leq \gamma$ with $\gamma \in
\mho$. Then $\mathbf{1}_{\alpha}$  belongs to $\mathbb{U}^\alpha$.
In particular, $\mathbf{1}_{\a}$ belongs to $\mathbb{U}^\a$ for
any torsion class $\alpha$.
\end{lem}
\noindent \textit{Proof.} If $\alpha$ is the class of a vector
bundle then $\alpha \in \mho$ and there is nothing to prove. So we
may assume that $\alpha$ is the class of a torsion sheaf. Let us
first consider the case where $\alpha \not\geq \delta$ and
$\alpha$ is supported at some point $\lambda_i$, i.e there exists
$i \in \{1, \ldots, N\}$ and $\eta_j \in \mathbb{N}$ with
$\eta_j=0$ for at least one $j$, such that $\alpha=\sum_j
\eta_j\alpha_{(i,j)}$. In this case $Q^\alpha$ is described in
terms of the category of representations of the cyclic quiver, and
the product corresponds to the usual Hall algebra product.
It is then easy to see that $\mathbf{1}_{\alpha}$
appears in a product of the form
$\mathrm{Ind}^{\alpha_{(i,j_r)},\ldots ,
\alpha_{(i,j_1)}}(\mathbf{1}_{\alpha_{(i,j_r)}}\boxtimes \cdots
\boxtimes \mathbf{1}_{\alpha_{(i,j_1)}})$. 
Now let us assume that $\alpha=l\delta$ for some $l>0$. 
Let us consider the product $A=\mathrm{Ind}^{\delta, \ldots,
\delta}(\mathbf{1}_{\delta} \boxtimes \cdots \boxtimes
\mathbf{1}_{\delta})$. The only subsheaves of class $\delta$ of a
sheaf $\mathcal{F}=\bigoplus_{t=1}^l\mathcal{O}_{x_t}$
corresponding to a point of $U^{l\delta}_{((1),\ldots,(1))}$ (so
that $x_1\ldots, x_l$ are distinct) are the $\mathcal{O}_{x_t}$.
It follows that the stalk of $A$ over such points is of rank $l!$,
and that the monodromy representation of
$\pi_1(U^{l\delta}_{((1),\ldots, (1))})$ on that space factors
through  the composition $\pi_1(U_{l\delta, ((1),\ldots, (1))})
\tto B_l \tto \mathfrak{S}_l$ and is equal to the regular
representation of $\mathfrak{S}_l$. In particular,
$\mathbf{IC}(U^{l\delta}_{((1),\ldots,(1))},\mathbf{1})=\mathbf{1}_{Q^{l\delta}}$
appears in $A$, and thus
belongs to $\mathbb{U}^\alpha$.\\
In the general case, we may write $\alpha=\sum_{i=1}^N \alpha_i +
l\delta$ where $\alpha_i \not\geq \delta$ is supported at
$\lambda_i$. We consider the product $B=\mathrm{Ind}^{\alpha_N,
\ldots, \alpha_1, l\delta}(\mathbf{1}_{l\delta} \boxtimes \cdots
\boxtimes \mathbf{1}_{{\alpha_1}})$. The set $U$ of points of
$Q^\alpha$ corresponding to sheaves isomorphic to
$\mathcal{F}\oplus \bigoplus_{i=1}^N \mathcal{F}_i$, for some
sheaf $\mathcal{F}$ supported on $\mathbb{P}^1\backslash \Lambda$
and for some torsion sheaves $\mathcal{F}_i$ of class
$[\mathcal{F}_i]=\alpha_i$, forms a dense open set. Observe that
$\mathcal{F}_i$ is the only subsheaf of class $\alpha_i$ of such a
sheaf. As before, this implies that
$\mathbf{IC}(U,\mathbf{1})=\mathbf{1}_{Q^\alpha}$ belongs to
$\mathbb{U}^\alpha$, as desired.~\qed

\vspace{.1in}

As for the induction product (see Lemma~\ref{L:12}), there exists a
canonical natural transformation
$\mathrm{Res}^{\gamma,\beta} \circ \mathrm{Res}^{\beta+\gamma,\alpha}
\simeq \mathrm{Res}^{\b,\a} \circ
\mathrm{Res}^{\gamma,\alpha+\beta}$ for any triple
$(\alpha,\beta,\gamma) \in K^+(X)^3$. This allows
us to define an iterated restriction functor
$$\mathrm{Res}^{\alpha_r,\ldots,\alpha_1}:
\mathbb{U}^{\alpha_1+\cdots+\alpha_r} \to \mathbb{U}^{\alpha_r}
\hat{\boxtimes}
\cdots \hat{\boxtimes} \mathbb{U}^{\alpha_1}.$$

\vspace{.1in}

\paragraph{\textbf{Remark.}} Fix $\a, \b,\gamma,\delta \in K_x$ such
that $\a+\b=\gamma+\delta$,
and put $\mathcal{I}=\{(\a_1, \a_2,
\b_1,\b_2) \in (K^+(X))^4\;| \gamma=\a_1+\b_1,
\delta=\a_2+\b_2,\alpha=\a_1+\a_2,
\b=\b_1+\b_2\}$.
Let $P_{\gamma} \in \mathcal{Q}_{G_n^\gamma}(Q_n^\gamma)$ and $P_{\delta} \in
\mathcal{Q}_{G_n^\delta}(Q_n^\delta)$. It is possible to show that,
for $m \gg n$,
\begin{equation*}
\begin{split}
\mathrm{Res}_m^{\b,\a}&(\mathrm{Ind}_{n,m}^{\delta,\gamma}(P_\delta
\boxtimes P_\gamma))\\
&=
\bigoplus_{(\a_i,\b_i) \in \mathcal{I}} \mathrm{Ind}_{n,m}^{\b_2,\b_1}
\times \mathrm{Ind}^{\a_2,\a_1}_{n,m}
\big(
\mathrm{Res}_n^{\b_2,\a_2}(P_\delta)\boxtimes
\mathrm{Res}_n^{\b_1,\a_1}(P_\gamma)
\big)[-2r_{\mathbf{a},\mathbf{b}}]
\end{split}
\end{equation*}
where $\mathbf{a}=(\a_1,\a_2)$ and $\mathbf{b}=(\b_1,\b_2).$

\vspace{.2in}

\paragraph{\textbf{5.4. The algebra $\mathfrak{U}_{\mathbb{A}}$.}}
Set $\mathbb{A}=\C[v,v^{-1}]$. Following Lusztig, we consider the
free $\mathbb{A}$-module $\mathfrak{U}_{\mathbb{A}}$ generated by
elements $\mathbf{b}_{\mathbb{P}}$ where $\mathbb{P}$ runs through
$\mathcal{P}=\bigsqcup_{\gamma} \mathcal{P}^{\gamma}$, modulo the
relation $\mathbf{b}_{\mathbb{P}[1]}=v
\mathbf{b}_{\mathbb{P}}$. The $\mathbb{A}$-module has a natural
$K^+_X$-gradation $\mathfrak{U}_{\mathbb{A}}= \bigoplus_{\gamma}
\mathfrak{U}_{\mathbb{A}}[\gamma]$. Let us say that an element
$\mathbf{b}_{\mathbb{P}}$ is of $h$-degree $n$ if
$\mathbb{P}=(P_m)_{m \in \Z}$ with $P_m=0$ for $m > -n$ and
$P_{-n} \neq 0$. We denote by
$\widehat{\mathfrak{U}}_{\mathbb{A}}$ the completion of
$\mathfrak{U}_{\mathbb{A}}$ with respect to the $h$-adic topology.
We also extend by linearity the notation $\mathbf{b}_{\mathbb{P}}$
to an arbitrary (admissible) complex $\mathbb{P} \in
\mathbb{U}^{\alpha},\; \alpha \in K^+(X)$. Finally, if $\a \in
K^+(X)$ and $\mathbf{1}_{\a} \in \mathcal{P}$ then
we set $\mathbf{b}_\alpha=\mathbf{b}_{\mathbf{1}_\a}$.

\vspace{.1in}

We endow the space $\widehat{\mathfrak{U}}_{\mathbb{A}}$ with an
algebra and a coalgebra structure as follows.
First, by Lemma~\ref{L:14} the functors $\mathrm{Ind}^{\alpha,\beta}$
and $\mathrm{Res}^{\a,\b}$ give rise to
well-defined $\mathbb{A}$-linear maps
$$m_{\b,\a}: \mathfrak{U}_{\mathbb{A}}[\b]
\otimes\mathfrak{U}_{\mathbb{A}}[\a] \to
\widehat{\mathfrak{U}}_{\mathbb{A}}[\ab],$$
$$\Delta_{\b,\a}: \mathfrak{U}_{\mathbb{A}}[\ab] \to
\widehat{\mathfrak{U}}_{\mathbb{A}}[\b] \hat{\otimes}
\widehat{\mathfrak{U}}_{\mathbb{A}}[\a].$$
Furthermore, these maps are continuous by Lemma~\ref{L:11} and
Lemma~\ref{L:13}, and we may extend them
to continuous maps
$$m_{\b,\a}: \widehat{\mathfrak{U}}_{\mathbb{A}}[\b]
\otimes\widehat{\mathfrak{U}}_{\mathbb{A}}[\a] \to
\widehat{\mathfrak{U}}_{\mathbb{A}}[\ab],$$
$$\Delta_{\b,\a}: \widehat{\mathfrak{U}}_{\mathbb{A}}[\ab] \to
\widehat{\mathfrak{U}}_{\mathbb{A}}[\b] \hat{\otimes}
\widehat{\mathfrak{U}}_{\mathbb{A}}[\a].$$
The associativity of $m=\bigoplus_{\a,\b} m_{\b,\a}$ follows from
Lemma~\ref{L:12}, and the coassociativity
of $\Delta=\bigoplus_{\a,\b} \Delta_{\b,\a}$ is proved in a similar way.
Note that $\Delta$ is an algebra morphism only after a suitable twist, as in \cite{L1}.

\vspace{.1in}

By definition, $\widehat{\mathfrak{U}}_{\mathbb{A}}$ comes
equipped with a (topological) basis
$\{\mathbf{b}_{\mathbb{P}}\}_{\mathbb{P} \in
\mathcal{P}}$. The following key result of the paper will
be proved in Sections 6-9. Let $K(X)^{tor}=\sum_{(i,j) \in \aleph}
\N \alpha_{(i,j)}$ be the set of classes of torsion sheaves.

\begin{theo}\label{T:5} i)The subalgebra
${\mathfrak{U}}_{\mathbb{A}}^{tor}=\bigoplus_{\alpha \in
K(X)^{tor}} \mathfrak{U}_{\mathbb{A}}[\alpha]$ is
generated by elements $\mathbf{b}_{l\alpha_{(i,j)}}$
and $\mathbf{b}_{l\delta}$ for $(i,j) \in \aleph$ and $l \geq 1$.\\
ii)Assume that $X$ is of genus zero. Then
$\widehat{\mathfrak{U}}_{\mathbb{A}}$ is topologically generated
by ${\mathfrak{U}}_{\mathbb{A}}^{tor}$ and the elements
$\mathbf{b}_{l[\mathcal{O}_X(n)]}$ for $n \in \Z$ and $l \in
\mathbb{N}$; i.e the subalgebra generated by
${\mathfrak{U}}_{\mathbb{A}}^{tor}$ and the elements
$\mathbf{b}_{l[\mathcal{O}_X(n)]}$ is dense in
$\widehat{\mathfrak{U}}_{\mathbb{A}}$.\end{theo}

\vspace{.1in}

We conjecture that the statement ii) holds with no restrictions on the
genus of the curve $X$ (see \cite{S2} for the case of an elliptic curve $X$).

\vspace{.2in}

\section{Harder-Narasimhan filtration}

\vspace{.1in}

\paragraph{}In this section we collect certain results pertaining to the
Harder-Narasimhan filtration for objects of $Coh_G(X)$, which will be used
in the proof of Theorem~5.1 ii). \textit{We assume that} $X \simeq \mathbb{P}^1$.

\vspace{.1in}

\paragraph{\textbf{6.1. The HN filtration.}} Define the slope of a
coherent sheaf $\mathcal{F} \in Coh_G(X)$ by
$$\mu(\mathcal{F})=\mu([\mathcal{F}])=\frac{|\mathcal{F}|}{rank(\mathcal{F})}
\in \mathbb{Q} \cup \{\infty\}.$$
We say that a sheaf $\mathcal{F}$ is semistable of slope $\mu$ if
$\mu(\mathcal{F})=\mu$ and if for any subsheaf
$\mathcal{G} \subset \mathcal{F}$ we have $\mu(\mathcal{G}) \leq
\mu(\mathcal{F})$. If the same thing holds with
``$\leq$'' replaced by ``$<$'' then we say that $\mathcal{F}$ is
stable. Clearly, any line bundle is stable and any torsion sheaf is semistable.
The following facts can be found in \cite{GL}, Section~5.

\begin{prop}[\cite{GL}]\label{P:71}
The following set of assertions hold:
\begin{enumerate}
\item[i)] If $\mathcal{F}_1$ and $\mathcal{F}_2$ are semistable and
$\mu(\mathcal{F}_1) >\mu(\mathcal{F}_2)$
then $\mathrm{Hom}(\mathcal{F}_1,\mathcal{F}_2)=0$,
\item[i')]  If $\mathcal{F}_1$ and $\mathcal{F}_2$ are semistable and
$\mu(\mathcal{F}_1) \geq \mu(\mathcal{F}_2)$
then $\mathrm{Ext}(\mathcal{F}_2,\mathcal{F}_1)=0$,
\item[ii)]Any indecomposable
vector bundle is stable,
\item[iii)] The full subcategory $\mathcal{C}_{\mu}$ of $Coh_G(X)$
consisting of the zero sheaf together with all
semistable sheaves of slope $\mu$ is abelian, artinian and closed
under extensions. Its objects are all of finite length and the simple
objects are formed by the stable sheaves of
slope $\mu$.
\end{enumerate}
\end{prop}

The assertion i) is true for an arbitrary curve $X$; i') follows from
Serre duality together with the fact that
$|\omega_X| < 0$.

\vspace{.1in}

Any coherent sheaf $\mathcal{F}$ possesses a unique filtration
$0 \subset \mathcal{F}_1 \subset \cdots \subset
\mathcal{F}_r=\mathcal{F}$ such that $\mathcal{F}_i/\mathcal{F}_{i-1}$
is semistable of slope $\mu_i$ for $i=1, \ldots, r$, and
$\mu_1>\mu_2 >\cdots >\mu_r$. Furthermore, by
Proposition~\ref{P:71} i') this filtration splits (noncanonically),
i.e we have $\mathcal{F} \simeq \bigoplus_i \mathcal{F}_i/
\mathcal{F}_{i-1}$. We put $HN(\mathcal{F})=([\mathcal{F}_1],
[\mathcal{F}_2/\mathcal{F}_1], \cdots,
[\mathcal{F}/\mathcal{F}_{r-1}]) \in (K^+(X))^r$ and call this
sequence the \textit{HN type} of $\mathcal{F}$.
For any $\mu \in \mathbb{Q} \cup \{\infty\}$ will write
$\mathcal{F}_{\geq \mu}=
\mathcal{F}_i$ if $\mu_i \geq \mu$ but $\mu_{i+1} <\mu$, and we put
$\mathcal{F}_{< \mu}=\mathcal{F}/\mathcal{F}_{\geq \mu}$.

\vspace{.1in}

The following result is an easy consequence of the existence of the
HN filtration and of Proposition~\ref{P:71} i).

\begin{lem}\label{L:70}
If $Q_n^{\alpha} \neq 0$ then $\mu(\alpha) \geq n$.
\end{lem}

\vspace{.1in}

For any $(\alpha_1, \ldots, \alpha_r) \in (K^+(X))^r$ with
$\mu(\alpha_1) > \cdots >\mu(\alpha_r)$ and $\sum_i \alpha_i=\alpha$
and for any $n \in \Z$, let us
consider the subfunctor
of $\underline{\mathrm{Hilb}}^0_{\mathcal{E}_n^\alpha,\alpha}$ defined by
\begin{equation*}
\begin{split}
\Sigma \mapsto \{(\phi: \mathcal{E}^\alpha_n \otimes
\mathcal{O}_{\Sigma} \tto \mathcal{F}) \in
&\underline{\mathrm{Hilb}}^0_{\mathcal{E}_n^\alpha,\alpha}\;|\\
&HN(\mathcal{F}_{|\sigma})=(\alpha_1, \ldots,
\alpha_r)\;\text{for\;all\;closed\;points}\;\sigma \in \Sigma\}
\end{split}
\end{equation*}
This subfunctor is represented by a subscheme $HN_n^{-1}(\alpha_1,
\ldots, \alpha_r) \subset Q_n^\alpha$.
It is clear that for any $(\alpha_1, \ldots, \alpha_r)$ the subscheme
$HN_n^{-1}(\alpha_1, \ldots, \alpha_r)$
is $G_n^{\alpha}$-invariant. We will simply denote by
$Q_n^{(\alpha)}$ the subscheme $HN^{-1}(\alpha) \subset Q_n^{\alpha}$, which is
open (see, e.g \cite{LP}, Prop. 7.9).

\vspace{.1in}

\begin{lem}\label{L:71} The subscheme $HN_n^{-1}(\alpha_1, \ldots,
\alpha_r)$ is constructible, empty for all but finitely
many values of $(\alpha_1, \ldots, \alpha_r)$, and we have
$$Q_n^\alpha=\underset{(\alpha_1, \ldots, \alpha_r)}{\bigsqcup}
HN_n^{-1}(\alpha_1, \ldots, \alpha_r).$$
\end{lem}
\noindent \textit{Proof.} Let us fix some decomposition
$\mathcal{E}_n^\alpha\simeq \mathcal{E}_n^{\alpha_1} \oplus \cdots
\oplus \mathcal{E}_n^{\alpha_r}$ and denote by $\iota:
Q_n^{\alpha_1} \times \cdots \times Q_n^{\alpha_r} \to Q_n^\alpha$
the corresponding closed embedding. Then, since $Q_n^{(\alpha_i)}$
is open in $Q_n^{\alpha_i}$, $HN_n^{-1}(\alpha_1, \ldots,
\alpha_r)$ is open in the image of the natural map
$$G_n^{\alpha} \underset{G'}{\times} (Q_n^{\alpha_1} \times \cdots
\times Q_n^{\alpha_r})\to Q_n^{\alpha},$$ where $G' \simeq
G_n^{\alpha_1} \times \cdots \times G_n^{\alpha_r} \subset
G_n^{\alpha}$. In particular,
$HN_n^{-1}(\alpha_1, \ldots, \alpha_r)$ is also constructible.\\
\hbox to1em{\hfill} By Lemma~\ref{L:70}, $Q_n^{(\alpha_i)}$ is empty
if $\mu(\alpha_i) <n$, i.e if
$|\alpha_i| < n \cdot rank(\alpha_i)$.
For any fixed $n$ there exists only finitely many decompositions
$\alpha=\alpha_1 + \cdots + \alpha_r$ where all $\alpha_i$
satisfy the requirement $\mu(\alpha_i) \geq n$. This proves the second part of the Lemma.
The final statement is now obvious.\qed

\vspace{.1in}

\noindent \textbf{Exemple.} If $\a$ is the class of a line bundle
then the decomposition $Q_n^\a=\bigsqcup Q_n^\a[\mathcal{L}]$ in
Section~5.1. is the HN stratification with
$Q_n^\a[\mathcal{L}]=HN_n^{-1}(\a-[\mathcal{L}],[\mathcal{L}])$.

\vspace{.1in}

A priori, the subschemes $HN_n^{-1}(\alpha_1, \ldots, \alpha_r)$
only provide a finite stratification of $Q_n^\alpha$ in the broad
sense, i.e the Zariski closure of a stratum in general may not be
a union of strata. Hence we cannot define a partial order on the
set of possible HN types directly using inclusion of strata
closures. Instead we use the following combinatorial order: we say
that $(\alpha_1, \ldots, \alpha_r) \succ (\beta_1, \ldots,
\beta_s)$ if there exists $l$ such that $\alpha_1=\beta_1, \ldots,
\alpha_{l-1}=\beta_{l-1}$ and $\mu(\alpha_l) < \mu(\beta_l)$ or
$\mu(\alpha_l)=\mu(\beta_l)$ and $|\alpha_l| < |\beta_l|$. Note
that $HN(\mathcal{F}) \succ HN(\mathcal{G})$ if and only if there
exists $t \in \mathbb{R}$ such that $[\mathcal{F}_{\geq t}] >
[\mathcal{G}_{\geq t}]$ and $[\mathcal{F}_{\geq t'}]
=[\mathcal{G}_{\geq t'}]$ for any $t'>t$. Finally observe that
$Q_n^{\a} \backslash HN_n^{-1}(\alpha)=\bigcup_{\underline{\beta}
\prec \alpha} HN_n^{-1}(\underline{\beta})$.

\vspace{.1in}

\begin{lem}\label{L:74} The following hold~:
\begin{enumerate}
\item [i)] Assume that $HN(\mathcal{F}) \prec (\alpha_1,
\ldots, \alpha_r)$ and suppose that $\mathcal{F} \subset
\mathcal{G}$. Then $HN(\mathcal{G}) \prec (\alpha_1, \ldots,
\alpha_r)$.
\item[ii)] Let $\mathcal{F}$ be a coherent sheaf and let $0 \subset \mathcal{F}_1
\subset \ldots \subset \mathcal{F}_r=\mathcal{F}$ be a filtration such that
$\mu(\mathcal{F}_i/\mathcal{F}_{i-1}) \geq \mu(\mathcal{F}_{j}/\mathcal{F}_{j-1})$
if $i <j$. Then $HN(\mathcal{F}) \preceq (\mu(\mathcal{F}_1),\ldots,
\mu(\mathcal{F}/\mathcal{F}_{r-1}))$ with equality if and only if
$\mathcal{F}_i/\mathcal{F}_{i-1}$ is semistable for all $i$.
\end{enumerate}
\end{lem}
\noindent
\textit{Proof.} Statement i) is obvious. We prove statement ii). Let us set
$\mathcal{G}_i=\mathcal{F}_i/\mathcal{F}_{i-1}$. If $\mathcal{G}_i$ is semistable for all $i$
then it is clear
that $HN(\mathcal{F})=(\alpha_1, \ldots, \alpha_r)$.
Suppose on the contrary that
$\mathcal{G}_i$ is semistable for $i \leq j-1$ but that
$\mathcal{G}_j$ isn't. In particular, we have
$\mathcal{G}_{j,>\mu(\alpha_j)} \neq 0$. There is a unique sheaf
$\mathcal{G}'_j$ with $\mathcal{F}_{j-1} \subset
\mathcal{G}'_j \subset \mathcal{F}_j$ and
$\mathcal{G}'_j/\mathcal{F}_{j-1} \simeq
\mathcal{G}_{j, >\mu(\alpha_j)}$.
It is easy to see that $HN(\mathcal{G}'_j) \prec \underline{\alpha}$.
Hence by i) we also have
$HN(\mathcal{F}) \prec \underline{\alpha}$ as wanted.\qed

\vspace{.1in}

Now let $\mathbb{P}=(P_n)_n$ be a simple object in
$\mathbb{Q}^{\alpha}$. There exists $n \in \Z$ such that
$P_n \neq 0$ is a simple $G_n^{\alpha}$-equivariant perverse sheaf on
$Q_n^{\alpha}$. By Lemma~\ref{L:71} there exists
a unique HN type $(\alpha_1, \ldots, \alpha_r)$ with $\sum_i
\alpha_i=\alpha$ such that $supp(P_n) \subset
\overline{HN_n^{-1}(\alpha_1, \ldots, \alpha_r)}$ while $supp(P_n)
\cap HN_n^{-1}(\alpha_1, \ldots, \alpha_r)$
is dense in $supp(P_n)$. It is easy to see that the sequence
$(\alpha_1, \ldots, \alpha_r)$ is independent of $n$. We
call this the \textit{generic HN type} of $\mathbb{P}$ and denote it
by $HN(\mathbb{P})$.

\vspace{.1in}

Finally, we record the following result for future use.

\begin{lem}\label{L:72} Let $\mathcal{F}$ be a coherent sheaf in
$Coh_G(X)$, and let $\mu \in \mathbb{Q} \cup \{\infty\}$.
If $\mathcal{G} \subset \mathcal{F}$ is a subsheaf satisfying
$[\mathcal{G}]=[\mathcal{F}_{\geq \mu}]$ then
$\mathcal{G}=\mathcal{F}_{\geq \mu}$.\end{lem}
\noindent
\textit{Proof.} We argue by induction on the number of indecomposable
summands of $\mathcal{G}$ (for all $\mathcal{F}$ and
$\mu$ simultaneously). If $\mathcal{G}$ is indecomposable then by
Proposition~\ref{P:71}, ii) $\mathcal{G}$ is semistable
of slope $\mu(\mathcal{F}_{\geq \mu}) \geq \mu$. There is a
(noncanonical) splitting $\mathcal{F} \simeq
\mathcal{F}_{\geq \mu} \oplus \mathcal{F}_{<\mu}$. By the HN
filtration we have $\mathrm{Hom}(\mathcal{G},\mathcal{F}_{<\mu})=0$.
Hence $\mathcal{G} \subset \mathcal{F}_{\geq \mu}$. But
$[\mathcal{G}]=[\mathcal{F}_{\geq \mu}]$ and so $\mathcal{G}=
\mathcal{F}_{\geq \mu}$.\\
\hbox to1em{\hfill} Now assume that the statement is true for all
$\mathcal{F}'$ and $\mathcal{G}'$, where $\mathcal{G}'$ has
at most $k-1$ indecomposable summands, and let us assume that
$\mathcal{G}=\mathcal{G}_1 \oplus \cdots
\oplus \mathcal{G}_k$ with $\mathcal{G}_i$ indecomposable and
$\mu(\mathcal{G}_1) \geq \cdots \geq \mu(\mathcal{G}_k)$. We have
$\mu(\mathcal{G}_1) \geq \mu(\mathcal{G}) \geq \mu$, so arguing as
above we get $\mathcal{G}_1 \subset \mathcal{F}_{\geq \mu}$.
Now set $\mathcal{G}'=\mathcal{G}/\mathcal{G}_1$ and
$\mathcal{F}'=\mathcal{F}/\mathcal{G}_1$. We have
$[\mathcal{F}'_{\geq \mu}]=[\mathcal{F}_{\geq
\mu}/\mathcal{G}_1]=[\mathcal{F}_{\geq
\mu}]-[\mathcal{G}_1]=[\mathcal{G}]-
[\mathcal{G}_1]=[\mathcal{G}']$ and $\mathcal{G}'$ has $k-1$
indecomposable summands. Hence by the induction hypothesis,
$\mathcal{G}' \subset \mathcal{F}'_{\geq \mu}$, and so $\mathcal{G}
\subset \mathcal{F}_{\geq \mu}$, and finally
$\mathcal{G}=\mathcal{F}_{\geq \mu}$.\qed

\vspace{.2in}

\paragraph{\textbf{6.2. Relation to root systems.}} By the McKay correspondence,
$G \subset SO(3)$
is associated to a Dynkin diagram of type $ADE$. In that case, $\mathcal{L}\g$ is the affine
Lie algebra of the same type. Recall that by Lemma~\ref{L:00}, for
any pair $(X,G)$ there is a canonical isomorphism
$h: K(X) \simeq \widehat{Q}$ restricting to an isomorphism $K^+(X)
\simeq \widehat{Q}^+$. The following result is
proved in \cite{Le} (see also \cite{S1}, Section 7.)~:

\begin{lem}\label{L:7new1}
There is an indecomposable sheaf $\mathcal{F}$ of class $\alpha \in
K^+(X)$ if and only if $h(\alpha)\in \widehat{Q}^+$
is a root of $\mathcal{L}\g$. Moreover, such a sheaf is unique up to
isomorphism if and only if $h(\alpha)$
is a real root.\end{lem}

\vspace{.1in}

 From the explicit description of $h$ we easily deduce that
$\mathcal{C}_{\mu}$ has only finitely many simple (i.e, stable)
objects if $\mu < \infty$. In addition, from $|\omega_X| <0$ and
Serre duality  we conclude that for
any two semistable sheaves $\mathcal{F}_1$ and $\mathcal{F}_2$
of same slope $\mu \in \mathbb{Q}$,
$\mathrm{dim\;Ext}(\mathcal{F}_1,\mathcal{F}_2)=\mathrm{dim\;Hom}(\mathcal{F}_2,
\mathcal{F}_1(\omega_X))=0$. Hence, for any $\mu < \infty$,
$\mathcal{C}_{\mu}$ is a semisimple category with finitely many
simple objects.

\begin{lem}\label{L:73} The connected components of $Q_n^{(\alpha)}$
are in one to one correspondence with the set of isoclasses of
semistable sheaves of class $\alpha$ generated by
$\{\mathcal{L}_s(n)\}_{s \in S}$. In addition, each of these
connected component is a homogeneous
$G_n^{\alpha}$-variety.\end{lem} \noindent \textit{Proof.} Let
$\mathcal{G}_1, \ldots, \mathcal{G}_k$ be a complete collection of
distinct simple objects in $\mathcal{C}_{\mu}$. Since
$\mathcal{C}_{\mu}$ is semisimple, we have  $\mathcal{F} \simeq
\bigoplus_{l} \mathcal{G}_l \otimes \mathrm{Hom}(\mathcal{G}_l,
\mathcal{F})$ for any $\mathcal{F} \in \mathcal{C}_{\mu}$, and
\begin{equation}\label{E:77}
\sum_l \mathrm{dim\;Hom}(\mathcal{G}_l, \mathcal{F}) \cdot
[\mathcal{G}_l]=\alpha.
\end{equation}
Consider the tautological sheaf $\underline{\mathcal{F}}$ on the
universal family $\mathcal{X}=X \times Q_n^{\alpha}$ and its
restriction to $\mathcal{X}'=X \times Q_n^{(\alpha)}$. By
definition we have $\underline{\mathcal{F}}_{|X \times t} \simeq
\mathcal{G}$ if $t=(\phi: \mathcal{E}_n^{\alpha} \tto
\mathcal{G})$. By the semicontinuity theorem (\cite{Ha}, III Thm.
12.8) the function $\mathrm{dim\;Hom}(\mathcal{G}_l,
\underline{\mathcal{F}}_{|X \times t})$ is upper semicontinuous in
$t$. Using (\ref{E:77}) we see that it is also lower
semicontinuous in $t$ for $t \in Q_n^{(\alpha)}$, and thus locally
constant on $Q_n^{(\alpha)}$. We deduce that
$\underline{\mathcal{F}}_{|X \times t} \simeq
\underline{\mathcal{F}}_{|X \times t'}$ if $t$ and $t'$ belong to
the same connected component of $Q_n^{(\alpha)}$. Finally, using
Lemma~\ref{L:3} we conclude that each connected component of
$Q_n^{(\alpha)}$ is formed by a single $G_n^{\alpha}$-orbit. \qed

\vspace{.1in}

\paragraph{} The next remark will be of importance for us.
\begin{lem}\label{L:fin71}
Let $\mathcal{H}$ be an indecomposable vector bundle, let $d \in
\N$ and set $\mathcal{F}=\mathcal{H}^{\oplus d}$. There exists
sheaves $\mathcal{G}$ and $\mathcal{K}$ with
$rank(\mathcal{G})<rank(\mathcal{F})$ or $rank(\mathcal{G})
=rank(\mathcal{F})$ and $|\mathcal{G}|<|\mathcal{F}|$, satisfying
the following property: there exists a unique subsheaf
$\mathcal{G}' \subset \mathcal{F}$ such that $\mathcal{G}' \simeq
\mathcal{G}$ and $\mathcal{F}/\mathcal{G}' \simeq \mathcal{K}$.
Moreover, if $\mathcal{H}$ is generated by $\{\mathcal{L}_s(n)\}$
and $\mathcal{H} \neq \mathcal{O}_X(n)$ then $\mathcal{G}$ and
$\mathcal{K}$ can also be chosen to be generated by
$\{\mathcal{L}_s(n)\}$.
\end{lem}
\noindent
\textit{Proof.} The statement is invariant under twisting by a line
bundle. If $\mathcal{H}=\mathcal{L}_i$ then we have $\mathrm{Hom}(\mathcal{O}_X,\mathcal{L}_i)=1$
and we may take $\mathcal{G}=\mathcal{O}_X^{\oplus d}$ and $\mathcal{K}=(S_i^{(0)})^{\oplus d}$. Thus the
Lemma is true for all line bundles. Now let us assume that
$\mathcal{H}$ is of rank at least two. By inspection of all the cases,
using Lemma~\ref{L:7new1} and the definition of $\widehat{Q}$, one
checks that, up to twisting by a line bundle, the class of
$\mathcal{H}$ can be written as
$$[\mathcal{H}]=rank(\mathcal{H})[\mathcal{O}_X] + \sum_{i,j} n_i^{(j)}
[\a_i^{(j)}],$$
with $n_i^{(j)} \geq 0$ and $n_i^{(0)}=0$ for all $i$. Note that
$[\mathcal{H}] \neq rank(\mathcal{H})[\mathcal{O}_X]$ since
$\mathcal{H}$ is
indecomposable, hence
$\mu(\mathcal{H})>0$. Therefore, using \cite{S1}, Proposition~5.1 we
have $\mathrm{dim\;Hom}(\mathcal{O}_X,\mathcal{H})=
\langle [\mathcal{O}_X],[\mathcal{H}]\rangle=rank(\mathcal{H})$. In
particular, the natural map $\phi:
\mathrm{Hom}(\mathcal{O}_X,\mathcal{H}) \otimes \mathcal{O}_X \to
\mathcal{H}$ is neither zero nor surjective. The same is true of the map
$\phi^{\oplus d}:\mathrm{Hom}(\mathcal{O},\mathcal{F}) \otimes \mathcal{O}_X \to
\mathcal{F}$. We claim that
$\mathcal{G}=\mathrm{Im}\;\phi^{\oplus d}$ and $\mathcal{K}=\mathrm{Coker}\;
\phi^{\oplus d}$ satisfy the requirements of the Lemma. Indeed, let
$$0 \to \mathcal{G} \stackrel{a}{\longrightarrow} \mathcal{F}
\stackrel{b}{\longrightarrow} \mathcal{K} \to 0$$
be an exact sequence. The map $a$ induces a linear map $c:
\mathrm{Hom}(\mathcal{O}_X,\mathcal{F}) \to \mathrm{Hom}(\mathcal{O}_X,
\mathcal{F})$ such that the diagram
$$\xymatrix{
\mathcal{G} \ar[r]^-{a} & \mathcal{F}\\
\mathrm{Hom}(\mathcal{O}_X,\mathcal{F}) \otimes \mathcal{O}_X
\ar[u]^-{can} \ar[r]^-{c}& \mathrm{Hom}(\mathcal{O}_X,\mathcal{F})
\otimes \mathcal{O}_X \ar[u]_-{\phi^{\oplus d}}}$$ is commutative,
and thus $\mathrm{Im}\;a \subseteq \mathrm{Im}\;\phi^{\oplus d}$.
But from $\mathrm{Coker}\;a \simeq \mathrm{Coker}\;\phi^{\oplus
d}$ we deduce that $\mathrm{Im}\;a = \mathrm{Im}\;\phi^{\oplus d}$
as wanted. The last assertion of the Lemma is clear from the above
construction. \qed

\vspace{.2in}

\section{Some Lemmas.}

\vspace{.1in}

\paragraph{} This section contains several technical results
needed later in the course of the proof of Theorem 5.1 ii).

\vspace{.1in}

\paragraph{\textbf{7.1.}} We begin with the following

\begin{lem}\label{L:after81} Let $\mathbb{V} \in \mathbb{U}^{\alpha_r} \hat{\boxtimes} \cdots
\hat{\boxtimes} \mathbb{U}^{\alpha_1}$, where $\mu(\alpha_1) \geq \cdots \geq \mu(\alpha_r)$, and
let us put $\underline{\alpha}=(\alpha_r, \ldots, \alpha_1)$.
Then
\begin{enumerate}
\item[i)] we have $supp(\mathrm{Ind}^{\underline{\alpha}}(\mathbb{V}))
\subset \bigcup_{\underline{\beta} \preceq \underline{\alpha}} HN^{-1}(\underline{\beta})$.
\item[ii)] If in addition $V_{-1}=0$ and $\mu(\alpha_i) \geq 0$ for all $i$ then
 $supp(\mathrm{Ind}^{\underline{\alpha}}(\mathbb{V}))
\subset \bigcup_{\underline{\beta} \prec \underline{\alpha}} HN^{-1}(\underline{\beta})$.
\end{enumerate}
\end{lem}
\noindent
\textit{Proof.} Consider the (iterated) induction diagram (see Section~4.1. for notations)
$$
\xymatrix{
Q_l^{\alpha_r} \times \cdots \times Q_l^{\alpha_1} & E' \ar[l]_-{p_1}
\ar[r]^-{p_2} & E'' \ar[r]^-{p_3} &
Q_l^{\alpha}}
$$
It is enough to show that $p_3(E'') \subset \bigcup_{\underline{\beta} \preceq \underline{\alpha}}
HN^{-1}(\underline{\beta})$. To see this, let $(\mathcal{G}_r,
\ldots, \mathcal{G}_1) \in
Q_l^{\alpha_r} \times \cdots \times Q_l^{\alpha_1}$ and
$(\underline{V}, \phi: \mathcal{E}_l^{\a} \tto \mathcal{F})
\in p_2p_1^{-1}(\mathcal{G}_r, \ldots, \mathcal{G}_1)$. Putting $\mathcal{F}_i=\phi(\underline{V}_i)$
we have $\mathcal{F}_1 \subset \cdots \subset \mathcal{F}_r=\mathcal{F}$ and
$\mu(\mathcal{F}_i/\mathcal{F}_{i-1}) \geq \mu(\mathcal{F}_{i+1}/\mathcal{F}_{i})$. Thus, from
Lemma~\ref{L:74} ii) we obtain $HN(\mathcal{F}) \preceq (\underline{\alpha})$. This proves i).
To prove ii), first note that by \cite{GL} Cor. 3.4 and
\cite{GL} Prop. 4.1, a sheaf $\mathcal{F}$ is generated by $\{\mathcal{L}_s(-1)\}$ if and only if
$\mathrm{Ext}(\bigoplus_{s} \mathcal{L}_s(-1) \oplus \mathcal{O}_X,\mathcal{F})=0$. In particular,
by Proposition~\ref{P:71} i'), any semistable sheaf of slope $\mu \geq 0$ is generated by
$\{\mathcal{L}_s(-1)\}$. From $\mu(\alpha_i) \geq 0$ and $V_{-1}=0$ we deduce that
$$supp(V_l) \cap Q_l^{(\alpha_r)} \times \cdots \times Q_l^{(\alpha_1)} =\emptyset$$
for any $l \in \Z$. But by Lemma~\ref{L:74} ii) again, we have
$p_2^{-1}p_{3}^{-1}(HN^{-1}(\underline{\alpha}))=p_1^{-1}(Q_l^{(\alpha_r)} \times \cdots \times
Q_l^{(\alpha_1)})$. It follows that $supp(\mathrm{Ind}^{\underline{\alpha}}(\mathbb{V})) \cap
HN^{-1}(\underline{\alpha})=\emptyset$, as wanted.\qed

\vspace{.2in}

\paragraph{\textbf{7.2.}} Let us now consider the following situation.
Let $\alpha,\beta, \gamma \in K^+(X)$ such that
$\alpha=\beta+\gamma$, and let $U_n^{\alpha} \subset
Q_n^{\a},U_n^{\b} \subset Q_n^{\b}$
and $U_n^{\gamma} \subset Q_n^{\gamma}$ be locally closed subsets
invariant under the corresponding group, and satisfying the following
properties~:
\begin{enumerate}
\item[i)] Denoting by $\xymatrix{ Q_n^{\alpha} & F \ar[l]_-{\iota}
\ar[r]^-{\kappa} & Q_n^{\gamma} \times Q_n^{\beta}}$ the restriction
diagram as in Section~4.2, we have $U_n^\a \cap F =\kappa^{-1}(U_n^\gamma
\times U_n^{\b})$,
\item[ii)] For any closed point $(\phi: \mathcal{E}_n^\a \tto
\mathcal{F}) \in U_n^{\alpha}$ there exists a unique subsheaf
$\mathcal{G}$
of $\mathcal{F}$ of class $\beta$, and we have
$\mathrm{Ext}(\mathcal{F}/\mathcal{G},\mathcal{G})=0$.
\end{enumerate}
The set $U^{\alpha}_n$ gives rise, via the induction functors
$\Xi_{m,n}$ and $\Xi_{n,m}$, to locally closed sets $U^\alpha_m$ for
any $m \in
\Z$. The same is true for $U^\b_n$ and $U^\gamma_n$, and properties i) and ii)
above are satisfied for all $m$.
Let us denote by $j_m: U_m^\a \hookrightarrow Q_m^\a$ the embedding.
\begin{lem}\label{L:7new5} Under the above conditions, we have, for
any $\mathbb{P}=(P_m)_m \in \mathbb{Q}^\alpha$ such that
$supp(P_m) \subset \overline{U_m^{\alpha}}$ for all $m$,
$$j^*(\mathrm{Ind}^{\gamma,\b} \circ
\mathrm{Res}^{\gamma,\b}(\mathbb{P})) \simeq j^*(\mathbb{P}).$$
\end{lem}
\noindent
\textit{Proof.} Consider the following commutative diagram, where we
have set $F^\vee=F \cap U_l^{\a}$, where the vertical
arrows are canonical embeddings and where the notations are otherwise
as in Section~4.2.~:
$$\xymatrix{
U_l^{\alpha} \ar[d]_-{g_1} & F^\vee \ar[l]_-{\iota'}
\ar[r]^-{\kappa'} \ar[d]_-{g_2} &
U_l^\gamma \times U_l^{\b} \ar[d]_-{g_3}\\
Q_l^\a & F \ar[l]_-{\iota} \ar[r]^-{\kappa} & Q_l^{\gamma} \times Q_l^{\b}}$$

By the assumption i), the rightmost square is Cartesian. We deduce
by base change that
\begin{equation}\label{E:71}
g_3^*\kappa_{!}\iota^* = \kappa'_{!}{\iota'}^*g_1^*.
\end{equation}

Next look at the diagram
$$\xymatrix{
U_l^{\a} \ar[d]_-{g_1} & {E''}^\vee \ar[l]_-{p'_3} \ar[d]_-{g_4}&
{E'}^\vee \ar[l]_-{p'_2} \ar[r]^-{p'_1} \ar[d]_-{g_5} & U^\gamma_l\times
U_l^\beta \ar[d]_-{g_3}\\
Q_l^\alpha & E'' \ar[l]_-{p_3} & E' \ar[l]_-{p_2} \ar[r]^-{p_1} &
Q_l^{\gamma} \times Q_l^{\beta}}$$
where ${E''}^\vee=p_3^{-1}(U_l^{\alpha})$ and
${E'}^\vee=p_2^{-1}({E''}^\vee)$. By assumption $i)$ the rightmost
square is Cartesian,
while by assumption ii)  $p'_3$ is an isomorphism . Hence, by base change
we get
\begin{equation}\label{E:72}
g_1^* p_{3!}p_{2\flat}p_1^* \simeq p'_{2\flat}{p'_1}^* g_3^*.
\end{equation}
Now let us set $\mathbb{R}=\mathrm{Ind}^{\gamma,\beta} \circ
\mathrm{Res}^{\gamma,\beta}(\mathbb{P})$.
Combining (\ref{E:71}) with (\ref{E:72}) and applying this to
$\mathbb{R}$ yields
\begin{equation}\label{E:73}
\begin{split}
g_1^*(R_m)=&g_1^* \big(\widetilde{\mathrm{Ind}}^{\underline{\alpha}}
\circ \widetilde{\mathrm{Res}}^{\underline{\alpha}}
(\mathbb{P})\big)_m[d_m]\\
=&\Xi_{l,m} p'_{2\flat}{p'_1}^*\kappa'_{!}{\iota'}^*g_1^*(P_l)[d_m]
\end{split}
\end{equation}
for any $ l \ll m$, where $d_m=\mathrm{dim}(G_m^\a)-\mathrm{dim}(
G_m^{\gamma}\times G_m^\b)$.

\vspace{.05in}

Recall that a choice of a subsheaf
$\mathcal{E}_l^{\beta} \subset \mathcal{E}_l^{\alpha}$
is implicit in the definition of $F$. Let us denote by $P \subset
G_l^{\a}$ its stabilizer. Let us also fix a splitting
$\mathcal{E}_l^{\alpha}=\mathcal{E}_l^{\beta} \oplus \mathcal{E}_l^{\gamma}$.
This choice defines a Levi subgroup $L \subset P$ isomorphic to
$G'= G_l^{\gamma}\times G_l^{\b}$, and induces a section $s:
Q_l^{\gamma} \times Q_l^{\b} \to F$ of $\kappa$.

\vspace{.05in}

We claim that $F^\vee \simeq P \underset{L}{\times} (U_l^\gamma \times
U_l^{\b})$.
Indeed, by assumption ii), if
$$\kappa((\phi: \mathcal{E}^\a_l \tto \mathcal{G}))=(
(\phi_2: \mathcal{E}^\gamma_l \tto \mathcal{F}_2),(\phi_1:
\mathcal{E}^\beta_l \tto\mathcal{F}_1)) \in U_l^{\gamma}
\times U_l^{\beta}$$
  then $\mathcal{G} \simeq
\mathcal{F}_1
\oplus \mathcal{F}_2$. Hence, using Lemma~\ref{L:3}, we have
\begin{equation*}
\begin{split}
F^\vee&= \bigg(G_l^{\alpha} \underset{L}{\times} s(U_l^{\gamma} \times
U_l^{\beta})\bigg)
\cap F\\
&=P \underset{L}{\times} s(U_l^{\gamma} \times U_l^{\beta}).
\end{split}
\end{equation*}
Similarly, it is easy to check that $U^\alpha_l \simeq G_l^\a
\underset{P}{\times} F^\vee$
and that ${E'}^\vee$ is a principal $G'$-bundle over $U_l^{\alpha}$.
To sum up, we have a commutative
diagram
\begin{equation}\label{D:71}
\xymatrix{
U_l^{\alpha} & F^\vee \ar[l]_-{\iota'} \ar[d]_-{\kappa'}\\
{E'}^\vee \ar[u]^-{p'_2} \ar[r]^-{p'_1} & U_l^{\gamma} \times U_l^{\beta}}
\end{equation}
and equivalences
\begin{align*}
{\iota'}^*:~& \mathcal{Q}_{G_l^{\a}}(U_l^\alpha) \stackrel{\sim}{\to}
\mathcal{Q}_P(F^\vee)[\mathrm{dim}(G_l^{\a}/P)],\\
{p'_2}^*:~&\mathcal{Q}_{G_l^{\a}}(U_l^\alpha) \stackrel{\sim}{\to}
\mathcal{Q}_{G' \times G_l^{\a}}({E'}^\vee)[-\mathrm{dim}(L)],\\
{\kappa'}^*:~& \mathcal{Q}_{G'}(U_l^{\gamma}
\times U_l^{\beta}) \stackrel{\sim}{\to}
\mathcal{Q}_P(F^\vee)[-\mathrm{dim}(P/L)],\\
{p'_1}^*:~& \mathcal{Q}_{G'}(U_l^{\gamma}
\times U_l^{\beta}) \stackrel{\sim}{\to} \mathcal{Q}_{G' \times
G_l^\a}({E'}^\vee)[-\mathrm{dim}(G_l^\a)].
\end{align*}
Let $\mathbb{T}$ be any object in $\mathcal{Q}_{G_l^\a}(U_l^\alpha)$.
By the above set of equivalences,
there exists $\mathbb{T}' \in \mathcal{Q}_{G'}(U_l^\gamma \times U_l^{\beta})
[\mathrm{dim}(G_l^\a/L)]$
such that ${\iota'}^*\mathbb{T} \simeq {\kappa'}^* \mathbb{T}'$.
Hence $\kappa'_! {\iota'}^* \mathbb{T}=\kappa'_!{\kappa'}^*
\mathbb{T}'=\mathbb{T}'[-2\mathrm{dim}(P/L)]$. Moreover, ${p'_1}^*
\mathbf{T}'={p'_2}^* \mathbb{T}$
so that $p'_{2\flat}{p'_1}^* \mathbb{T}'
\simeq \mathbb{T}$. All together we obtain $p'_{2\flat}{p'_1}^*
\kappa'_! {\iota'}^* \mathbb{T} \simeq \mathbb{T}
[-2\mathrm{dim}(P/L)]$. Applying
this to $g_1^*(P_l)$ yields the desired result. The Lemma is proved. \qed

\vspace{.1in}

\begin{cor}\label{C:72}
Let $\alpha,\alpha_1,\ldots,\alpha_r \in K^+(X)$ such that
$\alpha=\sum_i \alpha_i$. Assume given locally closed
subsets $U_n^\alpha \subset Q_n^\alpha$, $U_n^{\alpha_i} \subset
Q_n^{\alpha_i}$ satisfying the following conditions~:
\begin{enumerate}
\item[i)] Denoting by $\xymatrix{ Q_n^{\alpha} &
F_{\underline{\alpha}} \ar[l]_-{\iota} \ar[r]^-{\kappa} &
Q_n^{\alpha_r} \times
\cdots \times
Q_n^{\alpha_1}}$ the iterated restriction diagram, we have $U_n^\a
\cap F_{\underline{\alpha}} =\kappa^{-1}(U_n^{\alpha_r} \times
\cdots \times U_n^{\alpha_1})$,
\item[ii)] For any closed point $(\phi: \mathcal{E}_n^\a \tto
\mathcal{F}) \in U_n^{\alpha}$ there exists a unique filtration
$$\mathcal{G}_1 \subset \cdots \subset \mathcal{G}_r=\mathcal{F}$$
of $\mathcal{F}$ such that
$[\mathcal{G}_i/\mathcal{G}_{i-1}]=\alpha_i$, and we have
$\mathrm{Ext}(\mathcal{G}_j/\mathcal{G}_{j-1},\mathcal{G}_i/\mathcal{G}_{i-1})=0$
for any $j >i$.
\end{enumerate}
Then, for any $\mathbb{P}=(P_m)_m \in \mathbb{Q}^\alpha$ such that
$supp(P_m) \subset \overline{U_m^{\alpha}}$ and $supp(P_m) \cap
U_m^{\alpha} \neq \emptyset$ for any $m$, we have
$$j^*(\mathrm{Ind}^{\alpha_r,\ldots,\alpha_1} \circ
\mathrm{Res}^{\alpha_r,\ldots,\alpha_1}(\mathbb{P})) \simeq
j^*(\mathbb{P}),$$
where $j_l: U_l^{\alpha} \to Q_l^{\alpha}$ is the embedding.
\end{cor}
\noindent
\textit{Proof.} Use induction and Lemma~\ref{L:7new5}.\qed

\vspace{.2in}

\section{Proof of Theorem~\ref{T:5} $i)$}

\vspace{.1in}

\paragraph{\textbf{8.1.}}
Let $\mathcal{P}^\alpha_{\geq d}$ be the set of perverse sheaves
of $\mathcal{P}^\alpha$ with support contained in
$U^{\alpha}_{\geq d}$, and set
$\mathcal{P}^\alpha_d=\mathcal{P}^\alpha_{\geq d} \backslash
\mathcal{P}^\alpha_{\geq d-1}$. For simplicity, let us temporarily
denote by $\mathfrak{W}$ the subalgebra of
${\mathfrak{U}}_{\mathbb{A}}^{tor}$ generated by
$\mathbf{1}_{l\alpha_{(i,j)}}$ and
$\mathbf{1}_{l\delta}$ for $(i,j) \in \aleph$ and $l \in \N$. We will prove by induction
on $d$ that for any torsion class $\alpha$, any element
$\mathbf{b}_{\mathbb{P}}$ with $\mathbb{P} \in
\mathcal{P}^\alpha_d$ belongs to $\mathfrak{W}$. The following two
extreme cases will serve as a basis for our induction.

\begin{lem}\label{L:16} Any $\mathbb{P} \in\mathcal{P}^{\alpha}_0$ is
isomorphic to a
sheaf of the form $\mathbf{IC}(U^{l\delta}_{((1), \ldots,
(1))},\sigma)$ for some irreducible representation
$\sigma$ of $\mathfrak{S}_l$. \end{lem}
\noindent
\textit{Proof.} Since $U^{\alpha}_{\geq 0}=U^{\alpha}_{\geq 1}$ for
$\alpha\not\in \N\delta$ we may assume that $\alpha=l\delta$ for some
$l$. Next, it is clear that any perverse sheaf in $\mathcal{P}^\alpha$
appearing in $\mathrm{Ind}^{\alpha_r, \ldots, \alpha_1}( \mathbf{1}_{\alpha_r}
\boxtimes \cdots \boxtimes
\mathbf{1}_{\alpha_1})$ with $\alpha_k \in
\{l\alpha_{(i,j)}\}_{(i,j)\in \aleph}$ for at least one $k$ actually lies in
$\mathcal{P}^\alpha_{\geq 1}$. Therefore $\mathbb{P}$ appears in a
product
$A:=\widetilde{\mathrm{Ind}}^{\delta,\ldots,\delta}(\mathbf{1}_{\delta}
\boxtimes
\cdots \boxtimes \mathbf{1}_{\delta})$. By definition, we have
$A=p_{3!}(\mathbf{1}_{E''})$
where
\begin{equation*}
E''=\{(\phi: \mathcal{E}^{l\delta} \tto \mathcal{F}, \underline{V}_1
\subset \cdots \subset \underline{V}_{l}=\mathcal{E}^{l\delta}\;|\;
\underline{V}_i/\underline{V}_{i-1} \simeq \mathcal{E}^{\delta},
[\phi(\underline{V}_i)]=i\delta\}
\end{equation*}
and $p_3$ is the natural projection to $Q^{l\delta}$. Since by
assumption $U^{l\delta} \cap supp(\mathbb{P})$ is dense in
$supp(\mathbb{P})$ we
may restrict ourselves to the open set $p_3^{-1}(U^{l\delta})$.
We claim that
the map $p_3$ is semismall, with $U^{l\delta}_{((1),\ldots,(1))}$ being
the unique relevant strata. Indeed, let
$\underline{\lambda}=(\lambda^{(1)},\ldots, \lambda^{(r)})$ be
an $r$-tuple of partitions such that $\sum_{k} |\lambda^{(k)}|=l$. For any
point $\phi: \mathcal{E}^{l\delta} \tto \mathcal{F}$ in
$U^{l\delta}_{\underline{\lambda}}$ the sheaf $\mathcal{F}$ decomposes
as a direct sum $\mathcal{F}=\sum_k \mathcal{F}_k$ of torsion sheaves
supported at distinct points, and by Lemma~\ref{L:3} the dimension of
the corresponding
$G^{l\delta}$-orbit is equal to
$|\Gamma|l^2-\mathrm{dim}\;\mathrm{Aut}\;\mathcal{F}= |\Gamma|l^2 -
\sum_k \mathrm{dim\;Aut}\;\mathcal{F}_k$. Now, we may view a torsion
sheaf supported at a generic point as an element of a space
$\mathcal{N}^{(1)}=\mathcal{N}$ of nilpotent representations of the
quiver with one loop (or nilpotent cone). In particular, we have
\begin{equation*}
\begin{split}
\mathrm{dim\;Aut}\;\mathcal{F}_k&=\big(\sum_j
\lambda^{(k)}_j\big)^2-\mathrm{dim}\;\mathcal{N}_{\lambda^{(k)}}\\
&=\big(\sum_j
\lambda^{(k)}_j\big)^2-(\mathrm{dim}\;\mathcal{N}_{\sum_j\lambda^{(k)}_j}-
2\mathrm{dim}\;\mathcal{B}_{\lambda^{(k)}})\\
&=2\mathrm{dim}\;\mathcal{B}_{\lambda^{(k)}}+\sum_j\lambda^{(k)}_j,
\end{split}
\end{equation*}
where $\mathcal{N}_{\lambda^{(k)}} \subset
\mathcal{N}_{\sum_j\lambda^{(k)}_j}$ is the nilpotent orbit of
$\mathfrak{gl}_{\sum_j\lambda^{(k)}_j}$
associated to $\lambda^{(k)}$, and $\mathcal{B}_{\lambda^{(k)}}$ is the
Springer fiber over that orbit.
Since the isomorphism
classes of sheaves $\mathcal{F}$ of the above type forms a smooth
$k$-dimensional family, we obtain
$\mathrm{dim}\;U^{l\delta}_{\underline{\lambda}}=
|\Gamma|l^2 +r-\sum_k \mathrm{dim\;Aut}\;\mathcal{F}_k=|\Gamma|l^2 +
r-l-2\mathrm{dim}\; \mathcal{B}_{\lambda^{(k)}}$ and
\begin{equation}\label{E:61}
\mathrm{codim}\;U^{l\delta}_{\underline{\lambda}}=(l-r)+2\sum_k\mathrm{dim}\;\mathcal{B}_{\lambda^(k)}.
\end{equation}
Finally, there is a canonical finite morphism from the fiber
$p_3^{-1}(u)$ of any point $u\in U^{l\delta}_{\underline{\lambda}}$ to
$\prod_k \mathcal{B}_{\lambda^{(k)}}$, so that
$\mathrm{dim}\;p_3^{-1}(u)=
\sum_k\mathrm{dim}\;\mathcal{B}_{\lambda^(k)}$. The claim now follows
by comparing with (\ref{E:61}).\\
 From the Decomposition theorem of \cite{BBD} (see also \cite{BM}) we
thus obtain
\begin{equation}\label{E:62}
j^*A=\mathbf{IC}(U^{l\delta}_{((1),\ldots,(1))},V)=\bigoplus_{V_i \in Irrep
\;\mathfrak{S}_l} \mathbf{IC}(U^{l\delta}_{((1),\ldots,(1))},V_i)
\otimes V_i^*,
\end{equation}
where $j: U^{l\delta} \to Q^{l\delta}$ is the open embedding, and
$Irrep\; \mathfrak{S}_l$ denotes the set of irreducible
representations of $\mathfrak{S}_l$. The Lemma follows.~\qed

\vspace{.1in}

The above proof implies the following generalization of (\ref{E:62}).
\begin{cor} Fix integers $n_1, \ldots, n_r$ such that $\sum n_i=l$
and let $\sigma_i$ be a
representation of $\mathfrak{S}_{n_i}$. We have
\begin{equation}\label{E:625}
\begin{split}
j^*\widetilde{\mathrm{Ind}}^{n_r\delta, \cdots
,n_1\delta}(\mathbf{IC}&(U^{n_r\delta}_{((1),\ldots,(1))},\sigma_r)
\boxtimes \cdots \boxtimes
\mathbf{IC}(U^{n_1\delta}_{((1),\ldots,(1))},\sigma_1))\\
&=\mathbf{IC}(U^{l\delta}_{((1),\ldots,(1))},ind_{\mathfrak{S}_{n_r}\times
\cdots \times
\mathfrak{S}_{n_1}}^{\mathfrak{S}_l}(\sigma_r \otimes \cdots \otimes
\sigma_1)).
\end{split}
\end{equation}
\end{cor}

\vspace{.1in}

\begin{lem}\label{L:17} If $\mathbb{P} \in \mathcal{P}^\gamma_{|\gamma|}$ then
   $\mathbf{b}_{\mathbb{P}} \in \mathfrak{W}$.\end{lem}
\noindent
\textit{Proof.} By assumption we have
$$supp(\mathbb{P}) \subset Y_{\gamma}:=\{\phi:
\mathcal{E}^{\gamma} \tto \mathcal{F}\;|\; supp(\mathcal{F}) \subset
\Lambda\}.$$
Note that $Y_{\gamma}$ is a disjoint union of varieties $Y_{\gamma_1,
   \ldots, \gamma_N}$ with $\sum \gamma_i=\gamma$, where
$$Y_{\gamma_1, \ldots, \gamma_N}:=\{\phi: \mathcal{E}^{\gamma} \tto
\mathcal{F}\;|\; \mathcal{F}=\bigoplus_i \mathcal{F}_i,\;
supp(\mathcal{F}_i)=\lambda_i, [\mathcal{F}_i]=\gamma_i\}.$$
Since $\mathbb{P}$ is simple, we have $supp(\mathbb{P}) \subset
Y_{\gamma_1,\ldots,\gamma_N}$ for
some fixed $\gamma_1, \ldots, \gamma_N$.

Let us first assume that there exists $i$ such that $\gamma_i=\gamma$
(and $\gamma_j=0$ for $j \neq i$). Then
\begin{equation}\label{E:ff1}
Y_{(\gamma_j)} \simeq
G^{\gamma} \underset{\phi^i_{\gamma}(G_{\gamma})}{\times}
\mathcal{N}^{(p_i)}_\gamma.
\end{equation}
It is known that this space has finitely many $G^\gamma$-orbits and
that they are all simply connected.
Set
$$\mathfrak{U}_i^{tor}=\bigoplus_{\mathbb{P} \in \mathcal{P}(i)} \mathbb{A}\mathbf{b}_{\mathbb{P}},
\qquad
\text{where}\; \mathcal{P}(i)=\{\mathbb{P} \in \bigcup_{\gamma} \mathcal{P}^\gamma_{|\gamma|}\;|\;
supp(\mathbb{P})\subset Y_{\gamma_1, \ldots, \gamma_N}, \gamma_i=\gamma\}.$$
Following Lusztig (see \cite{L3}), consider the $\C(v)$-linear map to the Hall algebra
\begin{align*}
\tau_i:\;\mathfrak{U}_i^{tor} &\to \mathbf{H}_{p_i}\\
\mathbf{b}_{\mathbb{P}} \in \mathfrak{U}_i^{tor}[\gamma] &\mapsto
\big( x \mapsto
v^{\mathrm{dim}(G_{\gamma})-\mathrm{dim}(G^\gamma)}\sum_i
\mathrm{dim}_{\overline{\mathbb{Q}}_l}
\mathcal{H}^i_{|x}(\mathbb{P})v^{-i}\big),
\end{align*}
which is a (bi)algebra isomorphism (in the above, we view $\mathcal{N}^{(p_i)}_{\gamma}$ as a subvariety
of $Y_{(\gamma_j)}$ via (\ref{E:ff1})).
Furthermore, $\mathbf{H}_{p_i}$ is equipped with a canonical basis $\mathbf{B}_{\mathbf{H}_{p_i}}=
\{\mathbf{b}_{\mathbf{m}}\}$, which consists of the functions
$$x \mapsto \sum_i \mathrm{dim}_{\overline{\mathbb{Q}}_l}
\mathcal{H}^i_{|x}(\mathbf{IC}(O))v^{-i}$$ where $O$ runs through
the set of all $G_{\gamma}$-orbits in
$\mathcal{N}^{(p_i)}_{\gamma}$, for all $\gamma$. In particular,
up to a power of $v$, we have $\tau_i(\mathbf{b}_{\mathbb{P}}) \in
\mathbf{B}_{\mathbf{H}_{p_i}}$ for any $\mathbb{P} \in
\mathcal{P}(i)$.

We now prove by induction on $|\gamma|$ that
$\mathbf{b}_{\mathbb{P}} \in \mathfrak{W}$. If $\gamma \not>
\delta$ then $\tau_i(\mathbf{b}_{\mathbb{P}}) \in
\mathbf{U}_v^+(\widehat{\mathfrak{sl}}_{p_i}) \subset
\mathbf{H}_{p_i}$. If $\gamma=\delta$ then by definition either
$\mathbb{P}=\mathbf{1}_{\delta}$ or $\mathbb{P}$ appears in some
induction product of the sheaves $\mathbf{1}_{\alpha_{(i,j)}}$.
The first case is ruled out since $\mathbf{1}_{\delta} \not\in
\mathcal{P}^\delta_{|\delta|}$, and in the second case
$\tau_i(\mathbf{b}_{\mathbb{P}}) \in
\mathbf{U}_v^+(\widehat{\mathfrak{sl}}_{p_i})$ by \cite{L1},
Section~12.3. Now assume that the result is proved for all
$\gamma'$ with $|\gamma'| < |\gamma|$ and that $\gamma > \delta$.
By Lemma~\ref{L:14}, we have $\mathrm{Res}^{\b,\a}(\mathbb{P}) \in
\mathbb{U}^\b \boxtimes \mathbb{U}^\a$ for any
$0<\alpha,\beta<\gamma$ such that $\alpha+\beta=\gamma$. Using the
induction hypothesis, we deduce that
$$\mathbf{b}_{\mathrm{Res}^{\b,\a}(\mathbb{P})}=\Delta_{\b,\a}(\mathbf{b}_{\mathbb{P}})
\in \mathfrak{W} \otimes \mathfrak{W}.$$
We conclude using the following result, which is proved in the appendix.

\begin{lem}\label{C:8.1} If $\mathbf{b} \in \mathbf{B}_{\mathbf{H}_{p_i}}$ is of degree $\gamma$ with $\gamma >
\delta$ and if $\Delta_{\b,\a}(\mathbf{b}) \in
\mathbf{U}^+_v(\widehat{\mathfrak{sl}}_{p_i}) \otimes
\mathbf{U}^+_v(\widehat{\mathfrak{sl}}_{p_i})$ for all
$0<\alpha,\beta<\gamma$ such that $\alpha+\beta=\gamma$ then 
$\mathbf{b}\in \mathbf{U}^+_v(\widehat{\mathfrak{sl}}_{p_i})$.
\end{lem}

\vspace{.1in}

We now treat the general case. Let $\mathbb{P} \in
\mathcal{P}^{\gamma}_{|\gamma|}$ with $supp(\mathbb{P}) \subset
Y_{\gamma_1, \ldots, \gamma_N}$,
with arbitrary $\gamma_i$. It is clear that
$\mathbb{Q}:=\mathrm{Res}^{\gamma_N,\ldots,\gamma_1}(\mathbb{P})$ is
supported on the
disjoint union $\bigsqcup_{\beta_1^{(i)} + \cdots +
   \beta_N^{(i)}=\gamma_i} Y_{(\beta_j^{(N)})_j}\times \cdots \times
   Y_{(\beta_j^{(1)})_j}$.
Let $\mathbb{Q}_0$ be the restriction of $\mathbb{Q}$ to
$Y_{\gamma_N} \times \cdots \times
Y_{\gamma_1}$. By Lemma~\ref{L:14}, $\mathbb{Q}_0 \in \mathbb{U}^{\gamma_N}
\boxtimes \cdots \boxtimes \mathbb{U}^{\gamma_1}$, hence by the above
result we have $\mathbf{b}_{\mathbb{Q}_0} \in \mathfrak{W} \otimes \cdots
\otimes \mathfrak{W}$. But it is easy to see that , up to a shift,
$\mathbb{P}=\mathrm{Ind}^{\gamma_N, \ldots, \gamma_1}(\mathbb{Q}_0)$, so that
$\mathbf{b}_{\mathbb{P}} \in \mathfrak{W}$ as desired.\qed

\vspace{.1in}

\begin{lem}\label{L:18} Let $\mathbb{P} \in
\mathcal{P}^{\alpha+l\delta}_{|\alpha|}$. There exists complexes
$\mathbb{T} \in \mathcal{P}^\alpha_{|\alpha|}$ and $\mathbb{R} \in
\mathcal{P}^{l\delta}_0$ such that $\mathrm{Ind}^{l\delta,\a}(
\mathbb{R} \boxtimes \mathbb{T}) =\mathbb{P} \oplus \mathbb{P}_1$
with $\mathbb{P}_1$ being a sum of (shifts) of sheaves in
$\mathcal{P}^{\alpha+l\delta}_{>|\alpha|}$.
\end{lem} \noindent
\textit{Proof.} It is enough to prove the following weaker
statement~: there exists complexes $\mathbb{T}_i \in
\mathcal{P}^\alpha_{|\alpha|}$, $\mathbb{R}_i \in
\mathcal{P}^{l\delta}_0$ and integers $l_i \in \Z$ such that
$\mathrm{Ind}^{l\delta,\a}(\bigoplus_i \mathbb{R}_i \boxtimes
\mathbb{T}_i[l_i]) =\mathbb{P} \oplus \mathbb{P}_1$ with
$\mathbb{P}_1$ being a sum of (shifts) of sheaves in
$\mathcal{P}^{\alpha+l\delta}_{>|\alpha|}$. Indeed, since
$\mathbb{P}$ is simple there then exists $i$ such that
$\mathrm{Ind}^{l\delta,\a}(\mathbb{R}_i \boxtimes
\mathbb{T}_i[l_i]) =\mathbb{P} \oplus \mathbb{P}'_1$ with
$\mathbb{P}'_1$ satisfying the same conditions as $\mathbb{P}_1$.
The fact that $\mathrm{Ind}$ commutes with Verdier duality then
forces $l_i=0$. \\

We consider the diagram
$$\xymatrix{
Q_{\geq |\alpha|}^{\a+l\delta} & F_{\geq |\alpha|} \ar[l]_-{i}
\ar[r]^-{\kappa} & Q^{l\delta} \times Q^{\a}\\
& & Q^{l\delta} \times U^{\a} \ar[u]^-{j_1}}$$
where the top row is the restriction of the diagram for
$\mathrm{Res}^{l\delta,\a}$ to the closed
subvariety $Q^{\a+l\delta}_{\geq |\alpha|}$, and $j_1$ is the
embedding. By construction, the stalk of
$\widetilde{\mathrm{Res}}^{l\delta,\a}(\mathbb{P})=\kappa_!i^*(\mathbb{P})$
at a point $(\phi_2: \mathcal{E}^{l\delta} \tto
\mathcal{F}_2,\phi_1: \mathcal{E}^\alpha
\tto \mathcal{F}_1)$ is zero if $|[\mathcal{F}_{1,\Lambda}]|
+|[\mathcal{F}_{2,\Lambda}]| <|\alpha|$. Hence
$supp(j_1^*\kappa_!i^*(\mathbb{P})) \subset
U^{l\delta} \times Q^\alpha_{|\alpha|}$. Put $\mathbb{Q}=j_{1* !}(j_1^*\kappa_!
i^*(\mathbb{P}))$ (a direct summand of
$\widetilde{\mathrm{Res}}^{l\delta,\a}(\mathbb{P})$). Note that
$supp(\mathbb{Q}) \subset
Q^{l\delta} \times Q^\alpha_{|\a|}$. Also, by Lemma~\ref{L:14}, we have $\mathbb{Q} \in
\mathbb{U}^{l\delta} \boxtimes \mathbb{U}^{\a}$
so there exists $\mathbb{T}_i \in \mathcal{P}^\alpha_{|\alpha|}$,
$\mathbb{R}_i \in \mathcal{P}^{l\delta}_0$
and integers $d_i \in \Z$ such that $\mathbb{Q}=\bigoplus_i
\mathbb{R}_i \boxtimes
\mathbb{T}_i[d_i]$. We claim that, up to a shift, this $\mathbb{Q}$
satisfies the requirements of the Lemma.
To see this, first consider the diagram
$$\xymatrix{
Q^{\a+l\delta}_{|\a|} \ar[d]_-{j_4} & E''_1 \ar[l]_-{p'_3}
\ar[d]_-{j_3} & E'_1 \ar[l]_-{p'_2} \ar[d]_-{j_2}
\ar[r]^-{p'_1} &  U^{l\delta} \times Q^\a_{|\a|} \ar[d]_-{j_1} \\
Q^{\a+l\delta}_{\geq |\a|}  & E'' \ar[l]_-{p_3} & E' \ar[l]_-{p_2}
\ar[r]^-{p_1} &  Q^{l\delta}\times Q^\a_{|\a|}}$$
where the bottom row is the restriction of the induction diagram to
the closed subvariety
$Q^{l\delta} \times Q^\alpha_{|\a|} \subset Q^{l\delta} \times Q^\a$, where
$$E'_1=(p'_1)^{-1}(U^{l\delta}\times Q^\alpha_{|\alpha|})\subset E',$$
$$E''_1=\{(\phi,(V_s)_s)\;| \underline{V} \simeq \mathcal{E}^\a,
\phi(\underline{V})_{\Lambda}=\phi(\underline{V})=
\phi(\mathcal{E}^{\a+l\delta})_\Lambda\} \subset E''$$
and where the vertical arrows are the (open) embeddings.
Note that all the squares in the above diagram are cartesian. On the
other hand, if $(\phi: \mathcal{E}^{\alpha+l\delta}
\tto \mathcal{F}) \in Q_{|\alpha|}^{\a+l\delta}$ then there exists a
unique subsheaf $\mathcal{G} \subset \mathcal{F}$
such that $\mathcal{G}=\mathcal{G}_{\Lambda}$ and
$[\mathcal{G}]=\alpha$. Hence $p'_3$ is an isomorphism. Standard
base change arguments yield
\begin{equation}\label{E:66}
\begin{split}
j_4^* \widetilde{\mathrm{Ind}}^{\alpha,l\delta}(\mathbb{Q})&=j_4^*
p_{3!}p_{2\flat}p_1^*(\mathbb{Q})\\
&=p'_{3!}p'_{2\flat}{p'_1}^* j_1^*(\mathbb{Q}).
\end{split}
\end{equation}
Next, we look at the diagram
$$\xymatrix{
Q_{|\a|}^{\a+l\delta} \ar[d]_-{j_4} & F_{|\a|} \ar[l]_-{i'}
\ar[d]_-{h} \ar[r]^-{\kappa'} &
U^{l\delta}\times Q^\alpha_{|\a|} \ar[d]_-{j_1}\\
Q^{\a+l\delta}_{\geq |\a|} & F_{\geq |\a|} \ar[l]_-{i}
\ar[r]^-{\kappa} & Q^{l\delta}\times Q^\a}$$
where $F_{|\a|}=\kappa^{-1}(U^{l\delta}\times Q^{\a}_{|\a|})$.
Similar arguments as above give
\begin{equation}\label{E:67}
\begin{split}
j_1^*\kappa_!i^*(\mathbb{P})&=j_1^*(\mathbb{Q})\\
&=\kappa'_!{i'}^*j_4^*(\mathbb{P}).
\end{split}
\end{equation}
Finally, consider the diagram
$$\xymatrix{
Q_{|\a|}^{\a+l\delta} & F_{|\a|} \ar[l]_-{i'} \ar[d]_-{a}
\ar[r]^-{\kappa'} & U^{l\delta} \times Q^\a_{|\a|}\\
& E'_1 \ar[ul]_-{p''_2} \ar[ur]^-{p'_1}&}$$
where $a: F_{|\a|} \to E'_1$ is the canonical embedding and $p''_2=p'_3p'_2$.
Recall that we have fixed a subspace $\underline{V}
\subset \mathcal{E}^{\a+l\delta}$ in the definition of $F_{|\a|}$.
Let $P \subset G^{\a+l\delta}$ be the
stabilizer of $\underline{V}$ and let $U$ be its unipotent radical.
There are canonical identifications
$Q^{\a+l\delta}_{|\a|} \simeq E''_1 \simeq G^{\a+l\delta}
\underset{P}{\times} F_{|\a|},$ and
$E'_1 \simeq G^{\a+l\delta} \underset{U}{\times} F_{|\a|}$.
Furthermore, for any $(\phi_2: \mathcal{E}^{l\delta} \tto \mathcal{F}_2,
\phi_1: \mathcal{E}^\a\tto \mathcal{F}_1)
\in U^{l\delta} \times Q^\a_{|\a|}$
we have $\mathrm{Ext}(\mathcal{F}_1,\mathcal{F}_2)=0$ since
$\mathcal{F}_1$ and $\mathcal{F}_2$ have disjoint
supports. This implies that $F_{|\a|} \simeq P \underset{G^\a \times
G^{l\delta}}{\times} (U^{l\delta}\times Q^\a_{|\a|})$.

\vspace{.1in}

Now, ${i'}^*j_4^*(\mathbb{P})$ is $P$-equivariant hence there exists
$\mathbb{P}' \in \mathcal{Q}_{G^{l\delta}}(Q^{l\delta})\boxtimes  \mathcal{Q}_{G^\a}(Q^\alpha)$
such that
${i'}^*j_4^*(\mathbb{P})={\kappa'}^*(\mathbb{P}')$. Then
$\kappa'_!{i'}^*j_4^*(\mathbb{P})=\kappa'_!{\kappa'}^*(\mathbb{P}')=\mathbb{P}'[-2
d]$
where $d=\mathrm{dim}\;(P/( G^{l\delta}\times G^{\a}))$ and
\begin{equation*}
\begin{split}
a^*{p'_1}^*\kappa'_!{i'}^*j_4^*(\mathbb{P})&=a^*{p'_1}^*(\mathbb{P}')[-2d]\\
&={\kappa'}^*(\mathbb{P}')[-2d]\\
&={i'}^*j_4^*(\mathbb{P})[-2d]\\
&=a^*{p''_2}^*j_4^*(\mathbb{P})[-2d].
\end{split}
\end{equation*}
But $a^*: \mathcal{Q}_{G^{\a+l\delta}}(Q^{\a+l\delta}) \to
\mathcal{Q}_P(F_{|\a|})$ is an equivalence, hence
${p'_1}^*\kappa'_!{i'}^*j_4^*(\mathbb{P})={p''_2}^*j_4^*(\mathbb{P})[-2d]$ and
\begin{equation}\label{E:69}
p''_{2\flat}{p'_1}^*\kappa'_!{i'}^*j_4^*(\mathbb{P})=j_4^*(\mathbb{P})[-2d].
\end{equation}
Combining (\ref{E:69}) with (\ref{E:66}) and (\ref{E:67}) finally yields
$j_4^*\widetilde{\mathrm{Ind}}^{l\delta,\a}(\mathbb{Q}[2d]) =j_4^*
(\mathbb{P})$, i.e
$$\widetilde{\mathrm{Ind}}^{l\delta,\a}(\bigoplus_i \mathbb{R}_i
\boxtimes \mathbb{T}_i[d_i+2d])=\mathbb{P}
+ \mathbb{P}_1$$
with $supp(\mathbb{P}_1) \subset Q_{>|\a|}^{\a+l\delta}$. We are 
done.\qed

\vspace{.2in}

\paragraph{\textbf{8.2.}} We may now prove that
${\mathfrak{U}}_{\mathbb{A}}^{tor}=\mathfrak{W}$. We argue by
induction. Fix $\gamma \in K^+(X)$ and let us assume that
${\mathfrak{U}}_{\mathbb{A}}^{tor}[\gamma']=\mathfrak{W}[\gamma']$
for all $\gamma'<\gamma$. We will prove by descending induction on
$\alpha$ that $\mathbf{b}_{\mathbb{P}} \in \mathfrak{W}$ for all
$\mathbb{P} \in \mathcal{P}^{\gamma}_{\geq |\alpha|}$. If
$\alpha=\gamma$ this follows from Lemma~\ref{L:17}. So assume that
we have $\mathbf{b}_{\mathbb{P}} \in \mathfrak{W}$ for all
$\mathbb{P} \in \mathcal{P}^{\gamma}_{>|\alpha|}$ for some
$\alpha$ such that $\gamma=\alpha+l\delta$, and let $\mathbb{P}
\in \mathcal{P}^{\gamma}_{|\alpha|}$. Applying Lemma~\ref{L:18}
gives
$$\mathbf{b}_{\mathbb{P}}\in
\mathbf{b}_{\mathbb{R}}\mathbf{b}_{\mathbb{T}} +
\bigoplus_{\mathbb{P}' \in
\mathcal{P}^\gamma_{>|\alpha|}}\mathbb{A}
\mathbf{b}_{\mathbb{P}'}$$ for some $\mathbb{T} \in
\mathcal{P}^\alpha_{|\alpha|}$, $\mathbb{R} \in
\mathcal{P}^{l\delta}_0$. If $\alpha \neq 0$ then by the induction
hypothesis the right-hand side belongs to $\mathfrak{W}$ and hence
so does $\mathbf{b}_{\mathbb{P}}$. Now assume that $\alpha=0$, so
that by Lemma~\ref{L:16} we have $\mathbb{P}=
\mathbf{IC}(U^{l\delta}_{((1), \ldots,(1))}, \sigma)$ for some
$\sigma \in Irrep\;\mathfrak{S}_l$. It is well-known that the
Grothendieck group $K_0(\mathfrak{S}_l)$ is spanned by the class
of the trivial representation and the set of classes
$ind^{\mathfrak{S}_l}_{\mathfrak{S}_{n_r}\times \cdots \times
\mathfrak{S}_{n_1}}(\sigma_r \otimes \cdots \otimes \sigma_1)$ for
$\sigma \in K_0(\mathfrak{S}_{n_i})$ and $\sum_i n_i=l$ with $n_i
< l$. By the induction hypothesis,
$\mathbf{b}_{\mathbf{IC}(U^{n_i\delta}_{((1),\ldots,(1))},\sigma_i)}
\in \mathfrak{W}$ for all $n_i <l$. Using (\ref{E:625}) we
conclude that there exists $\mathbf{x} \in \mathfrak{W}$ such that
\begin{equation}\label{E:690}
\mathbf{b}_{\mathbf{IC}(U^{l\delta}_{((1),\ldots,(1))},\sigma)}-\mathbf{x} \in
\bigoplus_{\mathbb{P} \in \mathcal{P}^{l\delta}_{>0}} \mathbb{A}
\mathbf{b}_{\mathbb{P}}.
\end{equation}
By the induction hypothesis again, the right-hand side of
(\ref{E:690}) is in $\mathfrak{W}$, and therefore
so is
$\mathbf{b}_{\mathbf{IC}(U^{l\delta}_{((1),\ldots,(1))},\sigma)}$.
This closes the induction and concludes
the proof of Theorem~\ref{T:5} i) \qed

\vspace{.2in}

\section{Proof of Theorem 5.1.$ii$)}

\vspace{.1in}

\paragraph{\textbf{9.1.}} We start the proof of Theorem~\ref{T:5}
ii). Let us temporarily denote again by $\mathfrak{W}$ the
subalgebra of $\widehat{\mathfrak{U}}_{\mathbb{A}}$ generated by
${\mathfrak{U}}_{\mathbb{A}}^{tor}$ and the elements
$\mathbf{1}_{l[\mathcal{O}_X(n)]}$ for $n \in \Z$ and $l \in
\mathbb{N}$. We have to show the following statement~:

\vspace{.05in}

\noindent
\textit{a')\;For any} $\alpha \in K^+(X)$, $\mathbb{P}=(P_n)_n \in \mathcal{P}^{\alpha}$
\textit{and for
any} $m \in \Z$
\textit{there exists}
$\mathbb{S},\mathbb{S}' \in \mathbb{U}^\alpha$ \textit{such that}
$S'_m \oplus P_m \simeq S_m$ \textit{and} $\mathbf{b}_{\mathbb{S}}
-\mathbf{b}_{\mathbb{S}'} \in \mathfrak{W}$.

\vspace{.05in}

Observe that $\mathfrak{W}$ and $\mathcal{P}$ are invariant under
twisting by $\mathcal{L}_{\delta}$, and that we have, canonically,
$Q_n^{\alpha} \simeq Q_{n+1}^{\alpha+rank(\alpha)\delta}$.
Therefore, it is enough to show the following~:

\vspace{.05in}

\noindent
\textit{a)\;For any} $\alpha \in K^+(X)$ \textit{and for any}
$\mathbb{P}=(P_n)_n \in \mathcal{P}^{\alpha}$ such that $P_0 \neq 0$
\textit{there exists}
$\mathbb{S},\mathbb{S}' \in \mathbb{U}^\alpha$ \textit{such that}
$S'_{-1} \oplus P_{-1} \simeq R_{-1}$ \textit{and} $\mathbf{b}_{\mathbb{S}}
-\mathbf{b}_{\mathbb{S}'} \in \mathfrak{W}$.

\vspace{.05in}

Recall that
$O^{\mathcal{H}}_n$ denotes the $G_n^\a$-orbit in $Q_n^\alpha$
corresponding to a sheaf $\mathcal{H}$ of class $\alpha$ ( and
$O_n^{\mathcal{H}}=\emptyset$ if $\mathcal{H}$ is not generated by
$\{\mathcal{L}_s(n)\}$). For such a sheaf
$\mathcal{H}$, the collection of perverse sheaves
$\mathbf{IC}(O_n^{\mathcal{H}},\mathbf{1})$  defines a
simple object of $\mathbb{Q}^\alpha$, which we denote by
$\mathbf{IC}(O^{\mathcal{H}})$.
It will be convenient for us to prove at the same time as a) the following

\vspace{.05in}

\noindent
\textit{b)}\;For any
indecomposable vector bundle $\mathcal{F}$, we have $\mathbf{IC}(O^{\mathcal{F}})
\in \mathcal{P}$.

\vspace{.1in}

We argue by induction on the rank of $\alpha$. Moreover, by
Proposition~\ref{P:71} i), if $\mathbb{P} \in \mathcal{P}^{\alpha}$
and $P_0 \neq 0$ then $|\alpha|\geq 0$ and we may further argue by
induction on $|\alpha|$. The case of $rank(\alpha)=0$ is
the subject of Section~6. So let us fix $\alpha$ and let us assume
that $a)$ is proved for all $\beta$ with
$rank(\beta) < rank(\alpha)$ or $rank(\beta)=rank(\alpha)$ and
$|\beta| < |\alpha|$, and that $b)$ holds for all $\mathcal{F}$ with
$[\mathcal{F}]$ satisfying the same conditions.
Let $\mathbb{P} \in \mathcal{P}^\alpha$ such that $P_0 \neq 0$. We
will first show that $a)$ holds for
$\mathbb{P}$. We argue once again by induction, this time on the
generic HN type of $\mathbb{P}$ with respect to the
order ``$\prec$'' defined in Section~7.1. So we finally assume in addition
that $a)$ is proved for all $\mathbb{P}'$ with $HN(\mathbb{P}')
\prec HN(\mathbb{P})$. Let us write $HN(\mathbb{P})=(\alpha_1, \ldots,
\alpha_r)=\underline{\alpha}$ and $\underline{\alpha}'=(\alpha_r, \ldots, \alpha_1)$.
From $P_0 \neq 0$ and the fact that the
HN filtration splits it follows that $\mu(\alpha_i) \geq 0$ for all $i$.

\vspace{.1in}

\paragraph{} Let us first suppose that $r >1$. Consider the complex
$$\mathbb{R}=(R_m)_m=\mathrm{Ind}^{\underline{\alpha}'}
\circ \mathrm{Res}^{\underline{\alpha}'}(\mathbb{P}) \in
\mathbb{U}^{\alpha}.$$
By construction, $\mu(\alpha_1) > \ldots >\mu(\alpha_r)$. Thus, by Lemma~\ref{L:after81} i) we have
\begin{equation}\label{E:after81}
supp(R_m) \subset
\bigcup_{(\underline{\beta}) \preceq (\underline{\alpha})} HN_m^{-1}
(\underline{\beta})
\end{equation}
for any $m \in \Z$.
Let $j_m: HN_m^{-1}(\underline{\alpha}) \hookrightarrow Q_m^\alpha$
be the embedding.  Note that by repeated applications
of Lemma~\ref{L:72} it follows that for any $(\phi:
\mathcal{E}_m^{\alpha} \tto \mathcal{F}) \in
HN_m^{-1}(\underline{\alpha})$,
there exists a unique filtration $\mathcal{G}_1 \subset \cdots \subset \mathcal{G}_r=\mathcal{F}$
of $\mathcal{F}$ such that
$[\mathcal{G}_i/\mathcal{G}_{i-1}]=\alpha_i$. Moreover, we have
$\mathrm{Ext}(\mathcal{G}_j/\mathcal{G}_{j-1},\mathcal{G}_i/\mathcal{G}_{i-1})=0$
for any $j >i$ since $\mu(\alpha_j)>\mu(\alpha_i)$
Hence we may use Corollary~\ref{C:72} with
$U^{\alpha_i}_n=Q_n^{(\alpha_i)}$. We deduce that for
any $m\in\Z$ we have
\begin{equation}\label{E:mid71}
j_m^*(R_m) \simeq j_m^*(P_m).
\end{equation}
Now, by Lemma~\ref{L:14}, we have
$\mathrm{Res}^{\underline{\a}'}(\mathbb{P}) \in \mathbb{U}^{\alpha_r}
\hat{\boxtimes} \cdots
\hat{\boxtimes} \mathbb{U}^{\alpha_1}$. Note that for all $i$ there holds
$rank(\alpha_i) < rank(\alpha)$ or $rank(\alpha_i) = rank(\alpha)$ and
$|\alpha_i| < |\alpha|$. Hence by the first induction hypothesis,
there exists complexes
$\mathbb{S}, \mathbb{S}' \in  \mathbb{U}^{\alpha_r} \hat{\boxtimes} \cdots
\hat{\boxtimes} \mathbb{U}^{\alpha_1}$ such that
\begin{equation}\label{E:after91}
S'_{-1} \oplus \mathrm{Res}^{\underline{\a}'}(\mathbb{P})_{-1} \simeq S_{-1}
\end{equation}
and $\mathbf{b}_{\mathbb{S}'}-\mathbf{b}_{\mathbb{S}} \in
\mathfrak{W} \hat{\otimes} \cdots \hat{\otimes} \mathfrak{W}$.
 By (\ref{E:after91}) there exists complexes $\mathbb{V},\mathbb{V}'$ satisfying
$V_{-1} = V'_{-1}=0$ and
$\mathbb{S}' \oplus \mathrm{Res}^{\underline{\a}'}(\mathbb{P}) \oplus \mathbb{V}'
\simeq \mathbb{S} \oplus \mathbb{V}$.
Let us set $\mathbb{T}=\mathrm{Ind}^{\underline{\alpha}'}(\mathbb{V}),
\mathbb{T}'=\mathrm{Ind}^{\underline{\alpha}'}(\mathbb{V}')$. We have
$$\mathrm{Ind}^{\underline{\a}'}(\mathbb{S}') \oplus
\mathrm{Ind}^{\underline{\a}'}
\mathrm{Res}^{\underline{\a}'}(\mathbb{P}) \oplus \mathbb{T}'\simeq\mathrm{Ind}^{\underline{\a}'}
(\mathbb{S}) \oplus \mathbb{T},$$
and, by Lemma~\ref{L:after81} ii)
\begin{equation}\label{E:after92}
supp(\mathbb{T}) \subset \bigcup_{\underline{\beta} \prec \underline{\alpha}}
HN^{-1}(\underline{\beta}),
\qquad supp(\mathbb{T}') \subset \bigcup_{\underline{\beta} \prec \underline{\alpha}}
HN^{-1}(\underline{\beta})
\end{equation}

By (\ref{E:after81}) and (\ref{E:mid71}), we have
\begin{equation}\label{E:corlast}
\mathrm{Ind}^{\underline{\a}}
\mathrm{Res}^{\underline{\a}}(\mathbb{P})=\mathbb{R}\simeq
\mathbb{P} \oplus \mathbb{R}'
\end{equation}
 for some $\mathbb{R}'
\in \mathbb{U}^{\a}$ with $supp(\mathbb{R}') \subset
\bigcup_{\underline{\beta} \prec \underline{\a}}
HN^{-1}(\underline{\beta})$. By the second induction hypothesis
there exists $\mathbb{W}, \mathbb{W}' \in \mathbb{U}^{\alpha}$
such that $W_{-1} \oplus T'_{-1} \oplus R'_{-1}\simeq W'_{-1}
\oplus T_{-1}$ and
$\mathbf{b}_{\mathbb{W}'}-\mathbf{b}_{\mathbb{W}} \in
\mathfrak{W}$. Gathering terms, we finally obtain after
simplification
$$\mathrm{Ind}^{\underline{\a}'}(\mathbb{S}')_{-1} \oplus W'_{-1} \oplus P_{-1}
\simeq \mathrm{Ind}^{\underline{\a}'}(\mathbb{S})_{-1}
\oplus W_{-1},$$
where $\mathbf{b}_{\mathrm{Ind}^{\underline{\a}'}(\mathbb{S}')}
-\mathbf{b}_{\mathrm{Ind}^{\underline{\a}'}(\mathbb{S})}
+\mathbf{b}_{\mathbb{W}'}-\mathbf{b}_{\mathbb{W}} \in \mathfrak{W}$.
This closes the induction step when $r>1$.

\vspace{.2in}

\paragraph{\textbf{9.2.}} From now on, we suppose that $r=1$, i.e that $\mathbb{P}$ is
generically semistable~: $supp(P_0) \cap Q_0^{(\alpha)}$ is open
in $supp(P_0)$. Since $P_0$ is simple and $G_0^{\alpha}$-equivariant
it follows from Lemma~\ref{L:73} that
$supp(P_0) \cap Q_0^{(\alpha)}$ is a single $G_0^{\alpha}$-orbit
corresponding to some isoclass $\mathcal{F} \in
\mathcal{C}_{\mu(\alpha)}$. Since in addition $P_0$ is a $G_0^{\a}$-equivariant
perverse sheaf, we deduce that $P_0 \simeq \mathbf{IC}(O^{\mathcal{F}}_0,\mathbf{1})$,
and finally that $\mathbb{P}=\mathbf{IC}(O^\mathcal{F})$. Let us decompose $\mathcal{F}$ into isotypical
components $\mathcal{F}=\mathcal{G}_1 \oplus \cdots \oplus \mathcal{G}_s$ where $\mathcal{G}_i \simeq
\mathcal{H}_i^{\oplus d_i}$ and $\mathcal{H}_i$ are distinct indecomposable sheaves. Note that
$\mathrm{Ext}(\mathcal{H}_i,\mathcal{H}_j)=0$ for any $i,j$.

\vspace{.1in}

Let us first assume that $s>1$. In that case, we have $rank(\mathcal{G}_i) <
rank(\mathcal{F})$ or $rank(\mathcal{G}_i)=rank(\mathcal{F})$ and $|\mathcal{G}_i| < |\mathcal{F}|$
for all $i$. By induction hypothesis b), $\mathbf{IC}(O^{\mathcal{H}_i}) \in
\mathcal{P}^{[\mathcal{H}_i]}$. From $\mathrm{Ext}(\mathcal{H}_i,\mathcal{H}_i)=0$ it is easy to see that
$\mathbf{IC}(O^{\mathcal{G}_i})$ appears in the product
$\mathrm{Ind}^{[\mathcal{H}_i],\ldots [\mathcal{H}_i]}(\mathbf{IC}(O^{\mathcal{H}_i}) \boxtimes
\cdots \boxtimes \mathbf{IC}(O^{\mathcal{H}_i}))$. Hence we also have
$\mathbf{IC}(O^{\mathcal{G}_i}) \in \mathcal{P}^{[\mathcal{G}_i]}$. Now, by induction hypothesis a)
there exists complexes $\mathbb{S}_i,\mathbb{S}'_i$ such that
$\mathbf{b}_{\mathbb{S}_i}-\mathbf{b}_{\mathbb{S}'_i} \in \mathfrak{W}$ and $S_{i,-1} \oplus
\mathbf{IC}(O_{-1}^{\mathcal{G}_i}) \simeq S'_{i,-1}$. We can find complexes $\mathbb{V}_i,
\mathbb{V}'_i$ satisfying $V_{i,-1}=V'_{i,-1}=0$ and $\mathbb{S}_i \oplus \mathbb{V}_i \oplus
\mathbf{IC}(O^{\mathcal{G}_i}) \simeq \mathbb{S}_i' \oplus \mathbb{V}'_i$. For simplicity, put
$\beta_i=[\mathcal{G}_i]$ and $\underline{\beta}'=(\beta_s, \ldots, \beta_a)$. Note that $\mu(\beta_i)
=\mu(\beta_j)$ for any $i,j$. Applying the functor $\mathrm{Ind}^{\underline{\beta}'}$ we obtain
\begin{equation*}
\begin{split}
\mathrm{Ind}^{\underline{\beta}'}&(\mathbb{S}_s \boxtimes \cdots \boxtimes \mathbb{S}_1)
\oplus \mathrm{Ind}^{\underline{\beta}'}(\mathbf{IC}(O^{\mathcal{G}_s}) \boxtimes \cdots \boxtimes
\mathbf{IC}(O^{\mathcal{G}_1}))\oplus
\mathrm{Ind}^{\underline{\beta}'}(\mathbb{V}_s \boxtimes \cdots \boxtimes \mathbb{V}_1)\\
&\simeq \mathrm{Ind}^{\underline{\beta}'}(\mathbb{S}'_s \boxtimes \cdots \boxtimes \mathbb{S}'_1)
\oplus \mathrm{Ind}^{\underline{\beta}'}(\mathbb{V}'_s \boxtimes \cdots \boxtimes \mathbb{V}'_1).
\end{split}
\end{equation*}
Consider the induction diagram
$$\xymatrix{
Q_l^{\beta_s} \times \cdots \times Q_l^{\beta_1} & E' \ar[l]_-{p_1} \ar[r]^-{p_2} & E'' \ar[r]^-{p_3}&
Q_l^{\alpha}}.$$
Since any extension of the $\mathcal{G}_i$ is necessarily trivial and isomorphic to $\mathcal{F}$, we
have
$$p_3p_2p_1^{-1}(\overline{O_l^{\mathcal{G}_s}} \times \cdots \times \overline{O_l^{\mathcal{G}_1}})
\subset \overline{p_3p_2p_1^{-1}(O_l^{\mathcal{G}_s}\times \cdots \times O_l^{\mathcal{G}_1})} =
\overline{O_l^{\mathcal{F}}}.$$
Recall that $\overline{O_l^{\mathcal{G}_i}} \backslash O_l^{\mathcal{G}_i} =
\bigcup_{\underline{\gamma} \prec \beta_i} HN^{-1}(\underline{\gamma})$. Using Lemma~\ref{L:74} ii)
and the fact that $\mathrm{Hom}(\mathcal{G}_i, \mathcal{G}_j)=0$ for $i \neq j$ we deduce that
$p_3: p_3^{-1}(O_l^{\mathcal{F}}) \to O_l^{\mathcal{F}}$ is an isomorphism. Therefore
$p_{3!}p_{2\flat}p_1^*(\mathbf{IC}(O^{\mathcal{G}_s}) \boxtimes \cdots \boxtimes
\mathbf{IC}(O^{\mathcal{G}_1}))_{|O_l^{\mathcal{F}}}=\mathbf{1}$. Hence
$$\mathrm{Ind}^{\underline{\beta}'}(\mathbf{IC}(O^{\mathcal{G}_s}) \boxtimes \cdots \boxtimes
\mathbf{IC}(O^{\mathcal{G}_1}))=\mathbf{IC}(O^{\mathcal{F}}) \oplus \mathbb{R}$$
where
$supp(\mathbb{R}) \subset \overline{O^{\mathcal{F}}}\backslash O^{\mathcal{F}} \subset
\bigcup_{\underline{\gamma} \prec (\alpha)} HN^{-1}(\underline{\gamma})$.

\vspace{.05in}

Let us put $\mathbb{T}=
\mathrm{Ind}^{\underline{\beta}'}(\mathbb{V}_s \boxtimes \cdots \boxtimes \mathbb{V}_1),
\mathbb{T}'=
\mathrm{Ind}^{\underline{\beta}'}(\mathbb{V}'_s \boxtimes \cdots \boxtimes \mathbb{V}'_1)$.
By Lemma~\ref{L:after81} ii), we have
$supp(\mathbb{T}) \subset \bigcup_{\underline{\gamma} \prec (\alpha)}HN^{-1}(\underline{\gamma})$
and similarly for $\mathbb{T}'$. Finally, we get
\begin{equation}\label{E:after93}
\mathrm{Ind}^{\underline{\beta}'}(\mathbb{S}_s \boxtimes \cdots \boxtimes \mathbb{S}_1)
\oplus \mathbf{IC}(O^{\mathcal{F}}) \oplus \mathbb{R} \oplus \mathbb{T} \simeq
\mathrm{Ind}^{\underline{\beta}'}(\mathbb{S}'_s \boxtimes \cdots \boxtimes \mathbb{S}'_1)
\oplus \mathbb{T}'.
\end{equation}
By the induction hypothesis, there exists complexes $\mathbb{W}, \mathbb{W}'$ such that
$\mathbf{b}_{\mathbb{W}} -\mathbf{b}_{\mathbb{W}'} \in \mathfrak{W}$ and
$W'_{-1} \oplus R_{-1} \oplus T_{-1} \simeq W_{-1} \oplus T'_{-1}$. We deduce from (\ref{E:after93})
that
$$\mathrm{Ind}^{\underline{\beta}'}(\mathbb{S}_s \boxtimes \cdots \boxtimes \mathbb{S}_1)_{-1}
\oplus \mathbf{IC}(O^{\mathcal{F}})_{-1} \oplus W_{-1} \simeq
\mathrm{Ind}^{\underline{\beta}'}(\mathbb{S}'_s \boxtimes \cdots
\boxtimes \mathbb{S}'_1)_{-1} \oplus W'_{-1}.$$ This settles the
case $s>1$.

\vspace{.1in}

It remains to consider the situation $s=1$, which corresponds to
the case $\mathcal{F}=\mathcal{H}^{\oplus d}$ where $\mathcal{H}$
is indecomposable. If $\mathcal{H}=\mathcal{O}_X$ then
$\mathbb{P}=\mathbf{1}_{l[\mathcal{O}_X]}$ and
$\mathbf{b}_{\mathbb{P}} \in \mathfrak{W}$ by definition.
Otherwise let $\mathcal{G}, \mathcal{K}$ be the sheaves supplied
by Lemma~\ref{L:fin71}~: there is a unique subsheaf $\mathcal{G}'
\subset \mathcal{F}$ such that $\mathcal{G}' \simeq \mathcal{G}$
while $\mathcal{F}/\mathcal{G}' \simeq \mathcal{K}$. Arguing as
above, and using the induction hypothesis b) it is easy to see
that $\mathbf{IC}(O^{\mathcal{G}}) \in
\mathcal{P}, \mathbf{IC}(O^{\mathcal{K}}) \in
\mathcal{P}$. Consider as usual the induction
diagram
$$\xymatrix{ Q_l^{[\mathcal{H}]} \times Q_l^{[\mathcal{G}]} & E' \ar[l]_-{p_1} \ar[r]^-{p_2} & E'' \ar[r]^-{p_3}&
Q_l^{\a}}.$$
We have
$$p_3p_2p_1^{-1}(O^{\mathcal{K}} \times O^{\mathcal{G}}) \subset p_3p_2p_1^{-1}(\overline{O^{\mathcal{K}}}
\times \overline{O^{\mathcal{G}}}) \subset \overline{p_3p_2p_1^{-1}(O^{\mathcal{K}} \times O^{\mathcal{G}})},$$
thus $p_3p_2p_1^{-1}(\overline{O^{\mathcal{K}}} \times \overline{O^{\mathcal{G}}})$ is irreducible. Hence
$p_3p_2p_1^{-1}(\overline{O^{\mathcal{K}}} \times \overline{O^{\mathcal{G}}})\cap Q_l^{(\alpha)}$ is also
irreducible, and contains $O^{\mathcal{F}}$. By Lemma~\ref{L:73}, we deduce that $O^{\mathcal{F}}$
is a dense open subset of $p_3p_2p_1^{-1}(\overline{O^{\mathcal{K}}} \times \overline{O^{\mathcal{G}}})$.
In particular, $p_1$ being smooth and $p_2$ being a principal bundle, we have
$\mathrm{dim}\;p_2p_1^{-1}((\overline{O^{\mathcal{K}}}\backslash O^{\mathcal{K}}) \times
(\overline{O^{\mathcal{G}}}\backslash O^{\mathcal{G}})) < \mathrm{dim}\;
p_2p_1^{-1}(O^{\mathcal{K}} \times O^{\mathcal{G}})$. On the other hand, by construction, the map
$p_3^{-1}(O^{\mathcal{F}}) \cap p_2p_1^{-1}(O^{\mathcal{K}} \times O^{\mathcal{G}}) \to O^{\mathcal{F}}$
is an isomorphism, hence $\mathrm{dim}\;p_2p_1^{-1}(O^{\mathcal{K}} \times O^{\mathcal{G}})=\mathrm{dim}\;
O^{\mathcal{F}}$. Therefore,
$$p_3p_2p_1^{-1}((\overline{O^{\mathcal{K}}}\backslash O^{\mathcal{K}}) \times
(\overline{O^{\mathcal{G}}}\backslash O^{\mathcal{G}})) \cap O^{\mathcal{F}}=\emptyset.$$
It follows that
$\widetilde{\mathrm{Ind}}^{[\mathcal{K}],[\mathcal{G}]}(\mathbf{IC}(O^{\mathcal{K}}) \boxtimes
\mathbf{IC}(O^{\mathcal{G}}))_{|O^{\mathcal{F}}}=\mathbf{1}$ while
$$supp(\mathrm{Ind}^{[\mathcal{K}],[\mathcal{G}]}(\mathbf{IC}(O^{\mathcal{K}}) \boxtimes
\mathbf{IC}(O^{\mathcal{G}})))=\overline{O^{\mathcal{F}}}.$$
We deduce that
$$\mathrm{Ind}^{[\mathcal{K}],[\mathcal{G}]}(\mathbf{IC}(O^{\mathcal{K}}) \boxtimes
\mathbf{IC}(O^{\mathcal{G}})))=\mathbf{IC}(O^{\mathcal{F}}) \oplus \mathbb{R}$$
where $supp(\mathbb{R}) \subset
\overline{O^{\mathcal{F}}}\backslash O^{\mathcal{F}} \subset \bigcup_{\underline{\gamma} \prec (\alpha)}
HN^{-1}(\underline{\gamma})$. We may now conclude in exactly the same manner as in the case $s>1$ above.
Finally, note that this also closes the induction for statement b). Theorem~5.1 ii) is proved.\qed

\vspace{.2in}

\paragraph{\textbf{9.3. Corollaries.}} Let us put
$$\mathcal{P}^{tor}=\bigsqcup_{\a \in K(X)^{tor}}
\mathcal{P}^\a,$$
$$\mathcal{P}^{vb}=\{\mathbb{P} \in \bigsqcup_{\alpha}
\mathcal{P}^\alpha\;|HN(\mathcal{P})=(\a_1, \ldots, \a_r), \;
\mu(\a_1) \neq \infty\}$$ (in other words, $\mathcal{P}^{vb}$ is
the set of all simple perverse sheaves whose support generically
consists of vector bundles). We set
$\widehat{\mathfrak{U}}^{vb}_{\mathbb{A}}=\widehat{\bigoplus}_{\mathbb{P} \in
\mathcal{P}^{vb}} \mathbb{A} \mathbf{b}_{\mathbb{P}}$ (here $\widehat{\bigoplus}$ denotes
a possibly infinite, but admissible sum).

\vspace{.1in}

\begin{cor}\label{Co:90}
The following hold~:
\begin{enumerate}
\item[i)] For any simple perverse sheaf $\mathbb{P} \in \mathcal{P}$
which is neither in $\mathcal{P}^{vb}$ nor in $\mathcal{P}^{tor}$
there exists simple perverse sheaves $\mathbb{R} \in
\mathcal{P}^{vb}$ and $\mathbb{T} \in \mathcal{P}^{tor}$ such that
$$\mathrm{Ind}(\mathbb{R} \boxtimes \mathbb{T})=\mathbb{P} \oplus
\mathbb{P}_1$$
with $supp(\mathbb{P}_1) \subset\bigcup_{\underline{\beta} \prec HN(\mathbb{P})}
HN^{-1}(\underline{\beta})$.
\item[ii)] Multiplication defines an isomorphism of
$\mathbb{A}$-modules $\widehat{\mathfrak{U}}_{\mathbb{A}} \simeq
\widehat{\mathfrak{U}}^{vb}_{\mathbb{A}} \otimes
\mathfrak{U}_{\mathbb{A}}^{tor}$.
\item[iii)] Let $\mathcal{I}$ denote the set of all isoclasses of
vector bundles in $Coh_G(X)$. We have
$\mathcal{P}^{vb}=\{\mathbf{IC}(O^\mathcal{F})\;|\; \mathcal{F}
\in \mathcal{I}\}$.
\end{enumerate}
\end{cor}

\noindent \textit{Proof.} By Section~9.1 (\ref{E:corlast}), there
exists complexes $\mathbb{R}^i_1 \in \mathcal{P}^{vb}$,
$\mathbb{T}^i \in \mathcal{P}^{tor}$ and integers $l_i$ such that
$$\mathrm{Ind}(\bigoplus_i \mathbb{R}^i\boxtimes \mathbb{T}_i[l_i]) =
\mathbb{P} \oplus \mathbb{P}_1$$ with $\mathbb{P}_1$ satisfying
the same conditions as above. But since $\mathbb{P}$ is simple, it
appears in a product $\mathrm{Ind}(\mathbb{R}^i \boxtimes
\mathbb{T}_i[l_i])$ for some fixed $i$. Furthermore, from the fact
that $\mathrm{Ind}$ commutes with Verdier duality, we deduce that
$l_i=0$. This proves i). The surjectivity of the multiplication
map $\widehat{\mathfrak{U}}^{vb}_{\mathbb{A}} \otimes
\mathfrak{U}_{\mathbb{A}}^{tor}
\to\widehat{\mathfrak{U}}_{\mathbb{A}}$ now follows from i) and an
induction argument on the $HN$ type. To see that it is also
injective, let us assume on the contrary that there exists a
nontrivial relation $\sum_i \mathbf{b}_{\mathbb{R}_i}
\mathbf{b}_{\mathbb{T}_i}=0$ in degree $\a$, where $\mathbb{T}_i
\in \mathcal{P}^{\a_i}$ are simple and distinct. Reindexing if
necessary, we may assume that $|\alpha_1|=min\{|\alpha_i|\}$. Note
that we have $supp\; \mathrm{Ind}^{\a-\a_i,\a_i}(\mathbb{R}_i
\boxtimes \mathbb{T}_i)_n \subset \{ \phi:\mathcal{E}_n^\a \tto
\mathcal{F}\;|\; |\tau(\mathcal{F})| \geq |\alpha_i|\}$. Defining
$Q_{n,vb}^{\b}=\{\phi: \mathcal{E}_n^\b \tto \mathcal{F}\;|\;
\mathcal{F}\; \text{is\;a\;vector\;bundle}\}$ and denoting by
$h_\b: Q_{n,vb}^\b \to Q_n^\b$ the open embedding, we deduce that
up to a shift we have
\begin{equation}\label{E:999up}
\begin{split}
0&=(h^*_{\a-\a_1} \otimes 1)\;
\mathrm{Res}^{\a-\a_1,\a_1}
\big(\bigoplus_i\mathrm{Ind}^{\a-\a_i,\a_i}(\mathbb{R}_i \boxtimes \mathbb{T}_i)\big)\\
&=\bigoplus_{i; \alpha_i=\alpha_1} h_{\alpha-\a_i}^*(\mathbb{R}_i)
\boxtimes \mathbb{T}_i.
\end{split}
\end{equation}
 But since
$\mathbb{R}_i \in \mathcal{P}^{vb}$, $h_{\a-\a_i}^*(\mathbb{R}_i)
\neq 0$ for all $i$, in contradiction with (\ref{E:999up}).
This proves ii).\\
Finally, we show iii). Let $\mathbb{P} \in \mathcal{P}^{vb}$ be of
class $\alpha$ and put $HN(\mathbb{P})=\underline{\a}$. There
exists $m$ small enough such that $P_m \neq 0$. By Lemma~6.5 and
\cite{S1}, Lemma~7.1 there are only finitely many isomorphism
classes of vector bundles of class $\alpha$ generated by
$\{\mathcal{L}_s(m)\}$, i.e there are only finitely many
$G_m^\a$-orbits in $Q_{m,vb}^\a$. In particular, there is a unique
dense open orbit $O^\mathcal{F}_m$ in $supp(P_m)$. Since $P_m$ is
simple and $G_m^\a$-equivariant, we must have
$P_m=IC(O^{\mathcal{F}}_m)$ and hence
$\mathbb{P}=\mathbf{IC}(O^\mathcal{F})$. Conversely, by
Section~9.1, statement b), $\mathbf{IC}(O^\mathcal{G}) \in
\mathcal{P}$ for any indecomposable vector bundle $\mathcal{G}$.
It is easy to see that if $\mathcal{F}=\mathcal{G}_1 \oplus \cdots
\oplus \mathcal{G}_r$ with $\mathcal{G}_i$ indecomposable and
$\mu([\mathcal{G}_1]) \geq \cdots \geq \mu([\mathcal{G}_r])$ then
$\mathbf{IC}(O^\mathcal{F})$ appears in the product
$\mathrm{Ind}^{[\mathcal{G}_1], \ldots,
[\mathcal{G}_r]}(\mathbf{IC}(O^{\mathcal{G}_1}) \boxtimes \cdots
\boxtimes \mathbf{IC}(O^{\mathcal{G}_r}))$, from which it follows
that $\mathbf{IC}(O^{\mathcal{F}}) \in \mathcal{P}$ as claimed.
\qed

\vspace{.1in}

\paragraph{\textbf{Remark.}} Statements i) and ii) of the above Lemma (and their proofs)
are valid for an arbitrary pair $(X,G)$.

\vspace{.2in}

\section{Canonical basis of $\mathbf{U}_v(\mathcal{L}\n)$}

\vspace{.1in}

In this section, we use Theorem~5.2 to construct a (topological)
canonical basis of $\mathbf{U}_v(\mathcal{L}\mathfrak{n})$. We
first describe the precise link between
$\mathbf{U}_v(\mathcal{L}\mathfrak{n})$ and
$\widehat{\mathfrak{U}}_{\mathbb{A}}$.

\vspace{.1in}

\paragraph{\textbf{10.1. Cyclic quivers.}} Fix $i \in \{1, \ldots, N\}$. Recall the subalgebra $\mathfrak{U}_i^{tor}$
of $\mathfrak{U}^{tor}_{\mathbb{A}}$ corresponding to complexes of
sheaves supported at the exceptional point $\lambda_i \in \Lambda$
(see Lemma~\ref{L:17} for the precise definition), and the algebra
isomorphism
\begin{align*}
\tau_i:\;\mathfrak{U}_i^{tor} &\to \mathbf{H}_{p_i}\\
\mathbf{b}_{\mathbb{P}}\in \mathfrak{U}^{tor}_i[\gamma] &\mapsto
\big( x \mapsto
v^{\mathrm{dim}(G_{\gamma})-\mathrm{dim}(G^{\gamma})}\sum_i
\mathrm{dim}_{\overline{\mathbb{Q}}_l}
\mathcal{H}^i_{|x}(\mathbb{P})v^{-i}\big).
\end{align*}
Define elements $\zeta_l^{(i)} \in \mathbf{H}_{p_i}$ by the
relation $1+\sum_{l \geq 1} \zeta_l^{(i)} s^l= exp(\sum_{l\geq 1}
\mathbf{h}_{i,l}s^l)$. Keeping the notations of (\ref{E:Hall1}),
we have $\zeta_l^{(i)}=\sum_{\lambda,
|\lambda|=l}{1}_{O^{\lambda}}$ (see \cite{S1}, Section~6.3.).

\vspace{.1in}

We postpone the proof of this technical Lemma to the appendix~:

\begin{lem}\label{L:101}
There exists unique elements $u_l^{(i)} \in
\mathbf{U}_v^+(\widehat{\mathfrak{sl}}_{p_i}) \subset
\mathbf{H}_{p_i}$ for $l \geq 1$ satisfying the following set of
relations~:
\begin{equation}\label{E:equul}
1_{\mathcal{N}^{(p_i)}_{(l,\ldots,l)}}=\zeta_l^{(i)}+\zeta^{(i)}_{l-1}u_1^{(i)}+
\cdots + \zeta_1^{(i)}u_{l-1}^{(i)} + u_l^{(i)}.
\end{equation}
\end{lem}

\vspace{.2in}

\paragraph{\textbf{10.2. The morphism $\Psi: \mathbf{U}_v(\mathcal{L}\n) \to
\widehat{\mathfrak{U}}_{\mathbb{A}}$.}} We are at last in
position to give the perverse sheaf analogue of the main result in
\cite{S1}.
We define inductively some elements ${\boldsymbol{\xi}_l} \in
\widehat{\mathfrak{U}}_{\mathbb{A}}$ by the relations
\begin{equation}\label{E:1001}
\boldsymbol{\xi}_l=\mathbf{b}_{l\delta}-
\sum_{\underset{n<l}{n+n_1+\cdots +n_N=l}}
{\boldsymbol{\xi}_n}\tau_1^{-1}(u_{n_1}^{(1)}) \cdots
\tau_N^{-1}(u_{n_N}^{(N)}).
\end{equation}

\begin{theo}\label{T:1001} Let $(X,G)$ be arbitrary. For $\underline{\alpha}=(\alpha_1, \ldots, \alpha_n) \in (K(X)^{tor})^n$
we set $d(\underline{\alpha})=\langle [\mathcal{O}_X],\sum_j
\alpha_j \rangle$. The assignment
\begin{align}
E_{(i,j)} &\mapsto \mathbf{b}_{\alpha_{(i,j)}},\qquad \text{for\;}(i,j) \in \aleph,\\
\xi_l &\mapsto \boldsymbol{\xi}_l\qquad \text{for\;}l \geq 1,\\
E_{\star,t} & \mapsto \sum_{m \geq 0} \sum_{\underline{\alpha} \in
(K(X)^{tor})^m} (-1)^m v^{d(\underline{\alpha})}
{\mathbf{b}}_{[\mathcal{O}_X(t)]-\sum
\alpha_j}{\mathbf{b}}_{\alpha_m} \cdots {\mathbf{b}}_{\alpha_1}
\qquad \text{for\;} t \in \Z, \label{A:1001}
\end{align}
extends to an algebra homomorphism $\Psi:
\mathbf{U}_v(\mathcal{L}\n)\to
\widehat{\mathfrak{U}}_{\mathbb{A}}$.
\end{theo}
\noindent \textit{Proof.} The computations needed to prove that
$\Psi$ is well-defined, i.e that $\Psi(E_{(i,j)})$, $\Psi(\xi_l)$
and $\Psi(E_{\star,l})$ satisfy relations i)-v) in Section~1.3 are
a sheaf version of those done in \cite{S1}, Section~6, and may, in
principle, be checked directly. However, we may reduce ourselves
to the (easier) computations done in \cite{S1} by introducing the
following maps. Let $\mathfrak{U}^\a_{\mathbb{A},\geq n} \subset
\mathfrak{U}^\alpha_{\mathbb{A}}$ be the $\mathbb{A}$-span of
elements $\mathbf{b}_{\mathbb{P}}$ where $\mathbb{P}=(P_m)_{m \in
\mathbb{Z}}$ is such that $P_n \neq 0$, and consider the map
\begin{align*}
\tau_n^\a: \mathfrak{U}^\a_{\mathbb{A},\geq n} & \to
\mathbb{C}_{G_n^\a}(Q_n^\a) \otimes \mathbb{A},\\
\mathbf{b}_{\mathbb{P}} & \mapsto \big( x \mapsto
v^{-\mathrm{dim}\;G_n^\a} \sum_i
\mathrm{dim}_{\qlb}\mathcal{H}^i_{|x}(P_n)v^{-i}\big).
\end{align*}

\begin{lem} The map $\tau_n^\a$ is injective if $rank(\a) \leq 1$ or if $\alpha$ is
of the form $2[\mathcal{O}_X]+l\delta$ or $2[\mathcal{O}_X] +
l\delta+\alpha_{(i,1)}$ for some $l \in \Z$ and $i=1, \ldots, N$.
\end{lem}
\noindent \textit{Proof.} Let us first assume that $\a$ is a
torsion class and set $l=|\a|$. If $x \in \mathrm{Ker}\;\tau_n^\a$
we may, by Theorem~5.1.i), write
$$x=\sum_k P_k(\mathbf{b}_\delta, \cdots,
\mathbf{b}_{l\delta})T_k$$ where $P_k$ is a noncommutative
polynomial and $T_k$ belongs to the subalgebra generated by
$\mathbf{b}_{\a_{(i,j)}}$ for $(i,j) \in \aleph$. Set
$\mathfrak{U}'_{\mathbb{A}}=\bigoplus_{\mathbb{P} \in
\mathcal{P}'} \mathbb{A}\mathbf{b}_{\mathbb{P}}$ where
$\mathcal{P}'=\bigsqcup_{\beta} \mathcal{P}^\beta_{\geq 1}$. From
Lemma~\ref{L:18} we deduce that
$\mathfrak{U}'_{\mathbb{A}}=\sum_{(i,j)}
\mathfrak{U}_{\mathbb{A}}\mathbf{b}_{\a_{(i,j)}}$. By
Lemma~\ref{L:16}, we further have
$[\mathbf{b}_{l_1\delta},\mathbf{b}_{l_2\delta}]\in
\mathfrak{U}'_{\mathbb{A}}$ for any $l_1, l_2 \in \mathbb{N}$.
From this one sees that it is possible to reduce $x$ to a normal
form
$$x=\sum_k c_k \mathbf{b}_{l\delta}^{n_{l,k}} \cdots
\mathbf{b}_{\delta}^{n_{1,k}} T_k, \qquad c_k \in \mathbb{A}.$$
Rearranging the sum, we may further assume that $\sum_{h=1}^l h
\cdot n_{h,k}$ is maximal for $k=1, \ldots, m$, and we denote by
$l_0$ this common value. Now, let us fix a point $z \in
\mathbb{P}^1\backslash \Lambda$. Let $\mathcal{F}$ be a torsion
sheaf of length $l_0$ supported at $z$, and let $\mathcal{G}$ be a
torsion sheaf supported on the exceptional locus $\Lambda$.
Arguing as in Section~8.2 and using the maximality of $l_0$ one
checks that, up to a power of $v$,
\begin{equation}\label{E:10tsing}
0=\tau_n^\a(x)_{|O^{\mathcal{F} \oplus \mathcal{G}}}=\sum_{k=1}^m
c_k \tau_n^\a(\mathbf{b}_{l\delta}^{n_{l,k}} \cdots
\mathbf{b}_{\delta}^{n_{1,k}})_{|O^{\mathcal{F}}}\cdot
\tau_n^\a(T_k)_{|O^{\mathcal{G}}}.
\end{equation}
Define two closed subsets of $Q_n^\a$ by $\mathcal{N}_z=\{\phi:
\mathcal{E}^\a_n \tto \mathcal{F}\;|\; supp(\mathcal{F})=\{z\}\}$
and $ \mathcal{N}_{\Lambda}=\{\phi: \mathcal{E}^\a_n \tto
\mathcal{F}\;|\; supp(\mathcal{F})\subset \Lambda\}$. From
Proposition~\ref{P:HallS0} (for the quiver of type $A_0^{(1)}$)
and the description of $\mathcal{N}_z$ in Section~5.1 it follows
that the functions $\tau_n^\a(\mathbf{b}_{l\delta}^{n_{l,k}}
\cdots \mathbf{b}_{\delta}^{n_{1,k}})_{|\mathcal{N}_z}$ are all
linearly independent. Similarly, it is well-known (see e.g
\cite{L2}) that $\tau_n^\a(T_k)_{|\mathcal{N}_{\Lambda}}=0$ if and
only if $T_k=0$. These two facts combined with (\ref{E:10tsing})
imply that $x=0$.

\vspace{.05in}

Now we consider the case of a class $\alpha$ of rank one. The
arguments in Section~9.1 show that any $x \in
\mathfrak{U}_{\mathbb{A}}[\alpha]$ is equal to a sum
$x=\sum_{\mathcal{L}} \mathbf{b}_{[\mathcal{L}]} T_{\mathcal{L}}$
where $\mathcal{L}$ runs over $Pic(X,G)$ and $T_\mathcal{L} \in
\mathfrak{U}_{\mathbb{A}}[\a-[\mathcal{L}]]$. Note that for any
$\mathcal{L}$, $\a-[\mathcal{L}]$ is a torsion class. Fix $x \in
\mathrm{Ker}\;\tau_n^\a$ of the above form and, if $x \neq 0$,
choose among all $\mathcal{L}$ with $T_{\mathcal{L}} \neq 0$ one
which is maximal for the partial order on $K(X)$. Recall that by
Section~5.1, the set of $G_n^\a$-orbits in $Q_n^\a([\mathcal{L}])$
is in canonical bijection with the set of
$G_n^{\a-[\mathcal{L}]}$-orbits in $Q_n^{\a-[\mathcal{L}]}$, Hence
we may consider the restriction of $\tau_n^\a(x)$ to
$Q_n^\a([\mathcal{L}])$ as a $G_n^{\a-[\mathcal{L}]}$-invariant
function on $Q_n^{\a-[\mathcal{L}]}$. Using this identification
and the maximality of $\mathcal{L}$, we have, up to a power of
$v$,
$$0=\tau_n^\a(x)_{|Q_n^\a([\mathcal{L}])}=\tau_n^{\a-[\mathcal{L}]}(T_{\mathcal{L}}).$$
By the first part of the Lemma, this implies that
$T_\mathcal{L}=0$, a contradiction. Thus
$\mathrm{Ker}\;\tau_n^\a=\{0\}$ as wanted. Finally, the case of
$\alpha$ of rank two is solved in a similar way (observe that any
sheaf of class $\alpha$ is a sum of two line bundles and a torsion
sheaf- see \cite{S1}, Section~6.7). \qed

\vspace{.1in}

\noindent \textit{End of proof of Theorem~\ref{T:1001}.} Each
relation i)-v) is equivalent to certain equalities of semisimple
complexes $R_1=R_2$ in $\mathcal{Q}_{G_n^\a}(Q_n^\a)$ for $n \ll
0$ and $\a$ as in the above Lemma. Thus it is enough to prove the
corresponding equalities of functions
$\tau_n^\a(R_1)=\tau_n^\a(R_2)$. We first introduce some notation.
Recall that we fixed, in Section~2.1, an equivalence between the
category of torsion sheaves supported at some exceptional point
$\lambda_i$ and the category of nilpotent representations of the
quiver $A_{p_i-1}^{(1)}$. Let $C_i \subset Coh_G(X)$ be the full
subcategory consisting of all such torsion sheaves which
correspond to representations $(V,x)$ of $A_{p_i}^{(1)}$ where
$x_2, \ldots, x_{p_i}$ are all isomorphisms. Also, for $l \in \N$
and $n \in \Z$, we put
$$Z_{n}^l=\{\phi: \mathcal{E}_n^{l\delta} \tto
\mathcal{F}\;|\;\mathcal{F}_{|\lambda_i} \in C_i\;
\text{for\;all}\; i\}.$$ We claim that
\begin{equation}\label{E:10claim1}
\tau_n^{\a_{(i,j)}}({\mathbf{b}}_{\a_{(i,j)}})=v^{-1}1_{\a_{(i,j)}},
\qquad \tau_n^{l\delta}(\boldsymbol{\xi}_l)=1_{Z_n^l}, \qquad
\tau_n^{[\mathcal{O}_X(t)]}(\Psi(E_{\star,t}))=v^{-1}1_{O_n^{\mathcal{O}_X(t)}}.
\end{equation}
The first equality is obvious. The second equality is proved by
induction on $l$, using the defining property~\ref{E:equul} of
$u_l^{(i)}$. Finally, we compute
$$\tau_n(\Psi(E_{\star,t}))_{|O^{\mathcal{F}}}=v^{-1}\sum_{m \geq 0}
(-1)^m \sum_{(\alpha_0, \ldots, \alpha_m) \in A}
\mathrm{dim}_{\qlb}\;H^i(Quot_{\mathcal{F}}^{\alpha_0, \ldots,
\alpha_m}) v^{-1},$$ where $A=\{(\a_0, \ldots, \a_m)\;|\sum
\a_i=[\mathcal{O}_X(t)],\; rank(\a_0)=1\}$ and
$Quot_{\mathcal{F}}^{\a_0, \ldots, \a_m}$ denotes the hyperquot
scheme parameterizing filtrations of $\mathcal{F}$ with
subquotients of given classes $\alpha_0, \ldots, \alpha_m$. If
$\mathcal{F} \neq \mathcal{O}_X(t)$ then $\mathcal{F}=\mathcal{L}
\oplus \mathcal{T}$ where $\mathcal{L}$ is a line bundle and
$\mathcal{T}$ is a torsion sheaf. Now observe that
$Quot_\mathcal{F}^{\a_0, \ldots, \a_m} \simeq
Quot_{\mathcal{F}}^{[\mathcal{L}], \a_0-[\mathcal{L}], \ldots,
\a_m}$ and that therefore
$\tau_n(\Psi(E_{\star,t}))_{|O^{\mathcal{F}}}=0$. The last
equality in \ref{E:10claim1} follows. The set of equalities
$\tau^\a_n(R_1)=\tau_n^\a(R_2)$ now precisely amount to the
computations done in \cite{S1}. We leave the details to the
reader.\qed

\vspace{.1in}

\begin{theo}\label{T:1002} Let us assume that $X$ is of genus zero. Then the map $\Psi$
is injective and $\mathrm{Im}\;\Psi$ is dense in
$\widehat{\mathfrak{U}}_{\mathbb{A}}$.
\end{theo}
\noindent \textit{Proof.} We first show the density of
$\mathrm{Im}\;\Psi$. By construction, $\mathbf{b}_{l\a_{(i,j)}}$
belongs to $\mathrm{Im}\;\Psi$ for any $(i,j) \in \aleph$ and $l \geq 1$. 
From this and the definition of $\boldsymbol{\xi}$ it is also clear
that $\mathbf{b}_{l\delta} \in \mathrm{Im}\;\Psi$. Hence, by
Theorem~5.1.i), $\mathfrak{U}^{tor}_{\mathbb{A}} \subset
\mathrm{Im}\;\Psi$ and, according to Theorem~5.1 ii), it only
remains to check that $\mathbf{b}_{[\mathcal{O}_X(t)]} \in
\overline{\mathrm{Im}\;\Psi}$ for all $t \in \Z$. For this, we
define for any line bundle $\mathcal{L}'$ an element
$E_{\mathcal{L}'} \in \mathbf{U}_v(\mathcal{L}\n)$ as follows~:
$E_{\mathcal{O}_X(t)}=E_{\star,t}$ and if
$[\mathcal{L}'']=[\mathcal{L}'] + \alpha_{(i,j)}$ with
$\mathrm{dim\;Ext}(S_i^{(j)}, \mathcal{L}')=1$ then we put
$E_{\mathcal{L}''}=v^{-1}E_{(i,j)}E_{\mathcal{L}}-E_{\mathcal{L}}E_{(i,j)}$.
With this definition, we have, using \cite{S1}, Section~6.4,
\begin{equation}\label{E:10density}
\tau_n^{[\mathcal{L}']}(\Psi(E_{\mathcal{L}'}))=v^{-1}1_{O_n^{\mathcal{L}'}}\qquad
\text{for\;}n \ll 0.
\end{equation}
We claim that
\begin{equation}\label{E:10base}
\mathbf{b}_{[\mathcal{L}]}=\Psi(E_{\mathcal{L}}) +
\sum_{\mathcal{L}' < \mathcal{L}} v^{\langle
[\mathcal{L}],[\mathcal{L}]-[\mathcal{L}']\rangle}
\Psi(E_{\mathcal{L}'})\mathbf{b}_{[\mathcal{L}]-[\mathcal{L}']}.
\end{equation}
Indeed, it is enough to check the corresponding equality obtained
after applying $\tau_n^{[\mathcal{L}]}$, and this follows from
(\ref{E:10base}). Finally, note that the elements
$\mathbf{b}_{[\mathcal{L}]-[\mathcal{L}]'}$ are in
$\mathfrak{U}_{\mathbb{A}}^{tor} \subset \mathrm{Im}\;\Psi$. This
proves that
$\overline{\mathrm{Im}\;\Psi}=\widehat{\mathfrak{U}}_{\mathbb{A}}$.
In order to obtain the injectivity of $\Psi$, we fix some integer
$n \in \Z$ and consider subalgebra $\mathbf{U}_v^{\geq
n}(\mathcal{L}\n)$ of $\mathbf{U}_v^{\geq n}$ generated by $\{
E_{\star,l}, \mathbf{H}_{p_i}\;|l \geq n, i=1, \ldots, N\}$. Since
$\mathbf{U}_v(\mathcal{L}\n)$ is the limit of $\mathbf{U}_v^{\geq
n}(\mathcal{L}\n)$ as $n$ tends to $-\infty$, it is enough to show
that the restriction of $\Psi$ to $\mathbf{U}_v^{\geq
n}(\mathcal{L}\n)$ is injective. In fact, the following holds.

\vspace{.05in}

\noindent \textbf{Claim.} The composition of maps $\xymatrix{
\mathbf{U}^{\geq n}_v(\mathcal{L}\n) \ar[r]^-{\Psi} &
\widehat{\mathfrak{U}}_{\mathbb{A}} \ar[r]^-{\pi} &
\mathfrak{U}_{\mathbb{A}, \geq n}}$ is an isomorphism.\\
\noindent \textit{Proof of Claim.} First of all, by
Corollary~\ref{Co:90} we have $\mathfrak{U}_{\mathbb{A},\geq n} =
\mathfrak{U}_{\mathbb{A}, \geq n}^{vb} \otimes
\mathfrak{U}_{\mathbb{A}}^{tor}$ where
$\mathfrak{U}^{vb}_{\mathbb{A}, \geq n}=\bigoplus_{\mathbb{P} \in
\mathcal{P}_{\geq n}^{vb}} \mathbb{A}\mathbf{b}_{\mathbb{P}}$ with
$$\mathcal{P}^{vb}_{\geq n}=\{
\mathbf{IC}(O^{\mathcal{F}})\;|\mathcal{F}\;
\text{is\;a\;vector\;bundle\;generated\;by\;}\{\mathcal{L}_s(n)\}\}.$$
Next, by construction, the algebra $\mathbf{U}_v^{\geq
n}(\mathcal{L}\n)$ contains all elements $E_{\mathcal{L}}$ with
$\mathcal{L}$ being a line bundle generated by
$\{\mathcal{L}_s(n)\}$. By \cite{S1} Prop.~7.6, any indecomposable
vector bundle $\mathcal{F}$ generated by $\{\mathcal{L}_s(n)\}$
can be obtained as a successive extension of line bundles also
generated by $\{\mathcal{L}_s(n)\}$. We have
$\Psi(E_{\mathcal{L}})=\mathbf{b}_{\mathbb{P}_{\mathcal{L}}}$ for
some unique (virtual) complex $\mathbb{P}_{\mathcal{L}} \in
\mathbb{U}^{[\mathcal{L}]}$, and this complex satisfies
$(\mathbb{P}_{\mathcal{L}})_{|Q_{n,vb}^{[\mathcal{L}]}}=\qlb_{Q_{n,vb}^{[\mathcal{L}]}}
[\mathrm{dim\;}Q_n^{[\mathcal{L}]}]$. It follows that for any
$\mathcal{F}$ as above, $\mathbf{IC}(O^{\mathcal{F}})$ appears in
some induction product $\mathrm{Ind}^{[\mathcal{L}_1], \ldots,
[\mathcal{L}_r]}(\mathbb{P}_{\mathcal{L}_1} \boxtimes \cdots
\boxtimes \mathbb{P}_{\mathcal{L}_r})$, where all $\mathcal{L}_i$
are generated by $\{\mathcal{L}_s(n)\}$. An induction argument on
the rank similar to \cite{S1}, Prop.~7.6 then shows that
$\mathbf{IC}(O^{\mathcal{F}}) \in \pi\Psi(\mathbf{U}^{\geq
n}_v(\mathcal{L}\n))$ for any indecomposable vector bundle
$\mathcal{F}$ generated by $\{\mathcal{L}_s(n)\}$, and more
generally for any vector bundle $\mathcal{F}$ generated by
$\{\mathcal{L}_s(n)\}$. In particular
$\mathfrak{U}_{\mathbb{A},\geq n}^{vb} \subset
\pi\Psi(\mathbf{U}^{\geq n}_v(\mathcal{L}\n))$. Finally note that
$\Psi(\mathbf{U}^{\geq n}_v(\mathcal{L}\n))$ also contains
$\mathfrak{U}^{tor}_{\mathbb{A}}$. This implies that $\pi \Psi$ is
surjective. The injectivity is proved by a graded dimension
argument identical to the one in \cite{S1}, Section~7.5. \qed

\vspace{.1in}


\vspace{.1in}

\begin{cor}\label{C:1001} Assume $X$ is of genus zero. Then
the map $\Psi$ extends to an isomorphism
from a completion $\widehat{\mathbf{U}}_v(\mathcal{L}\n)$ of $\mathbf{U}_v(\mathcal{L}\n)$ to
$\widehat{\mathfrak{U}}_{\mathbb{A}}$.
\end{cor}

Again, we conjecture that this result remains valid for arbitrary data
$(X,G)$ (the case of $X$ being of genus one is the subject of
\cite{S2}).

\vspace{.2in}

\paragraph{\textbf{10.3. Canonical basis of
$\widehat{\mathbf{U}}_v(\mathcal{L}\n)$.}} We assume again that
$X$ is of genus zero. By Corollary~\ref{C:1001} we may use $\Psi$, to transport the basis 
$\{\mathbf{b}_{\mathbb{P}}\; \mathbb{P} \in \mathcal{P}\}$
of $\widehat{\mathfrak{U}}_{\mathbb{A}}$ to a (topological) basis of $\widehat{\mathbf{U}}_v(\mathcal{L}\n)$, whose
elements we still denote by $\mathbf{b}_{\mathbb{P}}$.
In particular, for $\alpha \in K(X)^{tor}$ we denote by
$\mathbf{b}_\a \in \mathbf{B}$ the canonical basis element
corresponding to $\mathbf{1}_{\a}$;  
likewise, by Corrolary~\ref{Co:90} iii), $\mathbf{IC}(O^\mathcal{F}) \in \mathcal{P}$ for any vector bundle $\mathcal{F}$
and we simply denote by $\mathbf{b}_{\mathcal{F}}$ the corresponding canonical basis element of 
$\widehat{\mathbf{U}}_v(\mathcal{L}\n)$. The following Lemma is clear.

\vspace{.1in}

\begin{lem} i) For any $i$, let $\mathbf{B}_i \subset \mathbf{B}_{\mathbf{H}_{p_i}}$
denote the usual canonical basis of
$\mathbf{U}_v^+(\widehat{\mathfrak{sl}}_{p_i}) \subset
\mathbf{H}_{p_i}$. Then $\mathbf{B}_i \subset \mathbf{B}$ and
$\mathbf{B}_i$ corresponds to the set of simple perverse sheaves
in $\mathcal{P}$ with support in $\{\phi: \mathcal{E}_n^\a \tto
\mathcal{F}\;|\;supp(\mathcal{F}) =\{\lambda_i\}\}$ for some
$\a$.\\
ii) For any line bundle $\mathcal{L} \in Pic(X,G)$ the element
$$\mathbf{b}_{\mathcal{L}}=E_{\mathcal{L}} + \sum_{\a \in
K(X)^{tor}} v^{\langle [\mathcal{L}], \a\rangle}
E_{[\mathcal{L}]-\a} \mathbf{b}_{\a}$$ belongs to the canonical
basis $\mathbf{B}$.
\end{lem}

\vspace{.1in}

Since Verdier duality commutes with all the functors
$\mathrm{Ind}^{\a,\b}_{n,m}$, it descends to a semilinear
involution of $\widehat{\mathfrak{U}}_{\mathbb{A}}$, and thus
gives rise to a semilinear involution $x \mapsto \overline{x}$ of
$\widehat{\mathbf{U}}_v(\mathcal{L}\n)$. From the definition of
$\Psi$ one gets
$$\overline{E_{(i,j)}}=E_{(i,j)}\;\qquad \text{for}\;(i,j)
\in \aleph,$$
$$\overline{\xi_l}=\mathbf{b}_{l\delta}-\sum_{\underset{n
<l}{n+n_1+\cdots+n_N=l}}
\overline{\xi}_n\overline{\tau_1^{-1}(u_{n_1}^{(1)})} \cdots
\overline{\tau_1^{-1}(u_{n_N}^{(N)})},$$
\begin{align*}
\overline{E_{\star,t}}&=\sum_{n \geq 0} (-1)^n
\sum_{\underline{\a} \in (K(X)^{tor})^n}v^{-\langle [\mathcal{O}],
\sum \a_j\rangle} \mathbf{b}_{[\mathcal{O}(t)]-\sum \a_j}
\mathbf{b}_{\a_n} \cdots
\mathbf{b}_{\a_1} \\
&= \sum_{n \geq 0} (-1)^n \sum_{\underline{\a} \in
(K(X)^{tor})^{n+1}} v^{d'(\underline{\a})} E_{[\mathcal{O}(t)]
-\sum \a_j} \mathbf{b}_{\a_{n+1}} \cdots \mathbf{b}_{\a_1},
\end{align*}
where $d'(\underline{\a})=\langle [\mathcal{O}],
\a_{n+1}-\sum_{j=1}^n \a_j \rangle -\langle \sum_{j=1}^n \a_j,
\a_{n+1} \rangle$. More explicit formulas seem difficult to obtain
in the general case.

\vspace{.1in}

\begin{prop}\label{P:last!}
The canonical basis $\mathbf{B}$ satisfies the following
properties~:
\begin{enumerate}
\item[i)] $\mathbf{b}'\mathbf{b}'' \in
\widehat{\bigoplus}_{\mathbf{b} \in \mathbf{B}}
\mathbb{N}[v,v^{-1}] \mathbf{b}$ for any $\mathbf{b}',
\mathbf{b}'' \in \mathbf{B}$.
\item[ii)] $\overline{\mathbf{b}}=\mathbf{b}$ for any $\mathbf{b}
\in \mathbf{B}$,
\item[iii)] Let $\kappa$ be the automorphism of
$\widehat{\mathbf{U}}_v(\mathcal{L}\n)$ mapping $E_{\star,t}$ to
$E_{\star, t+1}$ and all other generators $H_{\star,t}$,
$E_{(i,j)}$ to themselves. Then $\mathbf{b} \in \mathbf{B}$ if and
only if $\kappa(\mathbf{b}) \in \mathbf{B}$.
\item[iv)] For any $(i,j) \in \aleph$ and $n>0$ there exists subsets $\mathbf{B}_{(i,j), \geq n} \subset \mathbf{B}$
such that $\widehat{\mathbf{U}}_v(\mathcal{L}\n)\mathbf{b}_{n\alpha_{(i,j)}}=
\widehat{\bigoplus}_{\mathbf{b} \in \mathbf{B}_{(i,j), \geq n}} \mathbb{A} \mathbf{b}$.
\end{enumerate}
\end{prop}
\noindent \textit{Proof.} Statements i) and iii) are obvious. If
$\mathbb{P} \in \mathcal{P}^{vb}$ then
$\mathbb{P}=\mathbf{IC}(O^\mathcal{F})$ and thus
$D(\mathbb{P})=\mathbb{P}$. The same is true if $\mathbb{P} \in
\mathcal{P}^\a_{|\a|}$ (see Sections~5 and 8 for the notations).
On the other hand, if $\mathbb{P} \in \mathcal{P}^{l\delta}_0$
then by Lemma~\ref{L:16} we have
$\mathbb{P}=\mathbf{IC}(U^{l\delta}_{(1), \ldots, (1)}, \sigma)$
for some irreducible representation of the symmetric group
$\mathfrak{S}_l$. This is also self dual since any representation
of $\mathfrak{S}_l$ is self dual. Now let $\mathbb{P} \in
\mathcal{P}^{\a+l\delta}_{|\a|}$ be an arbitrary torsion simple
perverse sheaf. By Lemma~\ref{L:18} there exists simple perverse
sheaves $\mathbb{R} \in \mathcal{P}_0^{l\delta}$ and $\mathbb{T}
\in \mathcal{P}^\a_{|\a|}$ such that
$\mathrm{Ind}^{l\delta,\a}(\mathbb{R} \boxtimes \mathbb{T})=
\mathbb{P} \oplus \mathbb{P}_1$ with $\mathbb{P}_1$ being a sum of
shifts of complexes in $\mathcal{P}^{\a+l\delta}_{>|\a|}$. Since
Verdier duality commutes with $\mathrm{Ind}$, we obtain
$D(\mathbb{P})\oplus D(\mathbb{P}_1)=\mathbb{P} \oplus
\mathbb{P}_1$ from which it follows that
$D(\mathbb{P})=\mathbb{P}$. The case of an arbitrary $\mathbb{P}
\in \mathcal{P}$ now follows in the same manner using
Corollary~\ref{Co:90} i). This proves ii). Finally, iv) is obtained by
combining \cite{L1}, Theorem~14.3.2 with Corrolary~\ref{Co:90} i) and the description of
$Q_l^{n(\alpha_{(i,j)})}$ given in Section~5.1.
\qed

\vspace{.2in}

\section{Examples.}

\vspace{.2in}

\paragraph{\textbf{11.1. Canonical basis for $\widehat{\mathfrak{sl}}_2$}}. Let us
assume that $X=\mathbb{P}^1$ and $G=Id$, so that
$\mathbf{U}_v(\mathcal{L}\n)$ is the ``Iwahori'' positive part of
the quantum affine algebra
$\mathbf{U}_v(\widehat{\mathfrak{sl}}_2)$. In that situation, the
elements of the canonical basis corresponding to torsion sheaves
are indexed by partitions; let us identify the subalgebra of
$\mathbf{U}_v(\mathcal{L}\n)$ generated by $\{\xi_l, \; l\geq 1\}$
to Macdonald's ring of symmetric functions $\Lambda$ via the
assignment $\theta:\;\xi_l \mapsto p_l \in \Lambda$, where as
usual $p_l$ denotes the power sum symmetric function. From
Section~8.1 we deduce that
$$\mathbf{b}_{\lambda}:=\Psi^{-1}(\mathbf{IC}(U^{|\lambda|}_{(1),
\ldots, (1)}, \sigma_{\lambda}))=\theta^{-1}(s_{\lambda}),$$ where
$\sigma_{\lambda}$ and $s_{\lambda} \in \Lambda$ are respectively
the irreducible $\mathfrak{S}_{|\lambda|}$-module and the Schur
function associated to $\lambda$. In particular
$\xi_l=\mathbf{b}_{(l)}$ belongs to $\mathbf{B}$.

\vspace{.1in}

In rank one we have
$$\mathbf{b}_{\mathcal{O}(t)}=E_{t} + \sum_{l >0} v^l
E_{t-l} \xi_l$$ 
(for simplicity we omit the symbol $\star$). As an example of canonical basis elements
in rank two we give
$$\mathbf{b}_{\mathcal{O}(t)\oplus
\mathcal{O}(t)}=E_{t}^{(2)} + \sum_{\underset{t_2>0}{2t_1+t_2=2t}}
v^{2t_2}E_{t_1}^{(2)}\xi_{t_2} +
\sum_{\underset{t_1<t_2\; ; \; t_3\geq 0}{t_1+t_2+t_3=2t}} v^{1+t_2-t_1+2t_3}
E_{t_1}E_{t_2}\xi_{t_3},$$
\begin{equation*}
\begin{split}
\mathbf{b}_{\mathcal{O}(t)\oplus
\mathcal{O}(t+1)}=&v^2E_{t}E_{t+1} + \sum_{\underset{t_2>0}{2t_1+t_2=2t+1}}
v^{2t_2}E_{t_1}^{(2)}\xi_{t_2} \\
&\qquad +
\sum_{\underset{t_1<t_2\; ; \; t_3\geq 0}{t_1+t_2+t_3=2t+1}} v^{1+t_2-t_1+2t_3}
E_{t_1}E_{t_2}\xi_{t_3}.
\end{split}
\end{equation*}
The above examples of
canonical basis elements correspond to perverse sheaves
$\mathbf{1}_\a$ for some $\a \in K(X)$. Similar formulas give the
canonical basis elements corresponding to $\mathcal{O}(t)^{\oplus
n} \oplus \mathcal{O}(t+1)^{\oplus m}$ for any $n,m \in \N$.

\vspace{.1in}

Finally, the bar involution $x \mapsto \overline{x}$ is also easy
to describe explicitly~: it is given by $\overline{\xi_l}=\xi_l$
and
$$\overline{E_{t}}=E_{t} + \sum_{n>0}
(-1)^{n+1}\sum_{l_1, \ldots, l_n>0} v^{l_1-l_2-\cdots-l_n}
E_{t-\sum l_i} \xi_{l_1} \cdots \xi_{l_n}.$$

\vspace{.1in}

\paragraph{\textbf{Remark.}} From (\ref{A:1001}) and the relations 
\begin{align}
\Delta_{[\mathcal{O}(-n)],n\delta}
(\mathbf{b}_{[\mathcal{O}]})=&v^n \mathbf{b}_{[\mathcal{O}(-n)]} \otimes \boldsymbol{\xi}_n,\\
\Delta_{n\delta,[\mathcal{O}(-n)]}
(\mathbf{b}_{[\mathcal{O}]})=&v^{-n} \boldsymbol{\xi}_n \otimes\mathbf{b}_{[\mathcal{O}(-n)]},\\
\Delta_{n\delta,m\delta}(\boldsymbol{\xi}_{n+m})=&\boldsymbol{\xi}_n \otimes \boldsymbol{\xi}_m \label{E:deltaxi}
\end{align}
we may compute the coproduct on $\widehat{\mathfrak{U}}_{\mathbb{A}} = \widehat{\mathbf{U}}_v(\mathcal{L}\mathfrak{sl}_2)$.
By definition, $\Psi(E_{0})=\sum_{n \geq 0} v^n \mathbf{b}_{[\mathcal{O}(-n)]}\chi_n$ where
$\chi_n=\sum_{r >0} \sum_{l_1+\cdots + l_r=n}(-1)^r \boldsymbol{\xi}_{l_1} \cdots \boldsymbol{\xi}_{l_r}$. From 
(\ref{E:deltaxi}) we have
$\Delta_{n\delta,m\delta}(\chi_{n+m})=\chi_n \otimes \chi_m$ hence
\begin{equation*}
\Delta_{[\mathcal{O}(-l)],l\delta}(E_{0})=\sum_{n_1 \geq 0}
v^{l+n_1}\mathbf{b}_{\mathcal{O}(-l-n_1)}\chi_{n_1} \otimes (\boldsymbol{\xi}_l + 
\boldsymbol{\xi}_{l-1}\chi_1 + \cdots + \chi_l).
\end{equation*}
One finds that $\sum_{i+j=l} \boldsymbol{\xi}_i \chi_j=\delta_{l,0}$, so that
$\Delta_{[\mathcal{O}]-l\delta,l\delta}=\delta_{l,0}E_{0} \otimes 1$.
A similar computation gives
$$\Delta_{l\delta, [\mathcal{O}(-l)]}(E_{0})=\sum_{n_1 \geq 0}(v^{n_1-l} \boldsymbol{\xi}_{l} + v^{2+n_1-l} 
\boldsymbol{\xi}_{l-1}\chi_1 + \cdots + v^{n_1+l} \chi_l) \otimes \mathbf{b}_{[\mathcal{O}(-l-n_1)]} \chi_{n_1}.$$
Hence, setting $\theta_l=\sum_{k=0}^l v^{2k-l} \boldsymbol{\xi}_{l-k}\chi_k$ we obtain
$$\Delta(E_{0})=E_{0} \otimes 1 + \sum_{l \geq 0} \theta_l \otimes E_{-l}.$$
Finally, one checks that the elements $\theta_l$ may be characterized through
the relation $\sum_l \theta_l s^l=
exp((v^{-1}-v)\sum_{r \geq 1} [r]H_{(r)})$, and we see that the above coproduct coincides with Drinfeld's new 
coproduct for quantum affine algebras (\cite{Dr}). This was already implicitly mentionned in \cite{Kap} but not written
in the litterature.

\vspace{.2in}

\paragraph{\textbf{11.2. Action on the principal subspace.}} Let us again consider the case of $\widehat{\mathfrak{sl}}_2$.
Let $L_{\Lambda_0}$ be the irreducible highest weight $\mathbf{U}_v(\widehat{\mathfrak{sl}}_2)$-module generated by
the highest weight vector $v_{\Lambda_0}$ of highest weight $\Lambda_0$, where $\langle h_0,\Lambda_0\rangle=0, \langle c,
\Lambda_0 \rangle=1$. Following \cite{BS}, we consider the subalgebra $\mathbf{U}_v(\widehat{\mathfrak{n}}) \subset
\mathbf{U}_v(\widehat{\mathfrak{sl}}_2)$ generated by $E_{t}$ for $t \in \Z$ and put 
$W_0=\mathbf{U}_v(\widehat{\mathfrak{n}}) \cdot v_{\Lambda_0}$. Observe that by definition $E_{t} v_{\Lambda_0}=
H_{(t)}v_{\Lambda_0}=0$ for any $t \geq 0$, and that therefore $\widehat{\mathbf{U}}_v(\mathcal{L}\n)$ acts on 
$L_{\Lambda_0}$. Moreover, we also have $W_0=\widehat{\mathbf{U}}_v(\mathcal{L}\n)\cdot
v_{\Lambda_0}$. Let $I_0$ be the annihilator in $\widehat{\mathbf{U}}_v(\mathcal{L}\n)$ of $v_{\Lambda_0}$.

\vspace{.1in}

\begin{prop} The ideal $I_0\otimes \C(v)$ of $\widehat{\mathbf{U}}_v(\mathcal{L}\n) \otimes \C(v)$ is generated by the 
elements $\mathbf{b}_{\alpha}$ for $\alpha \in K(X)^{tor}$, 
$\mathbf{b}_{\mathcal{O}(t)}$ for $t \geq 0$ and the elements $\mathbf{b}_{\mathcal{O}(t) \oplus \mathcal{O}(t)}$,
$\mathbf{b}_{\mathcal{O}(t) \oplus \mathcal{O}(t+1)}$ for $t <0$.
\end{prop}
\noindent
\textit{Proof (sketch).} This result is a quantum analogue of \cite{BS}, Theorem~2.2.1. Let $I'$ be the ideal
of $\widehat{\mathbf{U}}_v(\mathcal{L}\n)$ generated by the above elements. For $t<0$ we have
$$\mathbf{b}_{\mathcal{O}(t) \oplus \mathcal{O}(t)}\cdot v_{\Lambda_0}= 
E_{\star,t}^{(2)}\cdot v_{\Lambda_0} +
\sum_{-t>s>0} v^{1+2s}
E_{\star,t-s}E_{\star,t+s}\cdot v_{\Lambda_0},$$
$$\mathbf{b}_{\mathcal{O}(t)\oplus
\mathcal{O}(t+1)}\cdot v_{\Lambda_0}=v^2E_{\star,t}E_{\star,t+1}\cdot v_{\Lambda_0} +
\sum_{-t-1>s>0} v^{2+2s}
E_{\star,t-s}E_{\star,t+1+s}\cdot v_{\Lambda_0}.$$
It may be explicitly checked that these expressions vanish, using (for instance) the description of the action of 
$\mathbf{U}_v(\widehat{\mathfrak{sl}}_2)$ on $L_{\Lambda_0}$ in \cite{VV} together with \cite{Kevin}, Proposition~4.4.
This implies that $I_0$ contains $I'$. On the other hand, $W_0$ has finite rank weight spaces, so using \cite{BS}
Theorem~2.2.1 and specializing to
$v=1$ we have for any weight $\alpha$
\begin{equation}\label{E:torion}
\mathrm{dim}(\widehat{\mathbf{U}}_v(\mathcal{L}\n)/I_0[\alpha])_{|v=1}=\mathrm{dim}\;W_0[\alpha]_{|v=1}=\mathrm{dim}\;
(\widehat{\mathbf{U}}_v(\mathcal{L}\n)/I'[\alpha])_{|v=1}.
\end{equation}
We deduce that $I'\otimes \C(v)=I_0\otimes \C(v)$ as desired. \qed

\vspace{.1in}
Set $\mathcal{C}_0=\{\mathcal{F}\;| \mathcal{F}\simeq \mathcal{O}(l_1) \oplus \cdots \oplus \mathcal{O}(l_n)
,\;l_1 <l_2 \cdots < l_n <0,\; l_{i+1}-l_i \geq 2\}$.

\begin{conj}\label{Conj:1} 
We have $I_0=\widehat{\bigoplus}_{\mathbb{P} \not\in \mathcal{P}_{W_0}} \mathbb{A} \mathbf{b}_{\mathbb{P}}$,
where 
$$\mathcal{P}_{W_0}=\{\mathbf{IC}(\mathcal{F})\;|\;\mathcal{F} \in \mathcal{C}_0\}.$$
\end{conj}

\vspace{.1in}

The above conjecture and Proposition~\ref{P:last!} iv) may be (perhaps optimistically) generalized as follows~: 

\begin{conj} Let $(X,G)$ be an arbitrary datum of finite type (i.e $X=\mathbb{P}^1$), and let $\alpha \in K^+(X)$ be such
that $\mathbf{1}_{\a} \in \mathcal{P}$. Then there exists a subset $\mathcal{P}_{\geq \a} \subset \mathcal{P}$ such that
$$(\widehat{\mathfrak{U}}_{\mathbb{A}} \cdot \mathbf{b}_{\a})\otimes \C(v)=
\widehat{\bigoplus}_{\mathbb{P} \in \mathcal{P}_{\geq \a}} \C(v) 
\mathbf{b}_{\mathbb{P}}.$$ 

\end{conj}

The same conjecture may be made in the context of quivers for the canonical basis constructed in \cite{L1}. It is easy to
check that it fails in general over $\mathbb{A}$.

\vspace{.2in}

Assuming Conjecture~\ref{Conj:1}, we may equip $W_0$ with a ``canonical basis'' 
$B_0=\{\mathbf{b}_{\mathcal{F}} \cdot 
v_{\Lambda_0}|\mathcal{F} \in \mathcal{C}_{0}\}$. Recall that $L_{\Lambda_0}$ is already equipped with a 
canonical basis $B^+$, which is obtained by application
on $v_{\Lambda_0}$ of the canonical basis $\mathbf{B}'$ of $\mathbf{U}^-_v(\widehat{\mathfrak{sl}}_2)$ (as in \cite{L1}).
More precisely, let $\Lambda^\infty \supset L_{\Lambda_0}$ be the Level one Fock space and let 
$\{b^+_{\lambda}\;|\lambda \in \Pi\}$ be its Leclerc-Thibon canonical basis (see e.g \cite{VV}), which is indexed by the
set $\Pi$ of all partitions. It is known that $L_{\Lambda_0}=\bigoplus_{\lambda \in \Pi'}\mathbb{A} b^+_{\lambda}$ where
$\Pi'=\{\lambda \in \Pi\;|\; \lambda_i \neq \lambda_{i+1}\;\text{for\;all\;}i\}$. 

One calculates directly that $\mathbf{b}_{\mathcal{O}(-l)} \cdot v_{\Lambda_0}=b^+_{(2l+1)}$. 
More computational evidence suggest that in fact

\begin{conj}\label{Conj:2} We have $B_0 \subset B^+$.
In particular, the canonical basis $B^+$ of $L_{\Lambda_0}$ is compatible with the principal subspace $W_0$.
\end{conj}

Note that the last statement in the conjecture is not obvious from the definition of $B^+$ since the subalgebra
of $\mathbf{U}^-_v(\widehat{\mathfrak{sl}}_2)$ generated by $\{E_{1,n}\;| \; n <0\}$ is \textit{not} compatible with
the canonical basis $\mathbf{B}'$ of $\mathbf{U}^-_v(\widehat{\mathfrak{sl}}_2)$.

\vspace{.1in}

The whole space $L_{\Lambda_0}$ may be recovered from the principal subspace $W_0$ by a limiting process as follows. The 
affine Weyl group $W_a=\mathbb{Z}_2 \rtimes \mathbb{Z}$ contains a subalgebra of translations $T_n$, $n \in \Z$. The weight
spaces of $L_{\Lambda_0}$ of weight $T_n(\Lambda_0)$ are one-dimensional and there exists a unique collection of vectors
$v_n \in L_{\Lambda_0}[T_n(\Lambda_0)]$ such that $\mathbf{b}_{\mathcal{O}(-2n+1)}\cdot v_{n-1}= E_{-2n+1}\cdot v_{n-1}
=v_n$ (see \cite{BS}; explicitly, we have $v_n=b^+_{(2n+1, 2n, \ldots, 1)}$ if $n >0$ and $v_n=b^+_{(-2n, 1-2n, \ldots, 1)}$
if $n <0$). We set $W_n=\widehat{\mathbf{U}}_v(\mathcal{L}\n)\cdot v_n=\mathbf{U}_v(\widehat{\n})\cdot v_n$. The subspaces
$W_n$ for $n \in \Z$ form an exhaustive filtration $\cdots W_{-1} \supset W_0 \supset W_1 \cdots$ of $L_{\Lambda_0}$.
Finally, if $n >0$, we set 
$\mathcal{C}_n=\{\mathcal{F}(-2n)\;|\mathcal{F} \in \mathcal{C}_0\}$ and 
$\mathcal{P}_{W_n}=\{\mathbf{IC}(\mathcal{G} \oplus \mathcal{O}(-2n+1) \oplus \mathcal{O}(-2n+3) \cdots \oplus 
\mathcal{O}(-1))\;|\mathcal{G} \in \mathcal{C}_n\}$.

\begin{lem}\label{L:111} Fix $n >0$. Then
\begin{enumerate}
\item[i)] For any $\mathcal{G} \in \mathcal{C}_n$, we have 
\begin{equation*}
\begin{split}
\mathbf{b}_{\mathcal{G}} \cdot \mathbf{b}_{\mathcal{O}(-2n+1)}&
\cdot \mathbf{b}_{\mathcal{O}(-2n+3)} \cdots \mathbf{b}_{\mathcal{O}(-1)} \\
&\in \mathbf{b}_{\mathcal{G} \oplus 
\mathcal{O}(-2n+1) \oplus \mathcal{O}(-2n+3) \oplus \cdots \oplus \mathcal{O}(-1)} \oplus 
\widehat{\bigoplus}_{\mathbb{P} \not\in \mathcal{P}_{W_0}} \mathbb{A} \mathbf{b}_{\mathbb{P}}.
\end{split}
\end{equation*}
\item[ii)] We have 
$$\widehat{\mathbf{U}}_v(\mathcal{L}\n)\cdot \mathbf{b}_{\mathcal{O}(-2n+1)}
\cdot \mathbf{b}_{\mathcal{O}(-2n+3)} \cdots \mathbf{b}_{\mathcal{O}(-1)}=
\widehat{\bigoplus}_{\mathbb{P} \not\in \mathcal{P}_{W_0}}\mathbb{A}\mathbf{b}_{\mathbb{P}}
\oplus \widehat{\bigoplus}_{\mathbb{P} \in \mathcal{P}_{W_n}}
\mathbb{A}\mathbf{b}_{\mathbb{P}}.$$
\end{enumerate}
\end{lem}
\noindent
\textit{Proof.} This is shown using arguments similar to those of Section~9.1. We omit the details.\qed

\vspace{.1in}

Assuming Conjecture~\ref{Conj:1}, it follows from the above Lemma that the basis $B_0$ restricts to bases $B_n$ of
$W_n$ for any $n >0$, and that moreover the basis $B_n$  is the result of the application of the canonical basis
$\mathbf{B}$ of $\widehat{\mathbf{U}}_v(\mathcal{L}\n)$ on the vector $v_n$ (i.e, $B_n=\{\mathbf{b}_{\mathbb{P}}\cdot v_n\;|
\mathbf{b}_{\mathbb{P}} \in \mathbf{B}\}$). Applying the automorphism $\kappa^{-2n}$ to
Lemma~\ref{L:111} (see Proposition~\ref{P:last!} iii)), we see that the collection of bases 
$B_n:=\{\mathbf{b}_{\mathbb{P}}\cdot v_n\;|\mathbf{b}_{\mathbb{P}} \in \mathbf{B}\}$ for
$n \in \Z$ are all compatible, and thus give rise to a well-defined basis $B_{-\infty}$ of $L_{\Lambda_0}$.
The natural generalization of Conjecture~\ref{Conj:2} is

\begin{conj} The basis $B_{-\infty}$ coincides with the Leclerc-Thibon basis $B^+$.\end{conj}

\vspace{.2in}

\section{Appendix}

\vspace{.2in}

\paragraph{\textbf{A.1. Proof of Lemma~\ref{C:8.1}.}} For simplicity we write $p$ for $p_i$.
Recall that the canonical basis $\mathbf{B}_{\mathbf{H}_{p}}$ is
naturally indexed by the set $G_\gamma$-orbits in
$\mathcal{N}^{(p)}_\gamma$ for all $\gamma$. Such an orbit is in
turn determined by the collection of integers
$d_i^k(O)=\mathrm{dim\;Ker\;}x_{i-k+1}\cdots x_{i-1}x_i~: V_i \to
V_{i-k}$ for $i \in \Z/p\Z$ and $k \geq 1$, where $x=(x_i) \in O$
is any point. Put
$m_{i,l}(O)=d_i^l(O)-d_i^{l-1}(O)+d_{i+1}^{l}(O)-d_{i+1}^{l+1}(O)$
and $\mathbf{m}(O)=\sum_{i,l} m_{i,l} e_{i,l}$ where $e_{i,l}$ are
formal variables (the element $\mathbf{m}$ determines
$\{d_i^k\}$). Let us call \textit{aperiodic} an element
$\mathbf{m}$ (or the corresponding orbit $O_{\mathbf{m}}$) for
which the following statement is true: for every $t \geq 1$ there
exists $i \in \Z/p\Z$ such that $m_{i,t} =0$. By \cite{L4}, we
have
\begin{equation}\label{E:A2}
\mathbf{U}_v^+(\widehat{\mathfrak{sl}}_p) =\bigoplus_{O \in
\mathcal{O}^{ap}} \C[v,v^{-1}] \mathbf{b}_{O},
\end{equation}
where $\mathcal{O}^{ap}$ is the set of aperiodic orbits.
Similarly, we call $\mathbf{m}$ and $O_{\mathbf{m}}$
\textit{totally periodic} if $m_{i,l}=m_{j,l}$ for any $i,j,l$.
Now let $\mathbf{b}_{\mathbf{m}}$ be as in the claim. From
Proposition~\ref{P:HallS0} we deduce directly that
$\mathbf{b}_\mathbf{m} \in
\mathbf{U}^+_v(\widehat{\mathfrak{sl}}_p)$ if $\gamma \neq (r,r,
\ldots, r)$ for some $r \in \N$, and that there exists $c \in \C$
such that
\begin{equation}\label{E:A1}
\mathbf{b}_{\mathbf{m}}-cz_r \in
\mathbf{U}^+_v(\widehat{\mathfrak{sl}}_p)
\end{equation}
$\;\text{if\;} \gamma=(r,r,\ldots,r)$. Let $s: \mathbf{H}_p
\stackrel{\sim}{\to} \mathbf{H}_p$ be the automorphism of order
$p$ induced by the cyclic permutation of the quiver. It is clear
that $s(\mathbf{U}^+_v(\widehat{\mathfrak{sl}}_p))=
\mathbf{U}^+_v(\widehat{\mathfrak{sl}}_p)$ and we have
$s(z_r)=z_r$. Hence, if $\mathbf{b}_{\mathbf{m}}$ satisfies
(\ref{E:A1}) then either $c=0$ and $\mathbf{m}$ is aperiodic, or
$c \neq 0$ and $\mathbf{m}$ is fixed by $s$, hence totally
periodic. The set of integers $\{m_{i,l}\}$ for any \textit{fixed}
$i$ then forms a partition $\lambda$ independent of $i$. Let $l$
be the least nonzero part of $\lambda$, put
$\mathbf{m}'=\mathbf{m}-\sum_i e_{i,l}$ and
$\mathbf{d}'=(1,\ldots, 1), \mathbf{d}''=(r-1, \ldots, r-1)$.
Using \cite{S0}, Lemma~3.4 we have
$$\Delta_{\mathbf{d}'',\mathbf{d}'}(v^{\mathrm{dim}\;O_{\mathbf{m}}}1_{O_{\mathbf{m}}})_{|O_{\mathbf{m}'} \times \{0\}}
\in 1 + v\mathbb{N}[v],$$ and, for any $\mathbf{n}$
$$\Delta_{\mathbf{d}'',\mathbf{d}'}(v^{\mathrm{dim}\;O_{\mathbf{n}}}1_{O_{\mathbf{n}}})_{|O_{\mathbf{m}'} \times \{0\}}
\in \mathbb{N}[v].$$ But
$$\mathbf{b}_{\mathbf{m}}\in v^{\mathrm{dim}\;O_{\mathbf{m}}}1_{O_{\mathbf{m}}}\oplus
\bigoplus_{O_{\mathbf{n}} \subset \overline{O_{\mathbf{m}}}}
v^{\mathrm{dim}\;O_{\mathbf{n}}+1}\mathbb{N}[v]
1_{O_{\mathbf{n}}}.$$ We deduce
\begin{equation}\label{E:A3}
\Delta_{\mathbf{d}'',\mathbf{d}'}(\mathbf{b}_{\mathbf{m}})_{|O_{\mathbf{m}'}
\times \{0\}} \in 1 + v\mathbb{N}[v].
\end{equation}
By hypothesis
$\Delta_{\mathbf{d}'',\mathbf{d}}(\mathbf{b}_{\mathbf{m}})_{|O_{\mathbf{m}'}
\times \bullet} \in \C_{G}(\mathcal{N}_{\mathbf{d}'}^{(p)})$
belongs to $\mathbf{U}^+_v(\widehat{\mathfrak{sl}}_p)$ and is
$s$-invariant. But since the orbit under $s$ of any aperiodic
$\mathbf{n}$ is of size at least two, (\ref{E:A2}) together with
the fact that $\Delta(\mathbf{b}_{\mathbf{m}}) \in
\bigoplus_{\mathbf{n},\mathbf{n}'} \mathbb{N}[v,v^{-1}]
\mathbf{b}_{\mathbf{n}}\otimes \mathbf{b}_{\mathbf{n}'}$ implies
that all nonzero coefficients of
$\Delta_{\mathbf{d}'',\mathbf{d}'}(\mathbf{b}_{\mathbf{m}})_{|O_{\mathbf{m}'}
\times \{0\}}$ are greater than two. This contradicts
(\ref{E:A3}), and hence proves that $c=0$ in (\ref{E:A1}). The
claim follows. \qed

\vspace{.2in}

\paragraph{\textbf{A.2. Proof of Lemma~\ref{L:101}.}} We again drop the index $i$ throughout, and write
$1_l$ for $1_{\mathcal{N}^{(p)}_{(l,\ldots,l)}}$, and
$\Delta_{l,m}$ for $\Delta_{(l, \ldots, l),(m, \ldots, m)}$.
Finally, the computations below are conducted up to a global power
of $v$. We prove the Lemma by induction. The result is empty for
$r=0$. Assume that $u_1, \ldots, u_{r-1}$ satisfying
(\ref{E:equul}) are already defined and belong to
$\mathbf{U}^+_v(\widehat{\mathfrak{sl}}_p)$, and let us put
$u_r=1_r-\zeta_{r}-\zeta_{r-1}u_1 - \cdots -\zeta_1u_{r-1}$.
Substituting, we obtain
$$u_r=\sum_{(n_0, \ldots, n_r) \in \mathcal{J}} \left(\begin{matrix} n_1+ \cdots +n_r\\ n_1,\cdots,n_r\end{matrix}\right)
(-\zeta_1)^{n_1} \cdots (-\zeta_r)^{n_r} 1_{n_0},$$ where
$\mathcal{J}=\{(n_0, \ldots, n_r)\;|\;sum n_i=r\}$. A simple
calculation using the relations $\Delta_{l,m}(1_{l+m})=1_{l}
\otimes 1_m$ and $\Delta_{l,m}(\zeta_{l+m})= \zeta_l \otimes
\zeta_m$ yields $\Delta_{l,r-l}(u_r)=u_l \otimes u_{r-l} \in
\mathbf{U}_v^+(\widehat{\mathfrak{sl}}_p) \otimes
\mathbf{U}_v^+(\widehat{\mathfrak{sl}}_p)$ for all $l \geq 1$.
Proposition~\ref{P:HallS0} now implies that there exists $c \in
\C$ and $u' \in \mathbf{U}_v^+(\widehat{\mathfrak{sl}}_p)$ such
that $u_r=cz_r + u'$. It remains to prove that $c=0$. Put
$\mathbf{H}'_p=\sum_i E_i \mathbf{H}_p$. Since $z_r \not\in
\mathbf{H}'_p$ but $\zeta_l u_{r-l} \in \mathbf{H}'_p$ for all $l
\geq 1$, it is enough to show that
\begin{equation}\label{E:A4}
1_r-\zeta_r \in \mathbf{H}'_p.
\end{equation}
We argue once more by induction, this time on $p$. If $p=2$ then
$1_r-\zeta_r=1_A$ where $A=\{(x_1,x_2)\;|\; \mathrm{Ker}(x_2) \neq
\{0\}\} \subset \mathcal{N}_{(r, \ldots, r)}^{(2)}$. Let us set
$A_l=\{ (x_1, x_2)\;| \mathrm{dim\;Ker}\; x_2 =l\}$ so that
$A=\bigsqcup_{l \geq 1} A_l$. It is easy to see that, up to a
power of $v$, $1_{A_l}=1_{(l,0)} \cdot 1_{A'_l}$ where
$A'_l=\{(x_1,x_2)\;|\;x_2\;\text{is\; surjective}\} \subset
\mathcal{N}^{(2)}_{(r-l,r)}$. Since $1_{(l,0)}$ is proportional to
$E_1^l$, $E_1$ divides $1_{A_l}$, and hence also $1_r-\zeta_r$.
Next, assume that (\ref{E:A4}) is established for some $p_0 \geq
2$. Consider the algebra homomorphism $a: \mathbf{H}_{p_0} \to
\mathbf{H}_{p_0+1}$ induced by the collection of embeddings $a':
\mathcal{N}_{\underline{V}}^{(p_0)} \to
\mathcal{N}_{\underline{W}}^{(p_0+1)}$ defined by $W_1=W_2=V_1$,
$W_i=V_{i-1}$ for $2 \leq i \leq p_0+1$ and $a'(x_1, \ldots,
x_{p_0})=(x_1, Id, x_2, \ldots, x_{p_0})$. Note that
$a(\mathbf{H}_{p_0}') \subset \mathbf{H}'_{p_0+1}$ since
$a'(E_1)=v^{-1}E_2E_1-E_1E_2$ and $a'(E_i)=E_{i+1}$ for $i \neq
1$. Thus $1_{B}-\zeta_r= a(1_r-\zeta_r)\in \mathbf{H}'_{p_0+1}$,
where $B=\{(y_1, \ldots, y_{p_0+1})\;|\; \mathrm{Ker\;}y_2=\{0\}\}
\subset \mathcal{N}^{(p_0+1)}_{(r, \ldots, r)}$. It only remains
to check that $1_C=1_r-1_B$ belongs to $\mathbf{H}'_{p_0+1}$,
where $C=\{(y_1, \ldots, y_{p_0+1})\;|\;
\mathrm{Ker\;}y_2\neq\{0\}\}$. This is proved in the same way as
in the case $p=2$. Thus $1_r-\zeta_r \in \mathbf{H}'_{p_0+1}$ and
the induction is closed.\qed

\vspace{.2in}

\centerline{\textbf{Acknowledgements}}

\vspace{.1in}

I would like to thank Ga\"etan Chenevier and Eric Vasserot for
many very helpful discussions. In particular, the work \cite{BS} was pointed out to me
by Eric Vasserot. Parts of this work were done at
Yale University and Tsinghua University, and I warmly thank these two
institutions for their hospitality.

\small{}

\vspace{4mm}
\noindent
Olivier Schiffmann,\\
 D.M.A ENS Paris, 45 rue d'Ulm, 75230 Paris Cedex 05 FRANCE\\
email:\;\texttt{schiffma@dma.ens.fr}

\end{document}